\documentclass[12pt]{article}
\usepackage{amsmath,amsthm}
\usepackage{graphicx}
\usepackage{amssymb,latexsym}

\newtheorem{theorem}{Theorem}[section]
\newtheorem{lemma}[theorem]{Lemma}

\newtheorem{claim}[theorem]{Claim}

\newtheorem{corollary}[theorem]{Corollary}

\theoremstyle{remark}
\newtheorem{remark}[theorem]{Remark}

{\hspace*{\fill}$\rule{.3\baselineskip}{.35\baselineskip}$\end{trivlist}}

\newcommand{\bi}{\begin{itemize}}
\newcommand{\ei}{\end{itemize}}
\newcommand{\bt}{\begin{theorem}}
\newcommand{\et}{\end{theorem}}
\newcommand{\bp}{\begin{proof}}
\newcommand{\ep}{\end{proof}}
\newcommand{\be}{\begin{equation}}
\newcommand{\ee}{\end{equation}}
\newcommand{\ben}{\begin{enumerate}}
\newcommand{\een}{\end{enumerate}}

\renewcommand{\and}{\text{~~ and ~~}}
\renewcommand{\part}{\partial}

\newcommand{\sech}{{\rm sech~}}

\newcommand{\eps}{\varepsilon}
\newcommand{\mm}{{\mu\over 2}}

\renewcommand{\gg}{\gamma}
\newcommand{\ggh}{\hat\gamma}
\newcommand{\oalpha}{\overline\alpha}

\newcommand{\system}[2]{\left\{\begin{array}{#1} #2   \end{array}\right.}

\begin{document}

\title{Long-time limit studies of an obstruction in the $g$-function mechanism for semiclassical focusing NLS}

\author{Sergey Belov \footnote{ Teknikringen 64, Stockholm 11428, Sweden, e-mail: sbelov@gmail.com},
Stephanos Venakides \footnote{ Department of Mathematics, Duke University, Durham, NC 27708, e-mail: ven@math.duke.edu.
SV thanks NSF for supporting this work under grants NSF DMS-0707488 and NSF DMS-1211638.}}

\abstract{ We consider the long-time properties of the an
obstruction in the Riemann-Hilbert approach to one dimensional
focusing Nonlinear Schr\"odinger equation in the semiclassical
limit for a one parameter family of initial conditions. For
certain values of the parameter a large number of solitons in the
system interfere with the $g$-function mechanism in the steepest
descent to oscillatory Riemann-Hilbert problems. The obstruction
prevents the Riemann-Hilbert analysis in a region in $(x,t)$
plane. We obtain the long time asymptotics of the boundary of the
region (obstruction curve). As $t\to\infty$ the obstruction curve
has a vertical asymptotes $x=\pm \ln 2$. The asymptotic analysis
is supported with numerical results. }

\maketitle

\section{Introduction}

Consider the one dimensional focusing Nonlinear Schr\"odinger
equation (NLS) in the semiclassical limit $\eps\to 0$
\begin{equation}
i \eps q_t+ \eps^2q_{xx}+2\left|q\right|^2q=0,\quad t>0
\label{1.1:NLS}
\end{equation}
subject to a one parameter family of initial conditions
\begin{equation}
q(x,0,\eps)=A(x) e^{\frac{i\mu}{\eps}S(x)},
\label{1.1:IC_generic}
\end{equation}
where real valued $A(x)$ decays fast as $x\to\pm\infty$ and real valued
$S(x)$ converges to $S_\pm$ as $x\to\pm\infty$.

This is a well known type of problems of finding the leading
asymptotic behavior of the solution $q(x,t,\eps)$ in a singular
limit. In the case of the semiclassical focusing NLS (as opposed
to the defocusing semiclassical NLS), the problem is known for
modulational instability when a smooth initial profile breaks into
a seemingly disordered structure.

The first progress in analysis of semiclassical NLS was made by
Miller and Kamvissis \cite{MillerKamvissis_98} when in numerical
studies they observed some order. This has lead to a number of
results \cite{BMM, BTVZ_07, Lyng_07, TVZzero_04, TVZlong2_06}
based on the Riemann-Hilbert approach to this completely
integrable equation. The Riemann-Hilbert approach is based on
replacing the nonlinear PDE with a pair of linear operators (Lax
pair) first introduced by Lax for KdV equation \cite{Lax_68} and
later applied to NLS by Zakharov and Shabat \cite{ZSh_72}. This
reformulates the problem for a nonlinear PDE as a
scattering/inverse scattering problem for a linear operator. So
the asymptotic analysis of the NLS becomes an asymptotic analysis
of the spectral data of a linear operator where the initial data
for NLS plays the role of a potential. Then the problem is usually
further reformulated as a jump (factorization) problem on a
contour related to the spectrum of a linear operator - called
(oscillatory) Riemann-Hilbert problem (RHP).

Riemann-Hilbert problems are a natural object for the inverse
scattering as was noted by Shabat \cite{Shabat_76} who expressed
the hardest step - the inverse scattering as a multiplicative
matrix Riemann-Hilbert problem. A (local) RHP we define as the
following: find a matrix valued function $m(z)$, which is analytic
everywhere in the complex plane except on an oriented contour
$\Sigma$, where the function has a prescribed multiplicative
matrix jump. Additionally, the function must satisfy a
normalization condition at infinity. More precise description of
the RH approach can be found in \cite{Deift_OPbook_99, KMM_03}. A
simple example of a RHP is the jump matrix to be the identity
matrix ("free" case, no jump).

To extract the leading contribution, a $g$-function mechanism was introduced by
Deift, Venakides, Zhou \cite{DVZ1_97},
applicable to highly oscillatory RHPs. The method is
based on factoring out contributions until of the
remainder RHP has approximately constant (upto $L_2$ correction) jump matrix. Next the
"model" RHP with the constant jumps on finitely many intervals (finite genus)
is solved explicitly in terms of Riemann theta functions. Then one needs to show
that the remainder "error" RHP has a small ($L_2$) solution.
This method can be thought as a nonlinear steepest descent method.

The Riemann-Hilbert approach to asymptotic analysis has a wide range of
applications to a diverse array of problems including integrable systems
(sine-Gordon, Toda lattice, (m)KdV, (m)NLS, Benjamin-Ono), combinatorics
(longest increasing subsequences), Random matrices (GUE, GOE, beta ensembles),
and orthogonal polynomials (OPRL, discrete polynomials), to name some.

For NLS a number of initial conditions in the semiclassical limit
were analyzed \cite{BTVZ_07, Lyng_07, TVZzero_04, TVZ_07}. The
leading order solution of NLS was found in terms of Riemann theta
functions with underlying Riemann surfaces of finite genus. Other
existing results include long time analysis ($\eps=1$) along
straight lines $\frac{x}{t}=const$ in the $(x,t)$ plane for
several initial data \cite{BV_07, TVZlong2_06}. The analysis was
similar and lead to a finite genus ($0$, $1$, and $2$) regions in
the $(x,t)$ plane. Key ingredients in all these cases were Lax
pair operators, Riemann-Hilbert problems and the $g$-function
mechanism.

Consider the one parameter ($\mu\ge 0$) family of initial conditions
\begin{equation}
q(x,0,\eps)=-\sech(x)\ e^{\frac{i\mu}{\eps}\int_0^x \tanh s ds}.
\label{1.1:IC_original}
\end{equation}

Even the simplest case $\mu=0$ carries many features and difficulties in the
analysis. A common approach is to approximate the initial data without
disturbing the leading order of the solution. Moreover, choosing a special
sequence $\eps_n\to 0$ leads to purely multi-soliton solution which
is much simpler for numerical studies. This case was analyzed by Lyng, Miller
\cite{Lyng_07}.

For $\mu>0$, the family of initial conditions (\ref{1.1:IC_original}) combines
both solitonless initial data for $\mu\ge 2$ as well as radiation in the presence
of solitons for $0<\mu<2$. This makes it interesting from the point of view
of influence of a large number of solitons.

The solitonless case $\mu\ge 2$ has been analyzed completely for
all $x$ and $t>0$ values by Tovbis, Venakides, Zhou in
\cite{TVZzero_04}. In the semiclassical limit, the leading order
asymptotics is written in terms of Riemann theta functions with
parameters $\left\{\alpha_j\right\}$ which arise as branchpoints
and leading to a Riemann surface. They proved that there is a
curve $t=t_0(x)$ (called the first break) in the $(x,t)$ plane
such that for $0<t<t_0(x)$ the leading order of the solution
$q(x,t)$ depends only on $\alpha_0(x,t)$ and
$\alpha_1(x,t)=\overline{\alpha}_0(x,t)$. This can be seen as a
genus $0$ Riemann surface and the asymptotic solution has a WKB
type approximation. For $t>t_0(x)$ the leading order depends on
$\alpha_0(x,t)$, $\alpha_2(x,t)$, $\alpha_4(x,t)$ and their
complex conjugates (the genus is 2).  In the case of radiation
with solitons ($0<\mu<2$) the previous results were partial: only
for finite interval of $t$ values (not global in time)
\cite{TVparam_09}. In these studies, there was no information on
the region/boundary of rigorous applicability.

In this paper we consider the case $0<\mu<2$ with a similar semiclassical
approximation of (\ref{1.1:IC_original}), as it was done in \cite{TVZzero_04}. In
these studies the RH approach was completed for finite values of $t$ and was
not extendable globally for all
$t>0$. The main obstacle for $0<\mu<2$ came from a large number of solitons
(order of $O\left( \frac{1}{\eps}\right)$).
These solitons correspond to isolated poles of a reflection
coefficient of the underlying Lax operator. In the semiclassical
limit the isolated singularities accumulate and densely fill an
interval in the complex (spectral parameter or energy) plane. This
adds significantly to the difficulty of the asymptotic analysis,
which breaks as a leading contributing contour coming from
analyzing the oscillatory terms in the RHP, collides with these
accumulated poles and the error estimates become invalid.

This paper studies the boundary of the region of rigorous applicability of
the available asymptotic result. The boundary $t=t_s(x)$ (we call
it a singular obstruction curve) is a curve in $(x,t)$ plane. We
prove that the singular obstruction curve has a vertical
asymptotes $x=\pm \ln 2$. We find that the rate at which the
singular obstruction curve approaches these asymptotes is $x-\pm
\ln 2 = O\left(\frac{\ln t}{t}\right)$ as $t\to +\infty$. We also
provide the long-time asymptotics of all important quantities
including $\alpha_0$, $\alpha_2$, $\alpha_4$. It is conjectured
that for $|x|>\ln 2$ the solution maintains genus 2 asymptotics
for all $t>t_0(x)$ beyond the first break.

The paper is structured as follows: in section 2 we introduce the
main object of study - a scalar RHP on $g$-function. In section 3
we obtain the long time limit of the singular obstruction curve.
Section 4 provides numerical evidence supporting the asymptotic
analysis. In section 5 we discuss the results. Appendix is used
for all technical asymptotic computations.

\section{g-function problem}\label{RHP}

\begin{figure}
\begin{center}
\includegraphics[height=5cm]{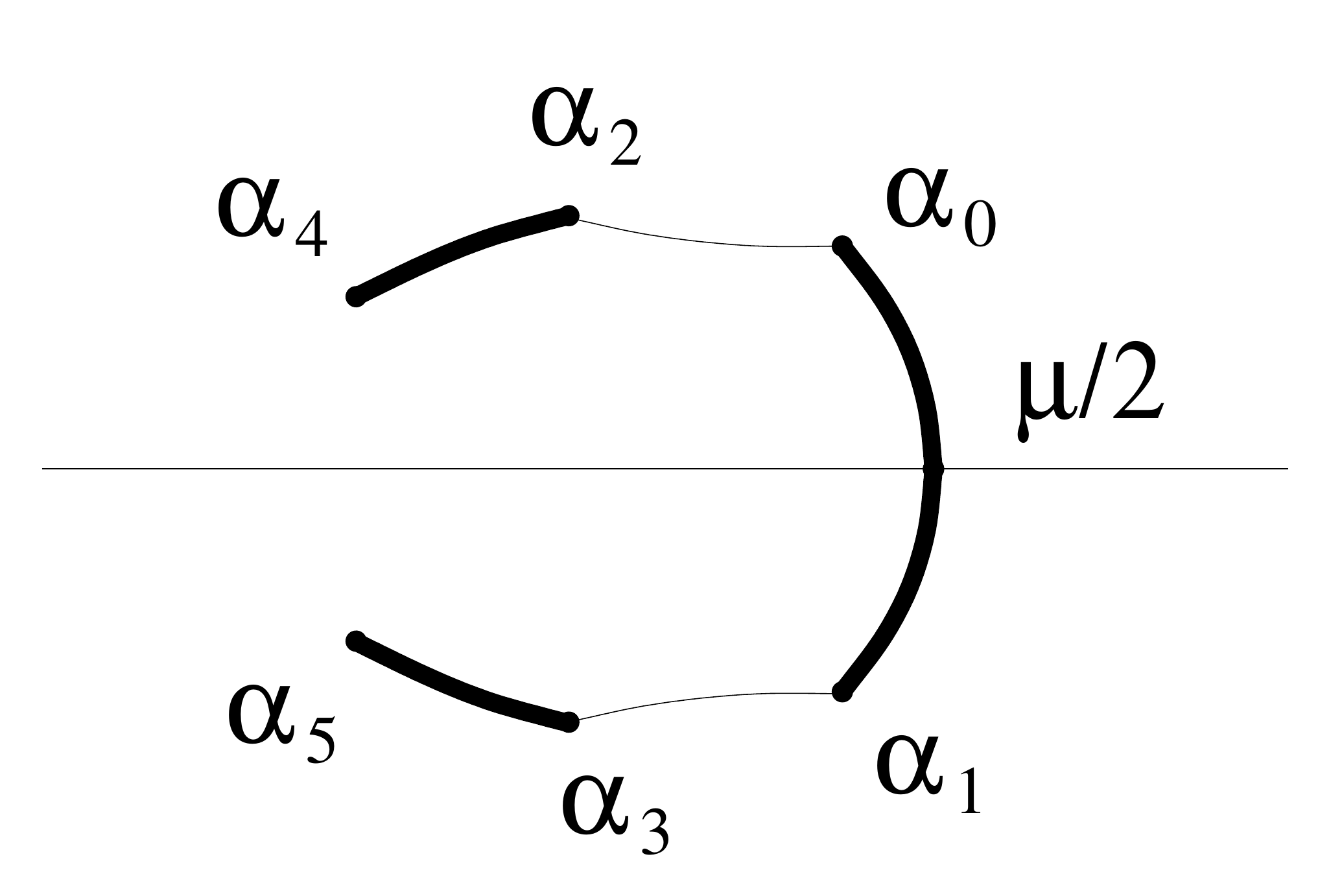}
\caption{\label{fig:main_compl_arcs} Main (bold) and complementary
(thin) arcs in the $g$-function problem. The exact positions of
the complementary arcs are not fixed.}
\end{center}
\end{figure}

The $g$-function mechanism was introduced in \cite{DVZ1_97} as a
method of extracting the leading order by factoring out an unknown
function $g(z)$ and setting up conditions to guarantee that this function
gives the leading order. Usually $g(z)$ is defined through
conditions on contours in the complex plane. For the semiclassical
NLS equation the contour is assumed
to be a union of so called "main" and "complementary" arcs on
which $g(z)$ is defined as the following in the case of genus 2 which is needed
for purposes of this paper. General setup of any finite genus is similar \cite{TVparam_09}.
\begin{enumerate}
\item Main arcs ($\gamma_0$, $\gamma_{m}$):
\begin{equation}
\left\{\begin{array}{ll} g_+ + g_--f=W & \mbox{on $\gamma_{m}$}\\
g_+ + g_--f=0 & \mbox{on $\gamma_0$}\\
\Im (2g_--f)< 0\quad& \mbox{right from the main arcs $\gamma_0$, $\gamma_{m}$}\\
\Im (2g_+-f)< 0\quad & \mbox{left from the main arcs $\gamma_0$, $\gamma_{m}$}
\end{array}\right. \label{1.3:main_arcs_ineq_g}
\end{equation}
\item Complementary arcs ($\gamma_{c}$):
\begin{equation}
\left\{\begin{array}{lll}
g_+ - g_-=\Omega & \mbox{on $\gamma_{c}$}\\
\Im (2g-f)> 0\quad& \mbox{on at least one side from $\gamma_{c}$}\\
\end{array}\right.\label{1.3:compl_arcs_ineq_g}
\end{equation}
\end{enumerate}
where $\Omega$, $W \in \mathbb{R}$. The main arcs $\gamma_{0}$, $\gamma_{m}$ form
a branch cut structure which defines $g(z)$ \cite{TVZzero_04}.

In general, the number of main and complementary arcs is
determined for each pair $(x,t)$, where $x$ and $t$ enter in RHP
for $g$-function as parameters through $f=f(z,x,t)$. The function
$f$ is assumed to be known and it comes from the initial condition
(\ref{1.1:IC_original}) through the logarithm of the reflection
coefficient $f=\frac{2i}{\eps}\log (r)$.

In this paper we use $f$ obtained by the semiclassical
approximation of the initial condition (\ref{1.1:IC_original}) as
in \cite{TVZzero_04}
\[
f(z,x,t)=\left(\frac{\mu}{2}-z\right)\left[\frac{\pi
i}{2}+\ln\left(\frac{\mu}{2}-z\right)\right]+\frac{z+T}{2}\ln\left(z+T\right)
+\frac{z-T}{2}\ln\left(z-T\right)
\]
\begin{equation}
-T\tanh^{-1}\frac{T}{\mu/2}-xz-2tz^2 +\frac{\mu}{2}\ln 2, \quad
\mbox{when}\ \ \Im z > 0, \label{eq:general f-function}
\end{equation}
\[
f(z,x,t)=\overline{f(\overline{z},x,t)} \quad \mbox{when}\ \ \Im z < 0,
\]
where the branch cuts in the logarithms are chosen as the following: from $\mm$ along the real
axis to $+\infty$, from $T$ to $0$ and along the real axis to $+\infty$, from
$-T$ to $0$ and along the real axis to $-\infty$.
\begin{equation}
f'(z,x,t)=-\frac{\pi i}{2}-\ln\left(\frac{\mu}{2}
-z\right)+\frac{1}{2}\ln\left(z^2-T^2\right) -x-4tz, \quad
\mbox{when}\ \ \Im z > 0. \label{eq2.1:general f'-function}
\end{equation}
In the limit to the real axis from the upper half plane $f'$ is
\begin{equation}
\Im f'(z+i0)=\system{ll}{\frac{\pi }{2}, & z<0 \\ -\frac{\pi}{2}, & 0<z<\mm \\ \frac{\pi}{2}, & z>\mm } \quad z\in\mathbb{R}
\end{equation}
so $f'$ has a jump on the real axis from Schwarz symmetry.

\begin{figure}
\begin{center}
\includegraphics[height=4cm]{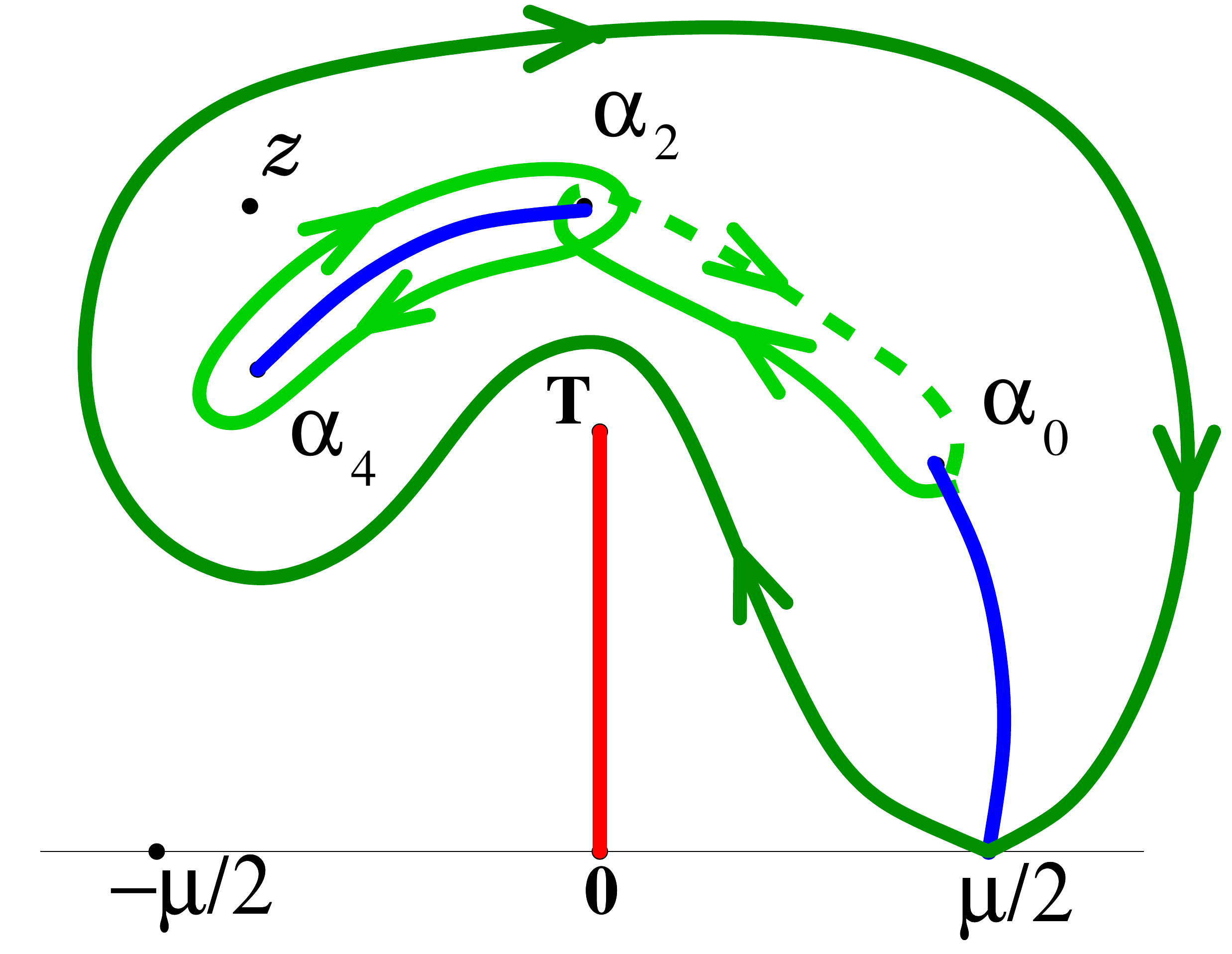}
\caption{\label{fig:h_fcn_cont} Contours of integration in the upper half plane for
$g(z)$, $h(z)$ and $B(z)$ in genus 2. In the lower half plane the contours are
symmetric. Large loop contour $\ggh$ and the small loop
contours $\ggh_m$, $\ggh_c$ are shown. For $h(z)$ the point $z$ is inside
of $\ggh$ and outside of $\ggh_m$, $\ggh_c$. For $g(z)$ the point $z$ is outside
of $\ggh$ and outside of $\ggh_m$, $\ggh_c$. For $B(z)$ the point $z$ is inside
of $\ggh$ and inside of $\ggh_m$, $\ggh_c$. As shown, the contours and $z$ suit to
compute $h(z)$.}
\end{center}
\end{figure}

The main and complementary arcs are described by the end points
$\{\alpha_j\}_{j=0}^5$:
$\gamma_0=\left[\alpha_1,\mm\right]\cup\left[\mm,\alpha_0\right]$,
$\gamma_m=\left[\alpha_2,\alpha_4\right]\cup\left[\alpha_5,\alpha_3\right]$,
and $\gg_c=[\alpha_0,\alpha_2]\cup[\alpha_3,\alpha_1]$. Because of
the Schwartz symmetry of the problem, $\alpha_1=\oalpha_0$,
$\alpha_3=\oalpha_2$, and $\alpha_5=\oalpha_4$. Introduce closed
loops $\ggh_0$, $\ggh_m$, and $\ggh_c$ around $\gamma_0$,
$\gamma_m$, and $\gamma_c$ respectively. The orientation of these
loop contours is clockwise (Fig. \ref{fig:h_fcn_cont}). The loop
$\ggh_0$ cannot be deformed away from $\mm$ since $f(z)$ is not
analytic at $z=\mm$. We also introduce a clockwise oriented closed
loop $\ggh$ enclosing all the main and complementary arcs
together. It is also passing through $z=\mm$.

The $g$-function have the following expression \cite{TVZzero_04} (Eq.(3.17), (3.18))
\begin{equation}
g(z)=\frac{R(z)}{2\pi
i}\left[\oint_{\ggh}\frac{f(\xi)}{(\xi-z)R(\xi)}d\xi +
\oint_{\ggh_c}\frac{\Omega}{(\xi-z)R(\xi)}d\xi+
\oint_{\ggh_m}\frac{W}{(\xi-z)R(\xi)}d\xi \right], \label{2.1:g_function}
\end{equation}
where point $z$ lies outside of a large loop $\ggh$ and outside
of small loops $\ggh_0$, $\ggh_c$ and $\ggh_m$.
The factor $R$ is
\[
R(\xi)=\sqrt{\left(\xi-\alpha_0\right)\left(\xi-\overline{\alpha}_0\right)
\left(\xi-\alpha_2\right)\left(\xi-\overline{\alpha}_2\right)
\left(\xi-\alpha_4\right)\left(\xi-\overline{\alpha}_4\right)},
\]
where $R(\xi)\sim -\xi^3$ as $\xi\to +\infty$ and it has branch cuts along the main arcs
$\gamma_0$, $\gamma_m$. The genus of
the Riemann surface of $R(\xi)$ is 2 since $\gamma_0$ and $\gamma_m$ form
three branch cuts.

The constants $W$ and $\Omega$ are solutions to the
system:
\begin{equation}
\left(
\begin{array}{cc}
\oint_{\ggh_m}\frac{1}{R(\xi)}d\xi &
\oint_{\ggh_c}\frac{1}{R(\xi)}d\xi\\
\oint_{\ggh_m}\frac{\xi}{R(\xi)}d\xi&
\oint_{\ggh_c}\frac{\xi}{R(\xi)}d\xi
\end{array}
\right)\left(
\begin{array}{c}
W\\
\Omega
\end{array}
\right)=\left(
\begin{array}{c}
-\oint_{\ggh}\frac{f(\xi)}{R(\xi)}d\xi\\
-\oint_{\ggh}\frac{\xi f(\xi)}{R(\xi)}d\xi
\end{array}
\right),
\end{equation}
which comes from the requirement $g(z)$ to be analytic at
infinity \cite{TVZzero_04}.

The branch points $\left\{\alpha_{0}, \alpha_{2},
\alpha_{4}\right\}$ are computed by solving the system
\begin{equation}
\left\{
\begin{array}{l}
B(\alpha_0)=0 \\
B(\alpha_2)=0 \\
B(\alpha_{4})=0
\end{array}
\right. , \label{2.1:B_system}
\end{equation}
where
\begin{equation}
B(z)=\oint_{\ggh}\frac{f(\xi)}{(\xi-z)R(\xi)}d\xi +
\oint_{\ggh_c}\frac{\Omega}{(\xi-z)R(\xi)}d\xi+
\oint_{\ggh_m}\frac{W}{(\xi-z)R(\xi)}d\xi,
\label{eq1.3:general_B_function}
\end{equation}
with $z$ being inside of the contours $\ggh$, $\ggh_{m}$ and $\ggh_{c}$
(see Fig. \ref{fig:h_fcn_cont}).

Additionally the 3 complex branch points satisfy a set of 4 real moment conditions \cite{TVZzero_04}:
\begin{equation}
\oint_{\ggh}\frac{\xi^j f'(\xi)}{R(\xi)}d\xi = 0, \quad j=0,1,2,3. \label{2.1:moment_conditions_g}
\end{equation}
These conditions come from expanding $g'(z)$ at infinity into a
power series. They are necessary but not sufficient to define
$\left\{\alpha_{0}, \alpha_{2}, \alpha_{4}\right\}$. The system
(\ref{2.1:B_system}) and formula (\ref{2.1:g_function}) imply
(\ref{2.1:moment_conditions_g}).

Introduce a convenient notation
\begin{equation}
h(z)=2g(z)-f(z)
\end{equation}
then
\[
h(z)=\frac{R(z)}{2\pi
i}\left[\oint_{\ggh}\frac{f(\xi)}{(\xi-z)R(\xi)}d\xi +
\oint_{\ggh_c}\frac{\Omega}{(\xi-z)R(\xi)}d\xi+
\oint_{\ggh_m}\frac{W}{(\xi-z)R(\xi)}d\xi \right],
\]
where point $z$ lies inside of a large loop $\ggh$ and outside
of small loops $\ggh_c$ and $\ggh_m$.

Then the $g$-function conditions
(\ref{1.3:main_arcs_ineq_g})-(\ref{1.3:compl_arcs_ineq_g}) are
rewritten in terms of $h(z)$ in genus 2:

\begin{enumerate}
\item Main arcs ($\gamma_0$, $\gamma_m$):
\begin{equation}
\left\{\begin{array}{ll} h_+ + h_-=0 & \mbox{on $\gamma_0$}\\
h_+ + h_-=2W & \mbox{on $\gamma_m$}\\
\Im h < 0\quad& \mbox{right from $\gamma_0$, $\gamma_m$}\\
\Im h < 0\quad & \mbox{left from $\gamma_0$, $\gamma_m$}
\end{array}\right. \label{1.3:main_arcs_ineq_h2}
\end{equation}
\item Complementary arcs ($\gamma_{c}$):
\begin{equation}
\left\{\begin{array}{lll}
h_+ - h_-=2\Omega & \mbox{on $\gamma_{c}$}\\
\Im h > 0 \quad & \mbox{left or right from $\gamma_{c}$}
\end{array}\right.\label{1.3:compl_arcs_ineq_h2}
\end{equation}
\end{enumerate}

This suggests the visualization: land ($\Im h>0$), sea ($\Im
h<0$), and sea-shore-lines or bridges ($\Im h=0$). A bridge has
$\Im h<0$ on both sides while a sea-shore-line has $\Im h>0$ on
one side $\Im h<0$ on the other. A main arc can be viewed as a
bridge connecting two land regions with sea on both sides. A
complementary arc is a land path with exact position being
unimportant as long as $\Im h\ge0$.

\begin{figure}
\begin{center}
\includegraphics[height=3.7cm]{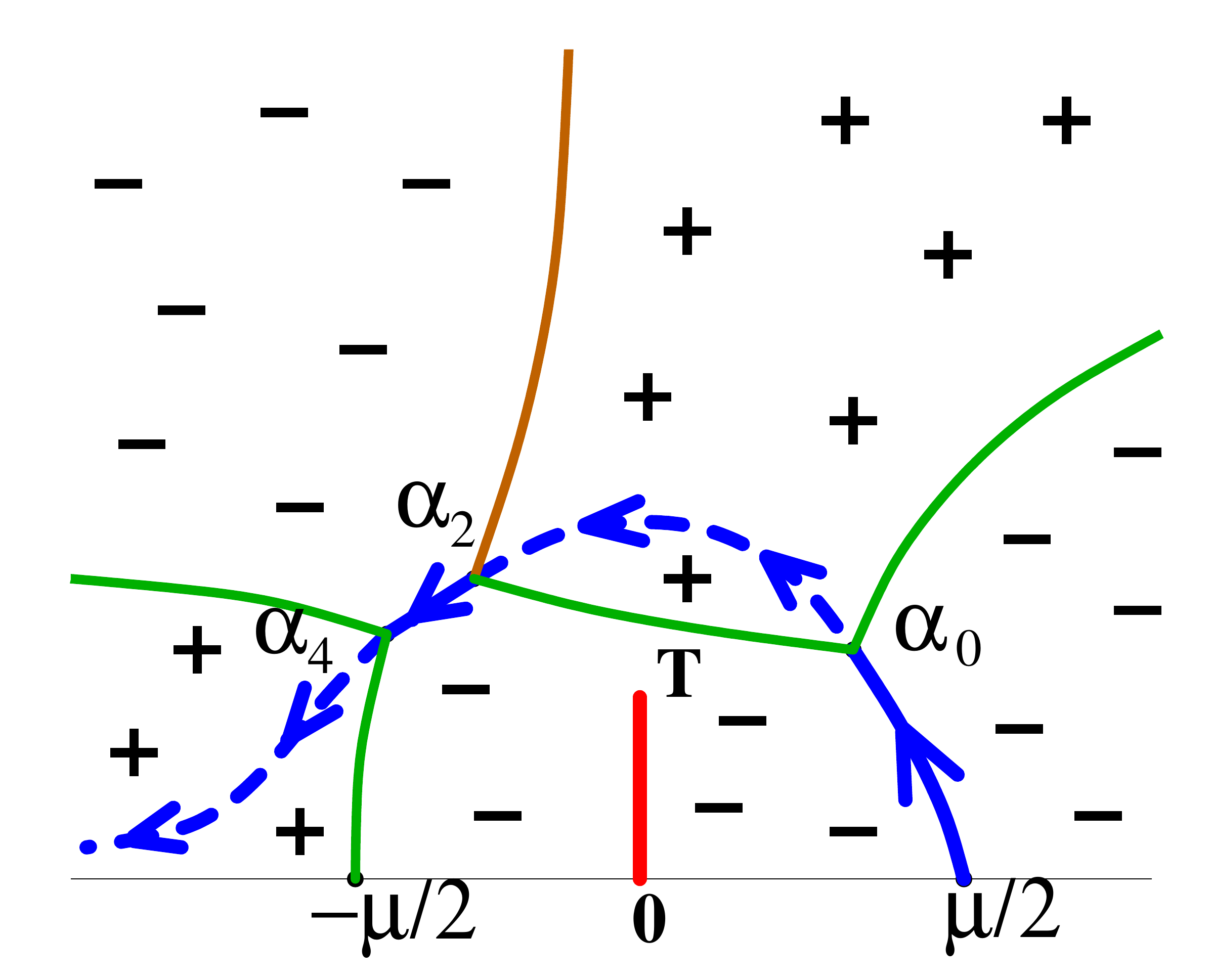}
\includegraphics[height=3.7cm]{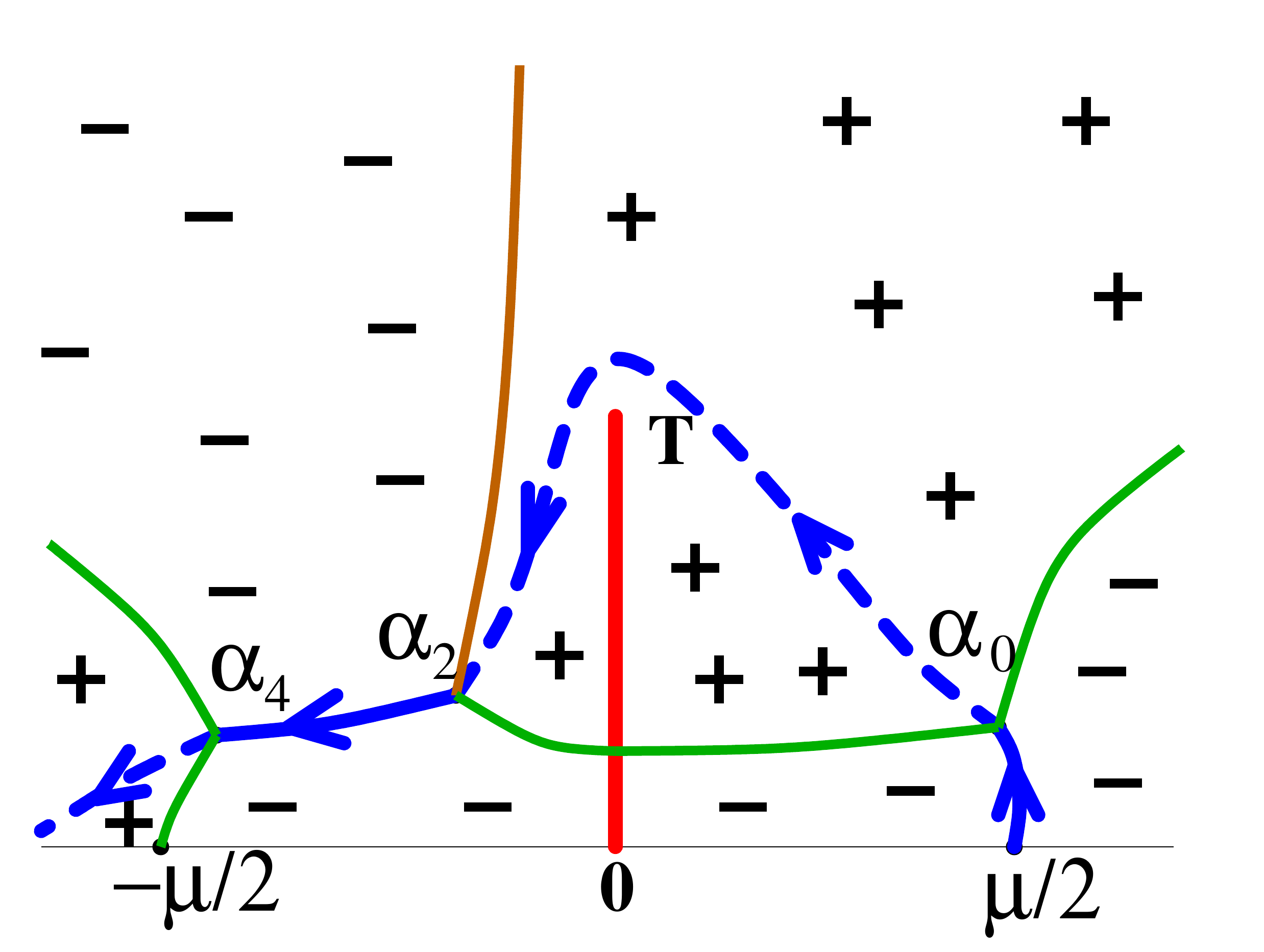}
\includegraphics[height=3.7cm]{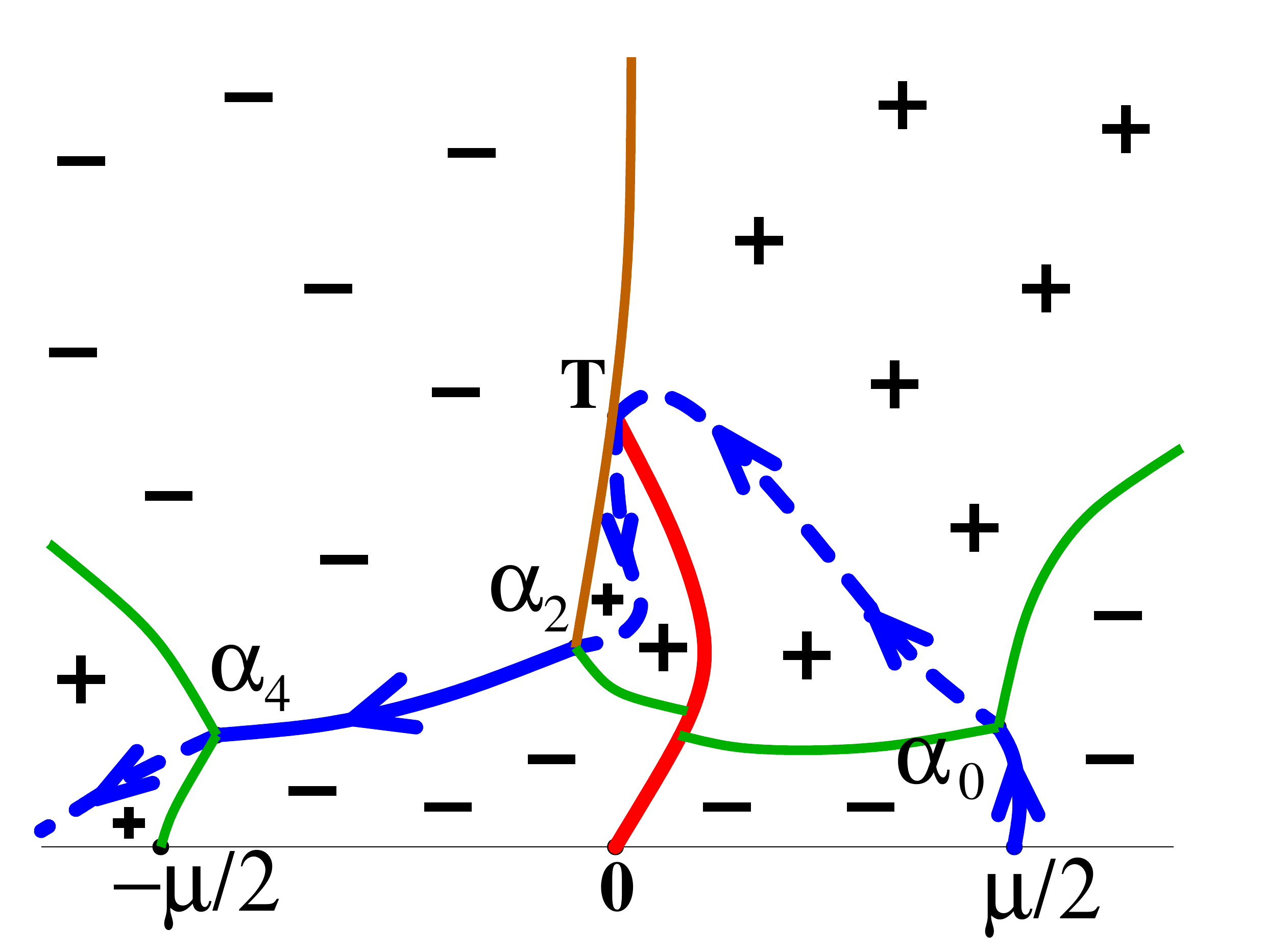}
\end{center}
\caption{ \label{fig:sing_br} Schematic illustration of the singular obstruction mechanism
in the solitons+radiation case ($0<\mu<2$) in the upper
half plane. The zero level curve from $\alpha_2$ going up
touches the top of the branch cut $[0, T]$ and closes the passage
to connect $\alpha_0$ and $\alpha_2$ (dashed).}
\end{figure}

The so called singular obstruction in the procedure occurs when the above
assumptions that there is contour connecting $\mm$ and $-\infty$
consisting of
main and complementary arcs exists is invalid. In particular as we show numerical
results in Section 4, this scenario occurs at finite values of $t$ for small values of $x$ ($|x|<\ln 2$). One of
the complementary arcs collides with the logarithmic branch cut
$[0,T]$. More precise assumption on the contours, that there is a
path connecting $\mm$ and $-\infty$ along which $\Im h(z)\ge 0$.
This condition could fail as shown on Figure \ref{fig:sing_br} for small $x$
values as time increases. We call such curve $x=x_c(t)$ in the $(x,t)$ plane
the singular obstruction curve defined by the condition
\begin{equation}
\Im h(T,x_c(t),t)=0, \label{eq:im_h_=_0}
\end{equation}
where $h(z=T,x,t)$ is understood in the sense of limit:
\begin{equation}
\Im h(T,x,t)=\Im\left( 2g(T,x,t)-\lim_{z\to T} f(z,x,t)\right).
\end{equation}
Equation (\ref{eq:im_h_=_0}) is an implicit condition for the singular obstruction curve $x=x_c(t)$
which we solve asymptotically below. Numerical investigation
suggests that this curve has a vertical asymptote in the $(x,t)$ plane.
In the next section we perform asymptotic analysis of the
long-time asymptotics of the singular obstruction curve.

\section{Long time asymptotic analysis}\label{asymp}

This section contains the core analysis of the paper. It is devoted to solving of
the equation (\ref{eq:im_h_=_0}) in the long time limit.

\begin{claim}\label{claim_1}
In genus 2 in the long time limit the branch points $\alpha_0$, $\alpha_2$
and $\alpha_4$ converge to $\mu/2$, $0$, $-\mu/2$ respectively. Convergence
of $\alpha_0\to \mu/2$ and $\alpha_4\to -\mu/2$ is exponentially fast in $t$ as $t\to\infty$.
\end{claim}

For $\mu\ge 2$, the result is proved in \cite{TVZlong2_06}. In the case
$0<\mu<2$, an additional logarithmic branch cut $[0,T]$ in the upper half plane
appears in equations from which $\alpha_0$, $\alpha_2$, $\alpha_4$ are determined.
The additional branch cut can be viewed as a small perturbation in the
limit $t\to\infty$ and does not affect the
leading behavior of the branch points in the claim.

\begin{corollary}\label{cor_1}
Claim \ref{claim_1} implies that the integrals
\begin{equation}
\int_{\mm}^{\alpha_0}\frac{f'(\xi)}{(\xi-z)R(\xi)}d\xi,\quad
\int_{\mm}^{\overline\alpha_0}\frac{f'(\xi)}{(\xi-z)R(\xi)}d\xi,
\end{equation}
\begin{equation}
\int_{\alpha_4}^{-\mm}\frac{f'(\xi)}{(\xi-z)R(\xi)}d\xi,\quad
\int_{\overline\alpha_4}^{-\mm}\frac{f'(\xi)}{(\xi-z)R(\xi)}d\xi
\end{equation}
are exponentially small in time as $t\to\infty$.
\end{corollary}

Proof of Corollary \ref{cor_1} is exactly the same as in \cite{TVZlong2_06}.

\begin{theorem}\label{thm_1}
Let $0<\mu<2$. Assuming the genus 2 region in the $(x,t)$ plane allows to send $t$ to infinity,
the long time larger solution of the equation
\[
\left\{
\begin{array}{l}
\Im h(T,x,t)=0,\\
x>0
\end{array}
\right.
\]
has the following long time asymptotics
\begin{equation}
x(t)=\ln 2-\frac{\left(8ta_2+\ln b_2\right)b_2^2}{2|T|^2}-C_3b_2^2+
O\left(a_2b_2\right), \ \ \ \ t\to\infty, \label{thmbr2a2b2}
\end{equation}
where $\alpha_2=a_2+ib_2$, $T=i\sqrt{1-\frac{\mu^2}{4}}$, $C_3=\frac{1-2\ln|T|}{4|T|^2}$ and
\begin{equation}
\left\{
\begin{array}{l}
a_2(x,t)=\left.\frac{A_1\ln t}{t}+\frac{A_2}{t}+\frac{A_3\ln t}{t^{2}}+O\left(\frac{1}{t^{2}}\right)\right. \\
b_2(x,t)=\frac{B_1}{\sqrt{t}}+\frac{B_2\ln t}{t^{3/2}}
 +\frac{B_3}{t^{3/2}}+O\left(\frac{\ln t}{t^{5/2}}\right)
\end{array},\quad t\to\infty
\right.\label{thma2b2}
\end{equation}
with
\begin{equation}
\left\{
\begin{array}{ll}
A_1=\frac{1}{8},         &   B_1=\frac{\sqrt\mu}{2}\\
A_2=\frac{1}{4}\ln{\frac{2|T|}{\sqrt\mu}}+\frac{\ln
2-x}{4},                  &   B_2=\frac{1}{16\sqrt{\mu}}\\
A_3=\frac{3}{32\mu},     &   B_3= \frac{2+\ln
\frac{2|T|}{\sqrt\mu} }{8\sqrt{\mu}}+\frac{\ln 2-x}{8\sqrt{\mu}}\\
\end{array}
\right.\label{thma2b2c}
\end{equation}
or explicitly in terms of $t$
\begin{equation}
x(t)=\ln 2+c_2\frac{\ln t}{t}+c_3\frac{1}{t}+ O\left(\frac{\ln
t}{t^{3/2}}\right), \ \ \ \ t\to\infty, \label{thmbr2}
\end{equation}
where
\begin{equation}
c_2=-\frac{\left(8A_1-\frac{1}{2}\right)B_1^2}{2|T|^2}=-\frac{\mu}{16|T|^2}
\end{equation}
and
\begin{equation}
c_3=\left(-\frac{8A_2-\ln
B_1}{2|T|^2}-C_3\right)B_1^2=-\frac{\mu\left(1+\ln\frac{4|T|^2}{\mu}
\right)}{16|T|^2}.
\end{equation}
\end{theorem}

\emph{Proof:}\\

The proof of the theorem has 4 steps and is presented in the subsections 3.1-3.4.

The following notation
\[
\system{l}{a_j=\Re(\alpha_j),\\ b_j=\Im(\alpha_j)}
\]
will be used throughout the rest of the paper.

\subsection{Simplification of $g'(z)$.}

We start with simplifying the expression for $g'(z)$ for $z$
values away from $0$ and $\pm\mu/2$
\begin{equation}
g'(z)=\frac{R(z)}{2\pi i}\int_{\gg_0\cup\gg_m}
\frac{f'(\xi)}{(\xi-z)R_+(\xi)}d\xi, \label{g'_main_arcs}
\end{equation}
where $z$ is off the main arcs $\gamma_0$ and $\gamma_m$.
To overcome the difficulty of $R$ in the denominator taking near-zero values
as the $\alpha_j$'s approach the real axis,
we transform (\ref{g'_main_arcs}) based on the following simple lemma
\begin{lemma} \label{lemma_moment}
Let $C$ be a closed rectifiable curve in the complex plane. Define
\begin{equation}
G(z)=\int_{C} \frac{F(\xi)}{\xi-z}d\xi,
\end{equation}
where $F(z)$ is any continuous (on $C$) function that satisfies the moment condition
\begin{equation}
\int_{C}F(\xi)d\xi=0. \label{eqn:lemma_moment}
\end{equation}
Then
\begin{equation}
G(z)=\frac{1}{z-z_0}\int_{C} \frac{(\xi-z_0)F(\xi)}{(\xi-z)}d\xi,
\end{equation}
where $z_0$ is off the contour $C$ and $z_0\neq z$.
\end{lemma}

The proof follows from the simple identity
\begin{equation}
\frac{1}{\xi-z}=\frac{\xi-z_0}{(z-z_0)(\xi-z)}-\frac{1}{z-z_0},
\end{equation}
which together with the moment condition (\ref{eqn:lemma_moment}), proves the lemma.

First we utilize the moment conditions (\ref{2.1:moment_conditions_g}) with the contour
of integration placed along the branchcuts $\gg_0\cup\gg_m$
\begin{equation}
\int_{\gg_0\cup\gg_m} \frac{\xi^k f'(\xi)}{R_+(\xi)}d\xi=0, \ \ \ \ \
\ k=0,1,2,3. \label{eq1.4:moments}
\end{equation}
These moment conditions can be rewritten in the form
\begin{equation}
\system{l}{ \int_{\gg_0\cup\gg_m} \frac{f'(\xi)}{R_+(\xi)}d\xi=0 \\
\int_{\gg_0\cup\gg_m} \frac{\Lambda(\xi)f'(\xi)}{R_+(\xi)}d\xi=0 \\
\int_{\gg_0\cup\gg_m} \frac{\Lambda(\xi)f'(\xi)}{R_+(\xi)}d\xi=0 \\
\int_{\gg_0\cup\gg_m} \frac{\Lambda(\xi)f'(\xi)}{R_+(\xi)}d\xi=0 }, \label{eq1.4:moments}
\end{equation}
where
\begin{equation}
\system{l}{ \Lambda(\xi)=\left(\xi-a_2\right), \\
\Lambda_2(\xi)=\left(\xi-a_0\right)\left(\xi-a_4\right), \\
\Lambda_3(\xi)=\left(\xi-a_0\right)\left(\xi-a_2\right)\left(\xi-a_4\right). }
\end{equation}

Thus by applying Lemma \ref{lemma_moment} to equation (\ref{g'_main_arcs}), we obtain
\begin{equation}
g'(z)=\frac{1}{2\pi i}\frac{R(z)}{\Lambda_3(z)}\int_{\gg_0\cup\gg_m}
\frac{\Lambda_3(\xi)f'(\xi)}{(\xi-z)R_+(\xi)}d\xi.
\end{equation}

\begin{figure}
\begin{center}
\includegraphics[height=4cm]{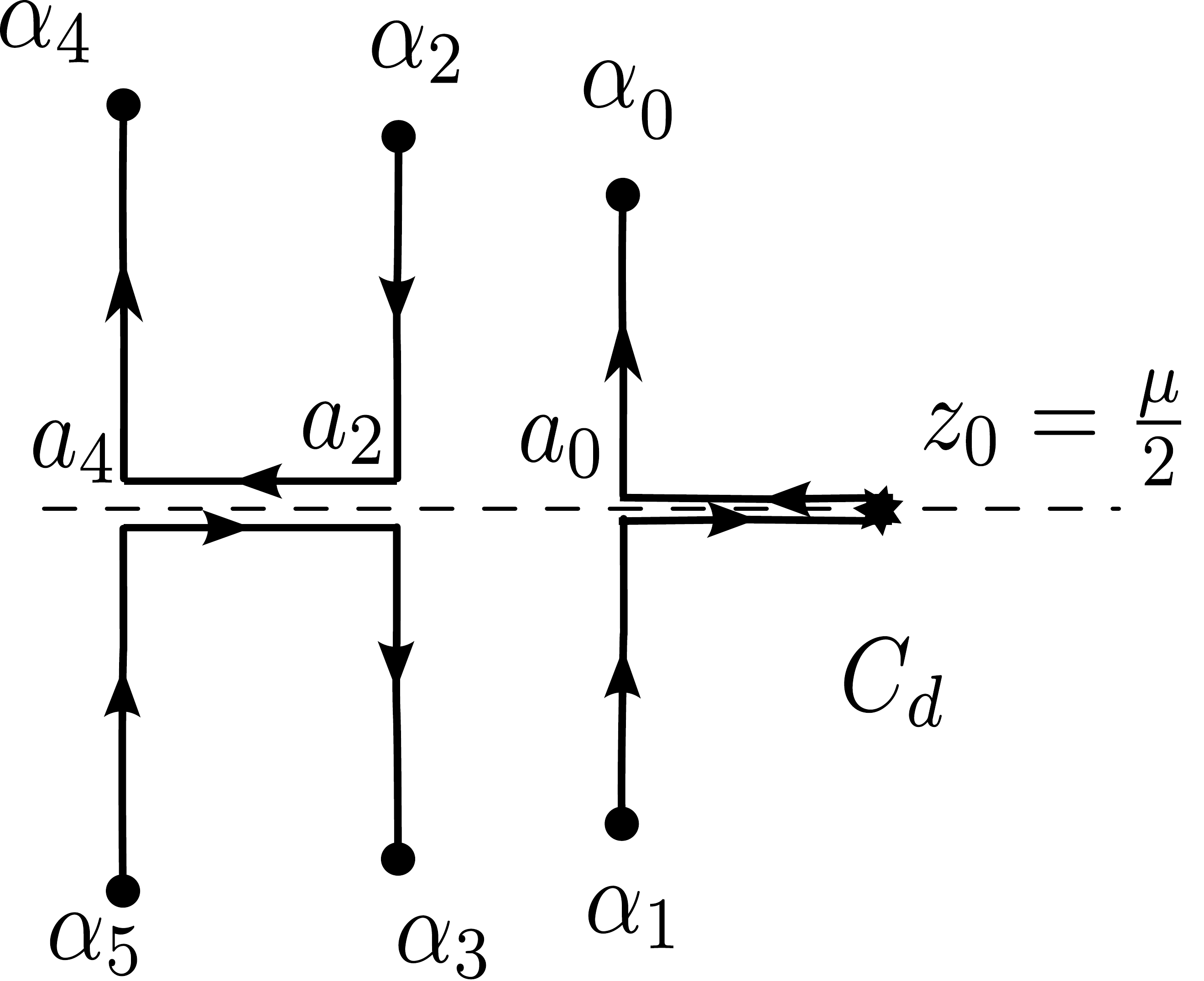}
\caption{\label{fig:deformed_contour} Deformed contour $C_d$.}
\end{center}
\end{figure}

Next we deform the contour of integration to the union of oriented arcs
in the upper half plane. Since $f$ in not analytic on the real axis,
$a_0,a_2, a_4$ as points in the complex plane, are understood as limits from
the upper half plane: $a_0+i0$, $a_2+i0$, $a_4+i0$.
Then $\gg_0\cup\gg_m$ is deformed into
\[
\left(\frac{\mu}{2},a_0\right) \cup \left(a_0,a_0+ib_0\right) \cup
\left(a_2+ib_2,a_2\right) \cup \left(a_2,a_4\right) \cup
\left(a_4,a_4+ib_4\right)
\]
and its complex conjugate (with the opposite orientation). We call the new
contour $C_d$ (see Fig. \ref{fig:deformed_contour}). Then
\begin{equation}
g'(z)=\frac{1}{2\pi i}\frac{R(z)}{\Lambda_3(z)}\int_{C_d}
\frac{\Lambda_3(\xi)f'(\xi)}{(\xi-z)R(\xi)}d\xi, \label{3.1:g'_C_d_reduction}
\end{equation}
and the moment conditions (\ref{eq1.4:moments}) become
\begin{equation}
\left\{\begin{array}{l}
\int_{C_d} \frac{\Lambda_2(\xi) f'(\xi)}{R(\xi)}d\xi=0\\
\int_{C_d} \frac{\Lambda_3(\xi) f'(\xi)}{R(\xi)}d\xi=0
\end{array}\right.\label{M14} .
\end{equation}

By Claim \ref{claim_1}, $\alpha_0$ and $\alpha_4$ converge to $\pm \mu/2$
exponentially fast in the long time limit. By Corollary \ref{cor_1}, the contributions from
the intervals $\left(\frac{\mu}{2},a_0\right)$, $\left(a_0,a_0+ib_0\right)$,
$\left(a_4,a_4+ib_4\right)$ and their conjugates in the contour of integration $C_d$
can be neglected, as they are exponentially small as $t\to\infty$,
and the interval $(a_2,a_4)$ can be replaced with $\left(a_2,-\frac{\mu}{2}\right)$.
Thus (\ref{3.1:g'_C_d_reduction}) reduces to
\begin{equation}
g'(z)=\frac{1}{2\pi i}\frac{R(z)}{\Lambda_3(z)}\int_{
\left(a_2,-\frac{\mu}{2}\right)\cup\overline{\left(-\frac{\mu}{2},a_2\right)}\cup\left(\alpha_2,\overline{\alpha}_2\right)}
\frac{\Lambda_3(\xi)f'(\xi)}{(\xi-z)R(\xi)}d\xi+exp.small,
\end{equation}
which, by recalling $f'_+-f'_-=\pi i$ on
$\left(a_2,-\frac{\mu}{2}\right)$, transforms into
\begin{equation}
g'(z)=\frac{R(z)}{2\Lambda_3(z)}\left[ \int_{a_2}^{-\frac{\mu}{2}}
\frac{\Lambda_3(\xi)}{(\xi-z)R(\xi)}d\xi + \frac{1}{\pi i}
\int_{\alpha_2}^{\overline{\alpha}_2}
\frac{\Lambda_3(\xi)f'(\xi)}{(\xi-z)R(\xi)}d\xi \right]+exp.small.
\end{equation}
Another simplifying observation is
\begin{equation}
\left\{
\begin{array}{l}
\frac{\Lambda_2(\xi)}{R(\xi)}=-\frac{1}{\sqrt{\left(\xi-\alpha_2\right)\left(\xi-\overline{\alpha}_2\right)}}+exp.small\\
\frac{\Lambda_3(\xi)}{R(\xi)}=-\frac{\xi-a_2}{\sqrt{\left(\xi-\alpha_2\right)\left(\xi-\overline{\alpha}_2\right)}}+exp.small,
\end{array}\right.
\end{equation}
Thus, $g'(z)$ in terms of real and imaginary parts of $\alpha_2=a_2+ib_2$ is
\[
g'(z)=-\frac{R(z)}{2\Lambda_3(z)}\left[ \int_{a_2}^{-\frac{\mu}{2}}
\frac{(\xi-a_2)}{(\xi-z)\sqrt{\left(\xi-a_2\right)^2+b_2^2}}d\xi
\right.
\]
\begin{equation}
\left.+\frac{1}{\pi i} \int_{a_2+ib_2}^{a_2-ib_2}
\frac{(\xi-a_2)f'(\xi)}{(\xi-z)\sqrt{\left(\xi-a_2\right)^2+b_2^2}}d\xi
\right]+exp.small. \label{DGas}
\end{equation}
Similarly for the last two moment conditions in (\ref{M14}) we
obtain
\begin{equation}
\left\{
\begin{array}{l}
\int_{a_2}^{-\frac{\mu}{2}}
\frac{1}{\sqrt{\left(\xi-a_2\right)^2+b_2^2}}d\xi +\frac{1}{\pi
i}\int_{a_2+ib_2}^{a_2-ib_2}
\frac{f'(\xi)}{\sqrt{\left(\xi-a_2\right)^2+b_2^2}}d\xi
+exp.small=0\\
\int_{a_2}^{-\frac{\mu}{2}}
\frac{\xi-a_2}{\sqrt{\left(\xi-a_2\right)^2+b_2^2}}d\xi
+\frac{1}{\pi i}\int_{a_2+ib_2}^{a_2-ib_2}
\frac{(\xi-a_2)f'(\xi)}{\sqrt{\left(\xi-a_2\right)^2+b_2^2}}d\xi
+exp.small=0\\
\end{array}\right.\label{M34as} .
\end{equation}

The expression for $g'(z)$ in (\ref{DGas}) and for the
moment conditions (\ref{M34as}) depend only on $\alpha_2=a_2+ib_2$
and not on $\alpha_0$, $\alpha_4$. All the
dependence on $\alpha_0$ and $\alpha_4$ is in the exponentially
small terms (abbreviated to "exp.small").

\subsection{Long time asymptotics of the branch point $\alpha_2$}

We now derive the asymptotics (\ref{thma2b2}). To simplify notations we call $a=a_2=o(1)$ and $b=b_2=o(1)$.
First we change the variable of integration in the second integral in (\ref{M34as})
$\xi=by+a$
\begin{equation}
\left\{
\begin{array}{l}
\left.
\ln\left(\xi-a+\sqrt{(\xi-a)^2+b^2}\right)\right|_a^{-\frac{\mu}{2}}
+\frac{1}{\pi i}\int_{i}^{-i} \frac{f'(by+a)}{\sqrt{y^2+1}}dy
+exp.small=0\\
\left.\sqrt{\left(\xi-a\right)^2+b^2}\right|_a^{-\frac{\mu}{2}}
+\frac{1}{\pi i}\int_{i}^{-i}
\frac{by\ f'(by+a)}{\sqrt{y^2+1}}dy+exp.small=0\\
\end{array}\right.\label{}
\end{equation}
with
\begin{equation}
f'(\xi)=-\frac{\pi
i}{2}-\ln\left(\frac{\mu}{2}-\xi\right)+\frac{1}{2}\ln\left(\xi^2+|T|^2\right)-x-4t\xi.
\label{df}
\end{equation}
The leading order of $f'$ as $t\to\infty$ comes from the last term $-4t\xi$ since $\xi\in[-i,i]$
\begin{equation}
f'(by+a)=-4t(by+a)+O(1),
\end{equation}
we obtain
\begin{equation}
\left\{
\begin{array}{l}
\ln\left(\sqrt{\left(\frac{\mu}{2}+a\right)^2+b^2}-\left(\frac{\mu}{2}+a\right)\right)
-\ln b+\frac{1}{\pi i}\int_{i}^{-i}
\frac{-4t(by+a)}{\sqrt{y^2+1}}dy+O(1)
=0\\
\sqrt{\left(\frac{\mu}{2}+a\right)^2+b^2}-b +\frac{1}{\pi
i}\int_{i}^{-i}
\frac{-4t(by+a)by}{\sqrt{y^2+1}}dy+O(b)=0\\
\end{array}\right.\label{}.
\end{equation}
By computing the integral $\int_{i}^{-i}
\frac{y}{\sqrt{y^2+1}}dy=0$, and expanding the square root for fixed $\mu>0$ and small $a$, $b$.
\begin{equation}
\left\{
\begin{array}{l}
\ln\left(\frac{b}{2(\frac{\mu}{2}+a)}+O(b^3)\right)
-\frac{4ta}{\pi i}\int_{i}^{-i} \frac{1}{\sqrt{y^2+1}}dy
+O(1)=0\\
\left(\frac{\mu}{2}+a\right) -\frac{4tb^2}{\pi i}\int_{i}^{-i}
\frac{y^2}{\sqrt{y^2+1}}dy+O(b)=0\\
\end{array}\right.\label{}.
\end{equation}
Then the integrals are evaluated explicitly $\int_{i}^{-i}
\frac{1}{\sqrt{y^2+1}}dy=-\pi i$, $\int_{i}^{-i}
\frac{y^2}{\sqrt{y^2+1}}dy=\frac{\pi i}{2}$ and the logarithm is expanded
\begin{equation}
\left\{
\begin{array}{l}
\ln b +4ta +O(1)=0\\
\frac{\mu}{2} - 2tb^2 +O(a+b)=0\\
\end{array}\right.\label{a2b2leading_2}.
\end{equation}
Solving the second equation for $b$ and plugging the result in the
first equation to solve for $a$, we obtain the leading order
asymptotics of $\alpha_2=a+ib$:
\begin{equation}
\left\{
\begin{array}{l}
a=\frac{\ln t}{8t}+O\left(\frac{1}{t}\right)\\
b=\sqrt{\frac{\mu}{4t}}+O\left(\frac{1}{t}\right)\\
\end{array}\right.\label{eq1.4:a2b2_lead}.
\end{equation}
As in the pure radiation case, $\alpha_2$ approaches to $0$ on
the scale $\frac{1}{\sqrt{t}}$ \cite{TVZlong2_06}.
\begin{remark}
The correction term in the second line of (\ref{eq1.4:a2b2_lead}) is of the order
$O\left(\frac{\ln t}{t^{3/2}}\right)$ since in
the second equation $O(a+b)$ is in fact $O(a)$ (see (\ref{b2_5}), (\ref{b2_6}) in Appendix).
\end{remark}
The higher order terms are obtained in the appendix:
\begin{equation}
\left\{
\begin{array}{l}
a=\Re(\alpha_2)=\frac{1}{8}\frac{\ln t}{t}+A_2\frac{1}{t}
 +A_3\frac{\ln t }{t^{2}}+O\left(\frac{1}{t^{2}}\right) \\
b=\Im(\alpha_2)=\sqrt{\frac{\mu}{4}}\frac{1}{\sqrt{t}}+B_2\frac{\ln
t}{t^{3/2}}
 +B_3\frac{1}{t^{3/2}}+O\left(\frac{\ln t }{t^{5/2}}\right).
\end{array}
\right.\label{a2b2_8}
\end{equation}
This result allows us to compute the long time asymptotics of $g'(z)$.

\subsection{Long time asymptotics of $g'(z)$}

Let $z$ be such that its distance $d(z)$ from $\alpha_0$ and the
interval $[a_2,a_4]$ satisfies $d(z)>>b_2$ as $t\to\infty$.
First, in (\ref{DGas}) we expand $\frac{1}{\xi-z}$ in powers of
$\frac{\xi}{z}$. This decomposition uniformly holds away from the branch points
$\alpha_2$ and $\overline{\alpha}_2$, i.e., for $|z|>|\alpha_2|$.
Then (\ref{DGas}) becomes
\[
g'(z)=\frac{\sqrt{\left(z-a\right)^2+b^2}}{2(z-a)}\left[ \int_{a}^{-\frac{\mu}{2}}
\frac{(\xi-a)}{(\xi-z)\sqrt{\left(\xi-a\right)^2+b^2}}d\xi \right.
\]
\begin{equation}
\left. -\frac{1}{\pi i z} \int_{a+bi}^{a-bi}
\frac{(\xi-a)f'(\xi)}{\sqrt{\left(\xi-a\right)^2+b^2}}\left(1+\frac{\xi}{z}
+\sum_{k=2}^{\infty}\left(\frac{\xi}{z}\right)^k \right)d\xi
\right]+exp.small.
\end{equation}
Taking into account the last moment condition in (\ref{M34as}) in
the form
\begin{equation}
\frac{1}{\pi i}\int_{a+bi}^{a-bi}
\frac{(\xi-a)f'(\xi)}{\sqrt{\left(\xi-a\right)^2+b^2}}d\xi =-
\int_{a}^{-\frac{\mu}{2}}
\frac{\xi-a}{\sqrt{\left(\xi-a\right)^2+b^2}}d\xi + exp.small
\end{equation}
we arrive at
\[
g'(z)=\frac{\sqrt{\left(z-a\right)^2+b^2}}{2(z-a)}\left[ \int_{a}^{-\frac{\mu}{2}}
\frac{(\xi-a)}{\sqrt{\left(\xi-a\right)^2+b^2}}\left(\frac{1}{\xi-z}+\frac{1}{z}\right)d\xi
\right.
\]
\begin{equation}
\left. -\frac{1}{\pi i z} \int_{a+bi}^{a-bi}
\frac{(\xi-a)f'(\xi)}{\sqrt{\left(\xi-a\right)^2+b^2}}\left(\frac{\xi}{z}
+\sum_{k=2}^{\infty}\left(\frac{\xi}{z}\right)^k \right)d\xi
\right]+exp.small.
\end{equation}
To simplify the expressions, for the rest of this section we reuse
(since there is no $\alpha_0$, $\alpha_4$ dependence anymore) notation
$R(z)=-\sqrt{\left(z-a\right)^2+b^2}$ and $\Lambda(z)=z-a$ in the integrals
\[
g'(z)=\frac{1}{2}\ \frac{R(z)}{\Lambda(z)}
\int_{a}^{-\frac{\mu}{2}} \frac{\Lambda(\xi)}{R(\xi)}\
\frac{\xi}{(\xi-z)z}d\xi
\]
\begin{equation}
 -\frac{1}{2\pi i z}\ \frac{R(z)}{\Lambda(z)}
\int_{a+bi}^{a-bi}
\frac{\Lambda(\xi)}{R(\xi)}f'(\xi)\left(\frac{\xi}{z}
+\sum_{k=2}^{\infty}\left(\frac{\xi}{z}\right)^k \right)d\xi
+exp.small. \label{eq1.4:g'_before_split}
\end{equation}
There are two main objects to analyze: $\frac{R(z)}{\Lambda(z)}$ and
$\frac{\Lambda(\xi)}{R(\xi)}$. We deal with these ratios
separately by a well known trick
\begin{equation}
f_n\int g_n\ h= (f_n-f)\int (g_n-g) \ h + (f_n-f)\int g \ h
+ f\int (g_n-g) \ h + f\int g \ h \label{eq1.4:fgh_trick}.
\end{equation}
Using as guidelines, for $t\to\infty$
\begin{equation}
\frac{R(z)}{\Lambda(z)}=-\frac{\sqrt{\left(z-a\right)^2+b^2}}{z-a}\to
-1 \quad \mbox{as} \ t\to\infty,
\end{equation}
which we think of as $f_n\to f$ in (\ref{eq1.4:fgh_trick}) and
\begin{equation}
\frac{\Lambda(\xi)}{R(\xi)}
=-\frac{\xi-a}{\sqrt{\left(\xi-a\right)^2+b^2}} \to 1 \quad
\mbox{in} \ L_2\left(a,-\mm\right) \quad \mbox{as} \ t\to\infty,
\end{equation}
which plays the role of $g_n\to g$.

We do not prove the above statements but rather use them as
suggestions in transforming the integrals in
(\ref{eq1.4:g'_before_split})
\[
g'(z)=-\frac{1}{2}\int_{a}^{-\frac{\mu}{2}}
\frac{\xi}{(\xi-z)z}d\xi
+\frac{1}{2}\left(\frac{R(z)}{\Lambda(z)}+1\right)
\int_{a}^{-\frac{\mu}{2}} \frac{\xi}{(\xi-z)z} d\xi
\]
\[
-\frac{1}{2}\int_{a}^{-\frac{\mu}{2}}
\left(\frac{\Lambda(\xi)}{R(\xi)} -1
\right)\frac{\xi}{(\xi-z)z}d\xi
+\frac{1}{2}\left(\frac{R(z)}{\Lambda(z)}+1\right)
\int_{a}^{-\frac{\mu}{2}} \left( \frac{\Lambda(\xi)}{R(\xi)} -1
\right) \frac{\xi}{\xi-z}d\xi
\]
\[
 + \frac{1}{2 \pi i z} \int_{a+bi}^{a-bi}
\frac{\Lambda(\xi)}{R(\xi)} f'(\xi)\ \frac{\xi}{z}d\xi
 -\frac{1}{2 \pi i z} \left(\frac{R(z)}{\Lambda(z)}+1\right)\int_{a+bi}^{a-bi}
\frac{\Lambda(\xi)}{R(\xi)} f'(\xi)\ \frac{\xi}{z}d\xi
\]
\begin{equation}
 - \frac{1}{2 \pi i z}\ \frac{R(z)}{\Lambda(z)} \int_{a+bi}^{a-bi}
\frac{\Lambda(\xi)}{R(\xi)} f'(\xi)\left(
\sum_{k=2}^{\infty}\left(\frac{\xi}{z}\right)^k \right)d\xi
+exp.small. \label{eq1.4:g'_after_split}
\end{equation}
For accounting we label these integrals as $I_{1-6}$ and the terms
in the last infinite sum of integrals are labeled as $H_k$,
$k=2,3,...$. We keep terms of the order up to and including
$O\left(\frac{1}{t}\right)$. See appendix for detailed calculations
of the following result:
\begin{equation}
\left\{
\begin{array}{l}
I_1=-\frac{1}{2}\ln\left(-\mm-z\right)+\frac{\mu}{4z}+\frac{1}{2}\ln(-z)
+O(a^2),\\
I_2= \frac{b^2}{2z^2}I_1+O(ab^2),\\
I_3 = \frac{1}{4z^2}b^2\ln b +\left(\frac{1-2\ln\mu}{8z^2}
+\frac{\ln\left(1+\frac{\mu}{2z}\right)}{4z^2} \right)b^2 +O(ab),\\
I_4= O(b^4\ln b)=O\left(\frac{\ln t}{t^2}\right),\\
I_5=\frac{1}{z^2}\left(2tab^2-\frac{C_1b^2}{4}\right)+O(ab),\quad C_1=\ln \frac{2|T|}{\mu}-x,\\
I_6= O(t ab^4)=O\left(\frac{\ln t}{t^2}\right),\\
H_2 = - \frac{3tb^4}{4 z^3} + O(b^3),\\
H_3 = O(tb^{5})=O\left(\frac{1}{t^{3/2}}\right), \\
H_k = O(tb^{k+2})=O\left(\frac{1}{t^{3/2}}\right), \quad k=4,5,\ldots.
\end{array}
\right.
\end{equation}

Recall: $b=O\left(\frac{1}{\sqrt{t}}\right)$ and
$a=O\left(\frac{\ln t}{t}\right)$ as $t\to\infty$.

Now putting together these results
\begin{equation}
g'(z)=I_1+I_2+I_3+I_4+I_5+I_6+H_2 + O(H_3)
\end{equation}
\[
=I_1 +O(a^2)+ \frac{b^2}{2z^2}I_1+O(ab^2) + \frac{1}{4z^2}b^2\ln b
+\left(\frac{1-2\ln\mu}{8z^2}
+\frac{\ln\left(1+\frac{\mu}{2z}\right)}{4z^2} \right)b^2
\]
\begin{equation}
+O(ab)+O(b^4\ln
b)+\frac{1}{z^2}\left(2tab^2-\frac{C_1b^2}{4}\right)+O(ab)
+O(t ab^4) - \frac{3tb^4}{4 z^3} + O(b^3)=
\end{equation}
\[
=I_1 + \frac{1}{4z^2}b^2\ln b +\frac{2tab^2}{z^2} -
\frac{3tb^4}{4 z^3}
\]
\begin{equation}
+\frac{1}{z^2}\left(\frac{1}{2}I_1 + \frac{1-2\ln\mu}{8}
+\frac{\ln\left(1+\frac{\mu}{2z}\right)}{4} -\frac{C_1}{4}
\right)b^2   +O(ab).
\end{equation}
After substitution of the expression for $I_1$ and $C_1$
the logarithmic terms cancel
\[
g'(z)= \frac{1}{2}\ln(z) -\frac{1}{2}\ln\left(\mm+z\right)
+\frac{\mu}{4z} + \frac{1}{z^2}\left(\frac{b^2\ln b}{4}
+2tab^2\right)
\]
\begin{equation}
+\frac{1}{z^2}\left( \frac{1-2\ln|T|}{8} +\frac{x-\ln 2}{4}
\right)b^2 +\frac{1}{z^3}\left( -\frac{3tb^4}{4}+\frac{\mu
b^2}{8}\right)+O(ab),
\end{equation}
as $t\to\infty$. Note: $O(ab)=O\left(\frac{\ln
t}{t^{3/2}}\right)$.

\subsection{Long time asymptotics of the singular obstruction curve}

In this section we asymptotically solve $\Im h(T,x,t)=0$ which requires
asymptotics of $g(z)$.
By integrating $g'(z)$ to obtain $g(z)$ and using the fact
$h(z)=2g(z)-f(z)$ we write
\begin{equation}
g(z)=\int_{r_0}^z g'(s)ds + g(r_0)
\end{equation}
for some $r_0\in\mathbb{R}$, $r_0\neq\mm$, which implies $\Im g(r_0)=0$, that is $g(r_0)$ is
real. In particular we can send $r_0\to +\infty$ since $g(z)$ is
analytic at infinity. Integrating $g'$ and evaluating the limit
\begin{equation}
\lim_{r_0\to\infty}\left[ \frac{1}{2}r_0\ln r_0
-\frac{1}{2}\left(r_0+\mm\right)\ln\left(r_0+\mm\right)
+\frac{\mu}{4}\ln r_0\right]
\end{equation}
\begin{equation}
=\lim_{r_0\to\infty} -\frac{1}{2}\left(r_0+\mm\right)
\left(\frac{\mu}{2r_0}+O\left(\frac{1}{r_0^2}\right)\right)
=-\frac{\mu}{4},
\end{equation}
we obtain
\[
g(z)=\frac{1}{2}z\ln z
-\frac{1}{2}\left(z+\mm\right)\ln\left(z+\mm\right)
+\frac{\mu}{4}\ln z - \frac{1}{z}\left(\frac{b^2\ln b}{4}
+2tab^2\right)
\]
\[
-\frac{1}{z}\left( \frac{1-2\ln|T|}{8} +\frac{x-\ln 2}{4}
\right)b^2 -\frac{1}{2z^2}\left( -\frac{3tb^4}{4}+\frac{\mu
b^2}{8}\right)
\]
\begin{equation}
+ O(ab) + \frac{\mu}{4} + g(\infty) .\label{eq1.4:g(z)_final}
\end{equation}
Then for $0<\mu<2$ and $T=i\sqrt{1-\frac{\mu^2}{4}}$
\[
\Im g(T)=\frac{1}{2}\Im\left(T+\mm\right)\left[\ln T
-\ln\left(T+\mm\right) \right] + \frac{1}{|T|}\left(\frac{b^2\ln
b}{4} +2tab^2\right)
\]
\begin{equation}
+\frac{1}{|T|}\left( \frac{1-2\ln|T|}{8} +\frac{x-\ln 2}{4}
\right)b^2 + O(ab).
\end{equation}
and the value of $f(T)$ is computed as the limiting value
\begin{equation}
\Im f(T)=
\Im \lim_{z\to
T}f(z)=\Im \left(\left(\frac{\mu}{2}-T\right)\left[\frac{\pi
i}{2}+\ln\left(\frac{\mu}{2}-T\right)\right]\right)+|T|\ln\left(2|T|\right)
-x|T|.
\end{equation}
Next we compute $\Im h(T)=2\Im g(T)-\Im f(T)$ and after some
algebra we arrive at
\[
\Im h(T)=|T|(x-\ln 2) + \frac{2}{|T|}\left(\frac{b^2\ln b}{4}
+2tab^2\right)
\]
\begin{equation}
+\frac{2}{|T|}\left( \frac{1-2\ln|T|}{8} +\frac{x-\ln 2}{4}
\right)b^2 + O(ab). \label{Im_h_at_T_asympt}
\end{equation}

Thus the equation
\begin{equation}
\Im h(T,x,t)=0
\end{equation}
is asymptotically solved as $t\to\infty$
\begin{equation}
x(t)=\ln 2 - \frac{1}{|T|^2}\left(\frac{b^2\ln b}{2}
+4tab^2\right)-\left( \frac{1-2\ln|T|}{4|T|^2}\right)b^2 + O(ab),
\end{equation}
which proves Theorem \ref{thm_1}.
\begin{remark}
An interesting observation that in terms of $\alpha_2=a_2+ib_2$ (\ref{a2b2_8})
\begin{equation}
x(t)=\ln 2+c_2\frac{\ln t}{t}+c_3\frac{1}{t}+ O\left(\frac{\ln
t}{t^{3/2}}\right), \ \ \ \ t\to\infty,
\end{equation}
where
\begin{equation}
c_2=-\frac{\left(8A_1-\frac{1}{2}\right)B_1^2}{2|T|^2},\quad \quad
c_3=\left(-\frac{8A_2-\ln B_1}{2|T|^2}-C_3\right)B_1^2.
\end{equation}
To compute $c_2$ only $A_1$ and $B_1$ are needed and
to compute $c_3$ only $A_2$ additionally required. This may indicate
that the terms of the order $O\left(\frac{\ln
t}{t^k}\right)$ and $O\left(\frac{1}{t^k}\right)$ could be combined
together in computations.

Then our conjecture is that the next terms are of the orders
$O\left(\frac{\ln t}{t^{3/2}}\right)$, $O\left(\frac{1}{t^{3/2}}\right)$ and
possibly $O\left(\frac{\ln^2 t }{t^{3/2}}\right)$ in the singular obstruction curve $x(t)$
long time asymptotics. To compute these terms only coefficients: $A_1-A_2$, $A_3-A_4$ and
$B_1$, $B_2-B_3$ in the asymptotics of $\alpha_2(t)$ would be utilized.
\end{remark}

\section{Numerical computations}

\begin{figure}
\begin{center}
\includegraphics[height=10cm]{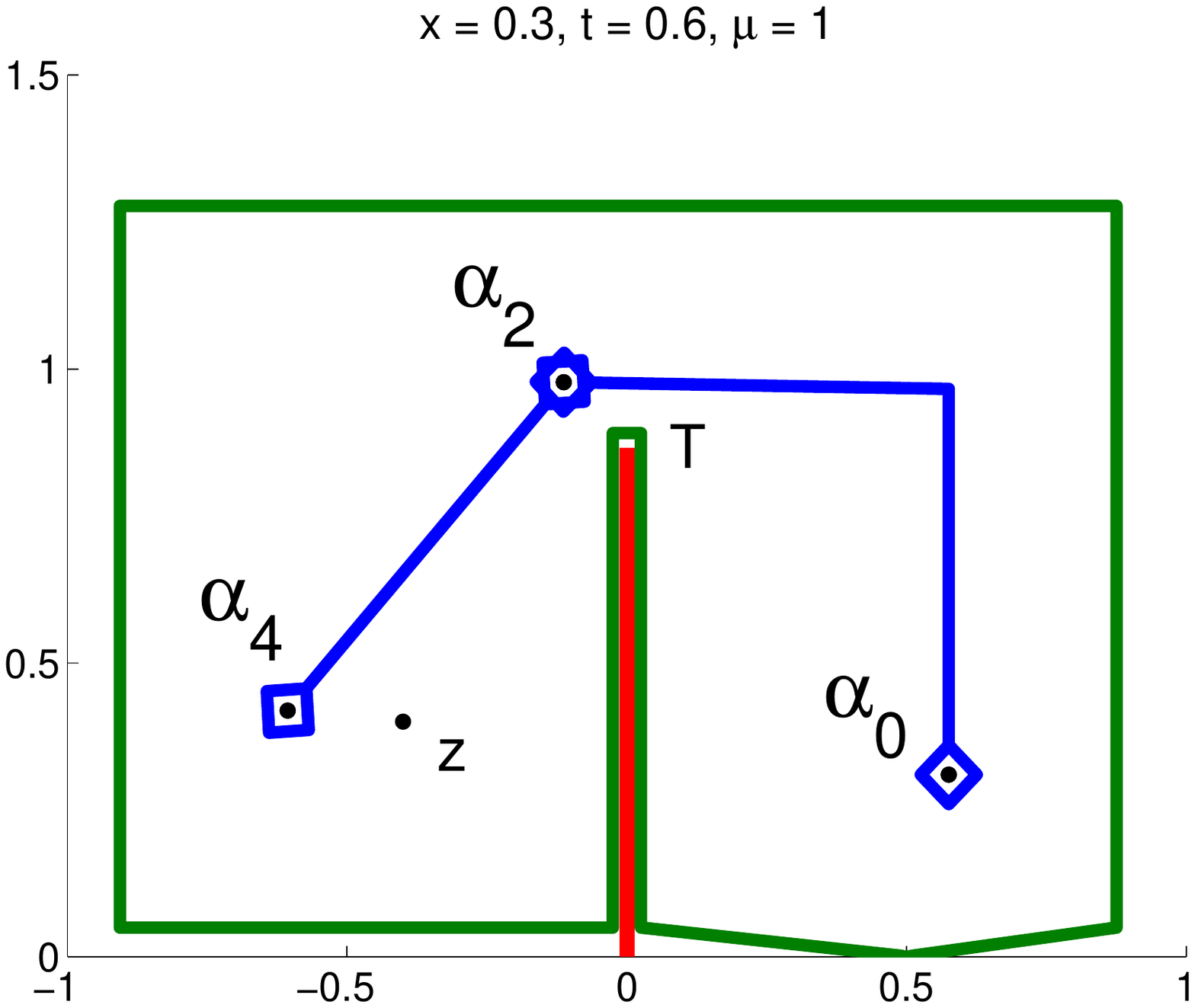}
\caption{\label{fig3:h_cont_g2} Matlab realization of contours of integration for $h(z)$
in genus 2 (see Fig. \ref{fig:h_fcn_cont}).
}
\end{center}
\end{figure}

\begin{figure}
\begin{center}
\includegraphics[height=4cm]{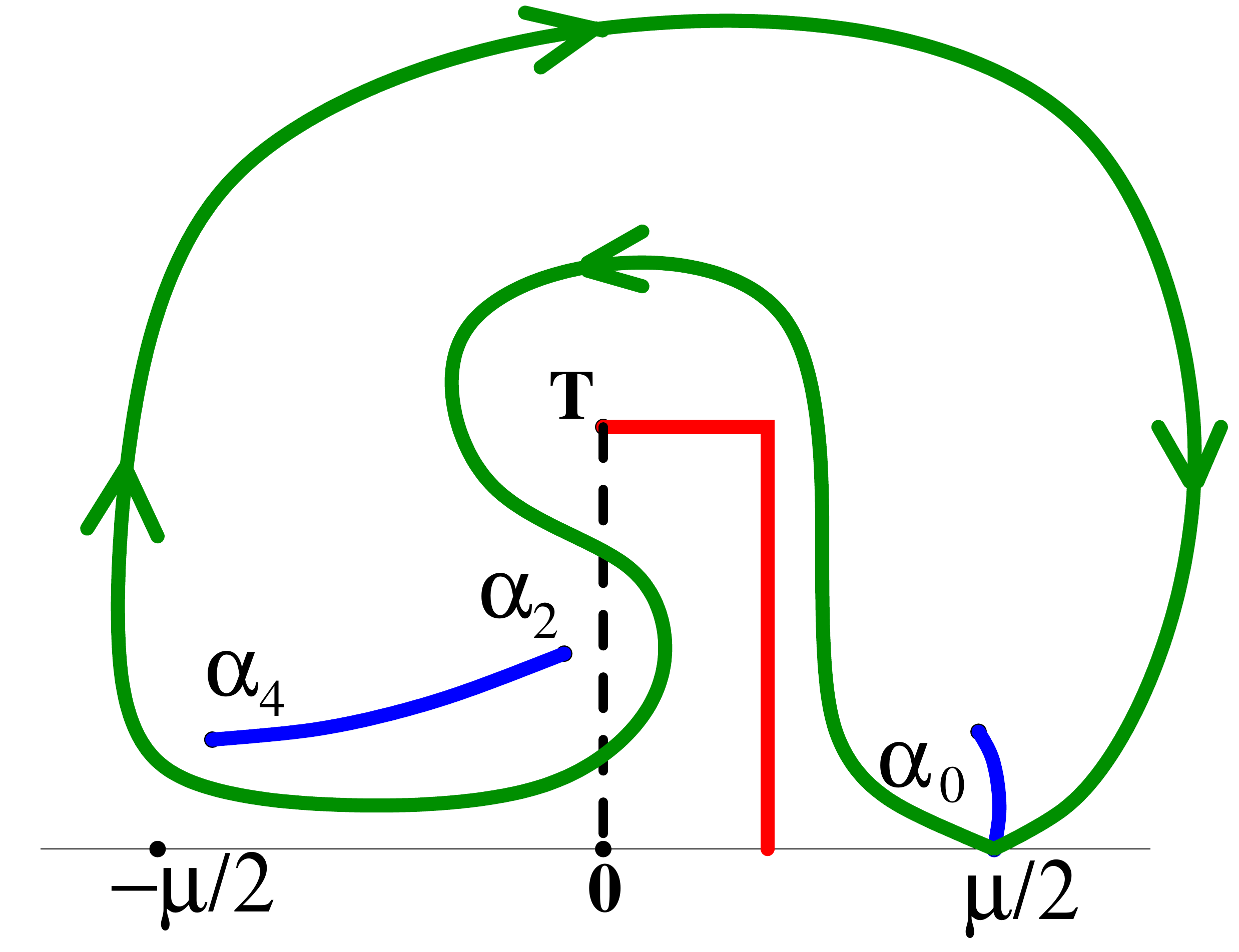}
\includegraphics[height=4cm]{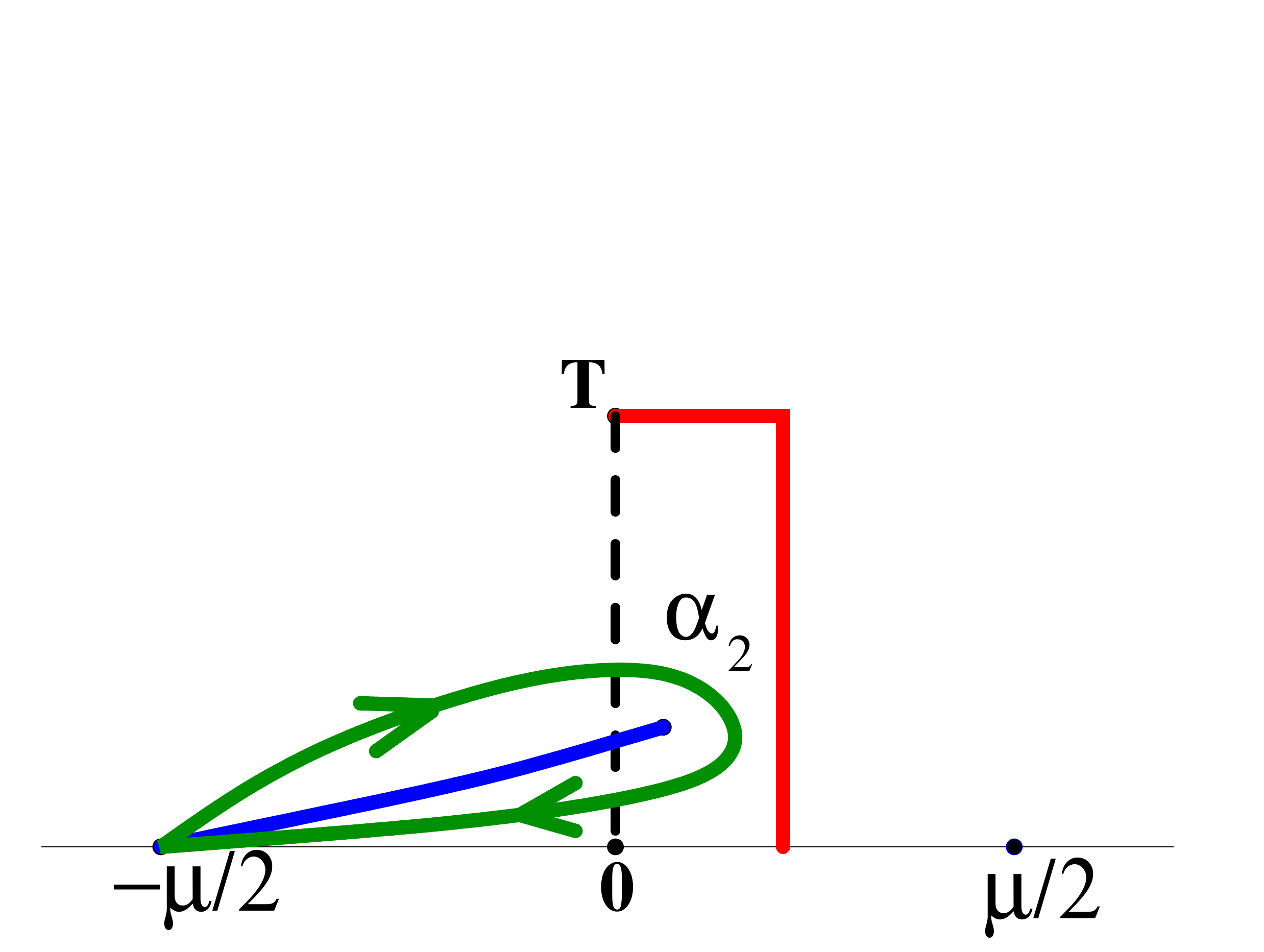}
\caption{\label{fig3:h_cont_g2LT} Contours of integration for
$h'(z)$ in genus 2. Large loop contour $\ggh$ is on the left.
Simplified contour of integration $\ggh_r$ for large values of $t$ is on
the right.}
\end{center}
\end{figure}

\begin{figure}
\begin{center}
\includegraphics[height=10cm]{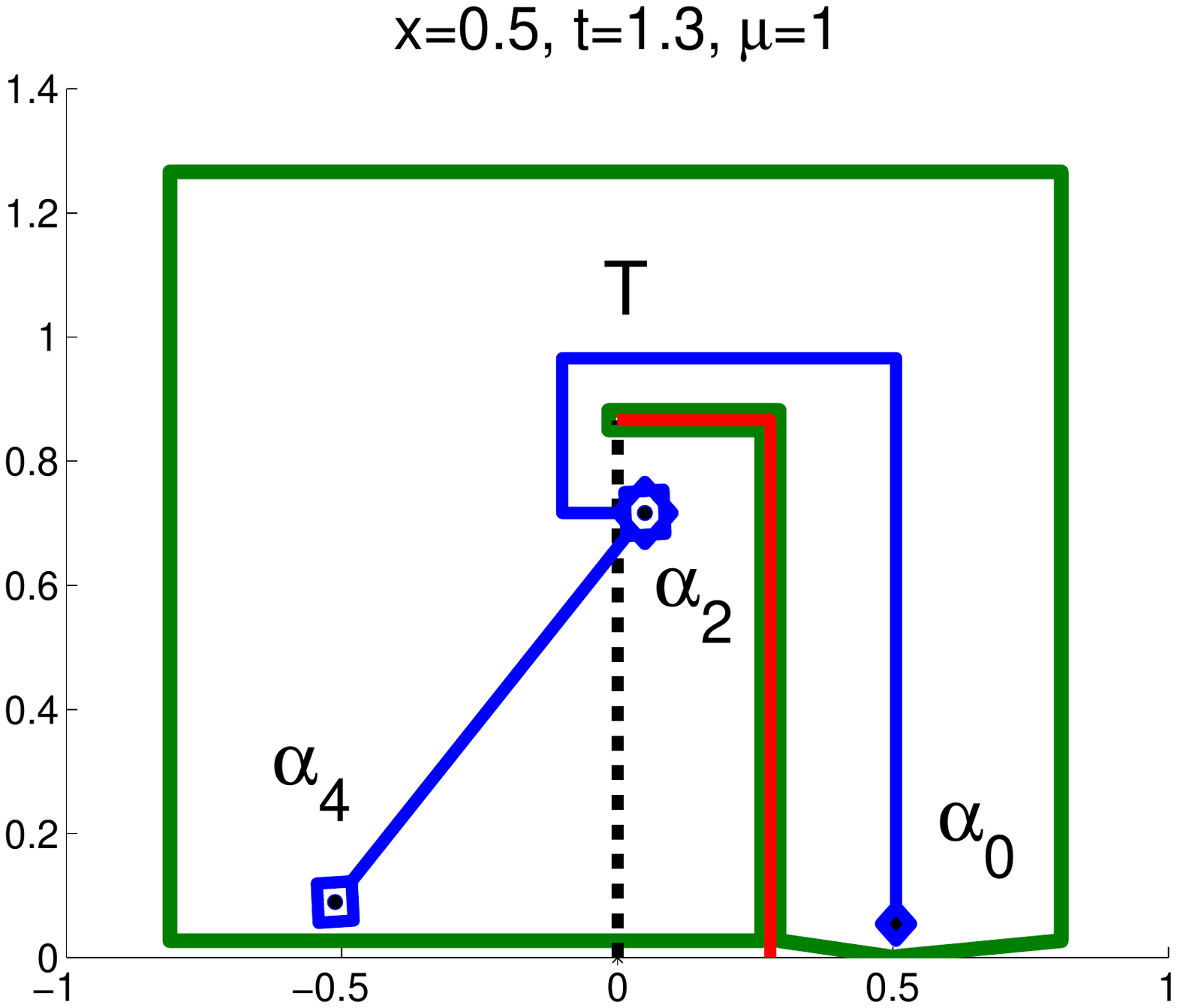}
\includegraphics[height=10cm]{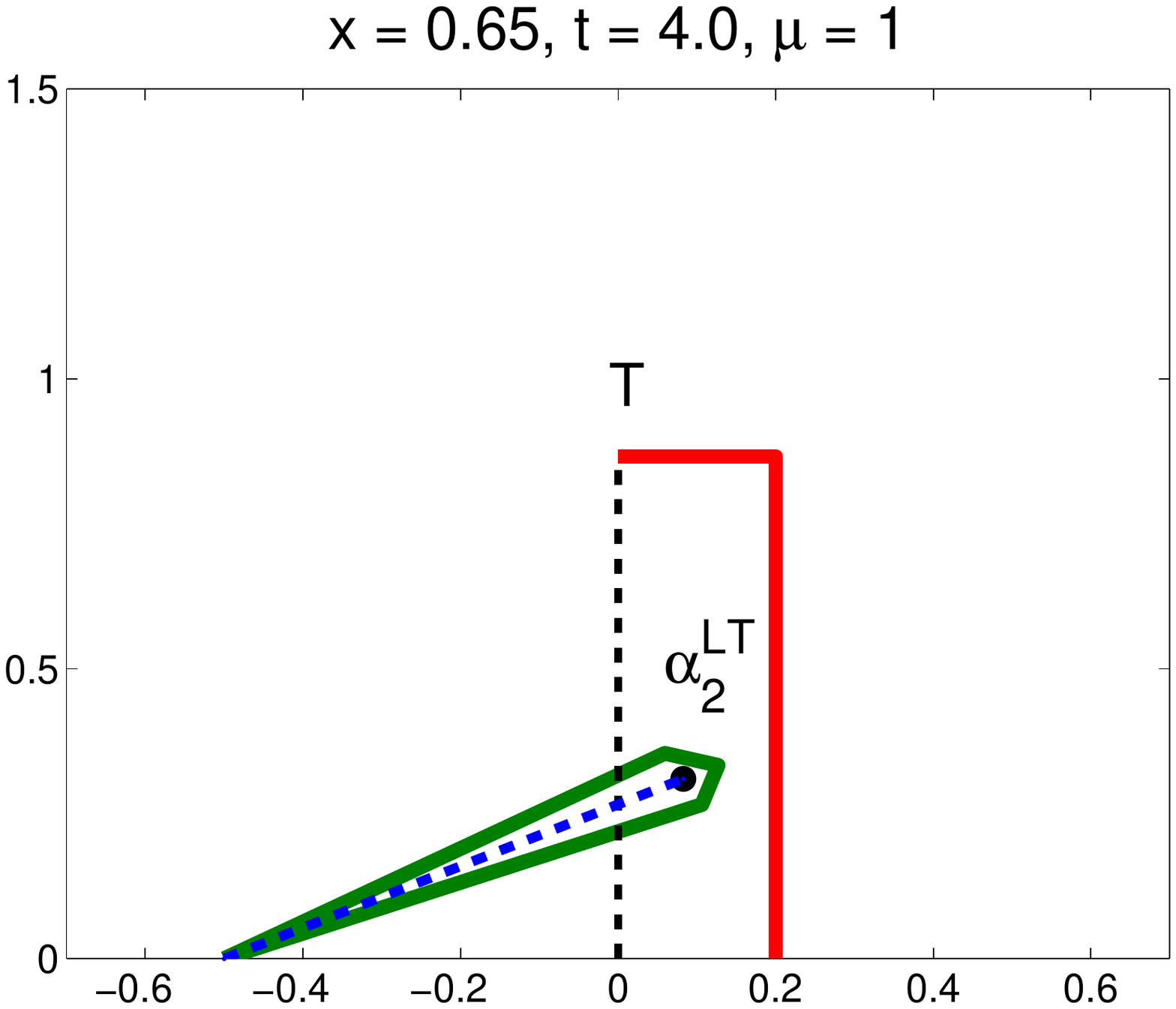}
\caption{\label{fig3:matlab_cont_dg2LT} Matlab realization of contours of
integration $\ggh$ for large $t$ values in genus 2 (left). Matlab realization of
contours of integration $\ggh_r$ for long time computations in
genus 2 (right).}
\end{center}
\end{figure}

\subsection{Numerical computation the singular obstruction curve}

\begin{figure}
\begin{center}
\includegraphics[height=10cm]{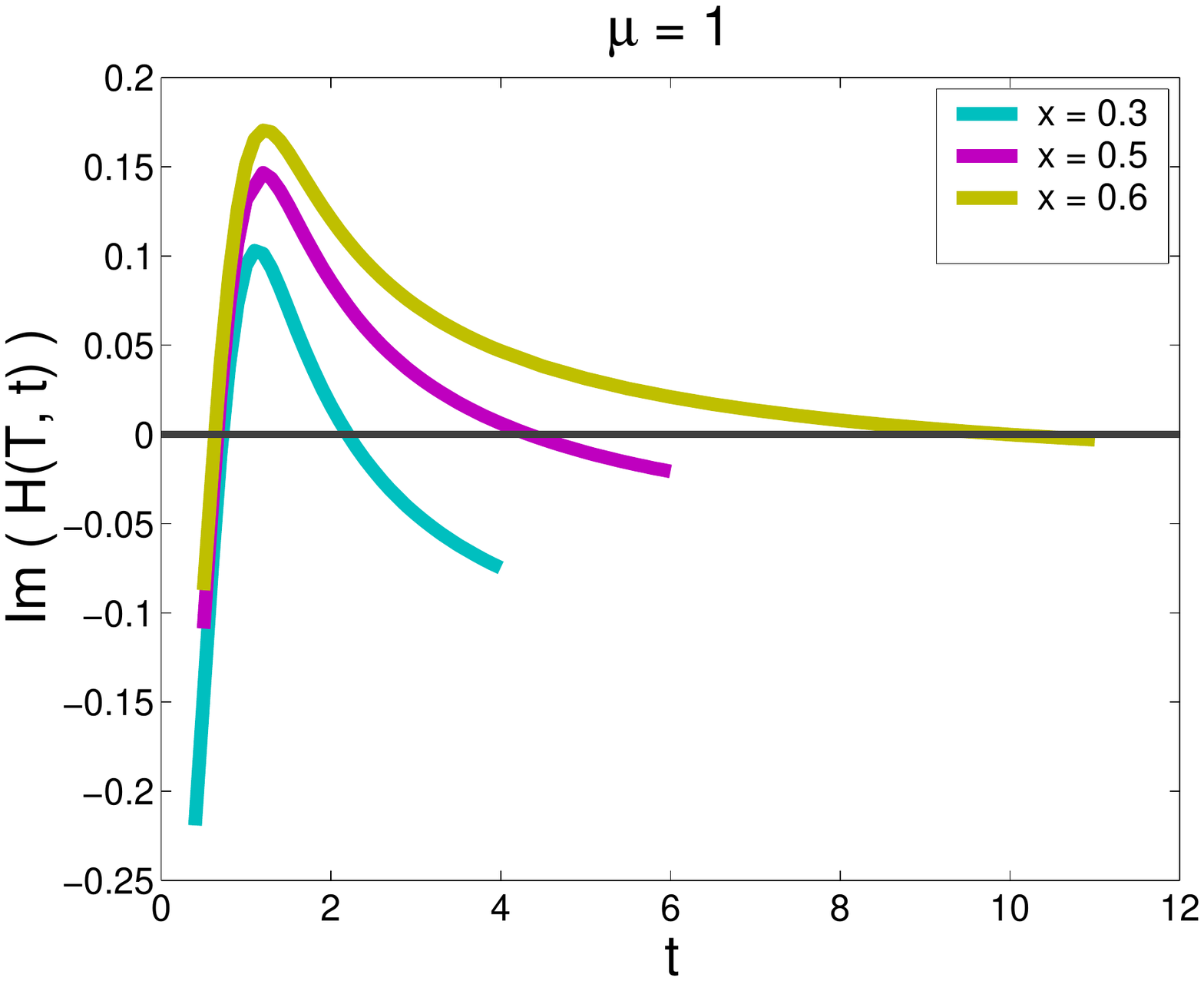}
\caption{\label{fig3:imag_h_at_T} $\Im h(T,x,t,\mu)$ as a function
of time $t$.}
\end{center}
\end{figure}

As discussed in Section 2 the singular obstruction for a fixed $x$ is defined as
one of the roots of the equation
\begin{equation}
\Im h(T, x, t_s(x)) = 0 \label{eq3.7:catastr_eqn}
\end{equation}
where the function $\Im h(T,t)$ is a nice function of $t$ (see
Fig. \ref{fig3:imag_h_at_T}).

From condition on $z$ being outside of the contour of integration
$\ggh$, $h(T)$ is not computable directly. Instead
\begin{equation}
\Im h(T)=\Im \left(2g(T)-f(T)\right)=0,
\end{equation}
where $f(T)$ is understood as limit, while function $g(z)$ is analytic at $z=T$.

For large $t$ values, the branch point $\alpha_2$ approaches the
imaginary axis below the point $z=T$ and hits a vertical branch cut
$[0,T]$ of function $f$. Using integration on a Riemann surface
this event is no special and $\alpha_2$ continues moving on
another sheet of the Riemann surface without any obstacles.

\subsection{Numerical long time computations}

Fast convergence of $\alpha_0$ and $\alpha_4$ to the real axis creates
challenges in numerical evaluation of $g$-function for large values of $t$.
We modify the computations by incorporating our asymptotic analysis.

The branch points $\alpha_0$, $\alpha_2$, and $\alpha_4$ are
singularities of the integrands (in the computation
of $h(z)$, $h'(z)$, $B(z)$), and need to be encircled by a contour
in the upper half plane. This requirement puts the branch points
close to the contour of integration causing falling in
accuracy. We utilize here that in the genus 2 region $\alpha_0$
and $\alpha_4$ are exponentially close to $\pm \frac{\mu}{2}$ respectively
as $t\to\infty$ and corresponding integrals are exponentially small.

The three complex branch points
$\alpha_0$, $\alpha_2$, and $\alpha_4$ satisfy the set of $4$
moment conditions, which are necessary but not sufficient to
compute the $\alpha$'s:
\begin{equation}
\left\{ \begin{array}{l}
\oint_{\ggh}\frac{ f'(\xi)}{R(\xi)}d\xi=0\\
\oint_{\ggh}\frac{\xi f'(\xi)}{R(\xi)}d\xi=0\\
\oint_{\ggh}\frac{\xi^2 f'(\xi)}{R(\xi)}d\xi=0\\
\oint_{\ggh}\frac{\xi^3 f'(\xi)}{R(\xi)}d\xi=0,
\end{array}\right.\label{eq:M_14}
\end{equation}
where
\begin{equation}
R(\xi)=\sqrt{\left(\xi-\alpha_0\right)\left(\xi-\overline{\alpha}_0\right)
\left(\xi-\alpha_2\right)\left(\xi-\overline{\alpha}_2\right)
\left(\xi-\alpha_4\right)\left(\xi-\overline{\alpha}_4\right)}.
\end{equation}
with the branch cuts chosen along the main arcs:
$\left[\mm,\alpha_0\right]$, $\left[\alpha_2,\alpha_4\right]$ and
their complex conjugates.

By considering linear combinations of the moment conditions (\ref{eq:M_14}), the
last two of these conditions can be written as
\begin{equation}
\left\{ \begin{array}{l}
\oint_{\ggh}\frac{\left(\xi-\mm\right)\left(\xi+\mm\right) f'(\xi)}{R(\xi)}d\xi=0\\
\oint_{\ggh}\frac{\xi \left(\xi-\mm\right)\left(\xi+\mm\right)
f'(\xi)}{R(\xi)}d\xi=0.
\end{array}\right. \label{eq:alpha_LT_system}
\end{equation}
Since $\alpha_0\to\mm$ and $\alpha_4\to-\mm$, we modify the contour of
integration accordingly. The large loop $\ggh$
is reduced to a smaller loop around $\left[-\mm,\alpha_2\right]$
and its complex conjugate which we call $\ggh_{r}$ (see Fig \ref{fig3:h_cont_g2LT}).
This leads to improvement of speed and stability and simplifies the system
(\ref{eq:alpha_LT_system})
\begin{equation}
\left\{ \begin{array}{l}
\oint_{\ggh_r}\frac{\left(\xi-\mm\right)\left(\xi+\mm\right) f'(\xi)}{R(\xi)}d\xi+exp.small=0\\
\oint_{\ggh_r}\frac{\xi \left(\xi-\mm\right)\left(\xi+\mm\right)
f'(\xi)}{R(\xi)}d\xi+exp.small=0,
\end{array}\right. \label{eq:gg_r_syst}
\end{equation}
where
\begin{equation}
R(\xi)=\left(\xi-\mm\right)\left(\xi+\mm\right)\sqrt{
\left(\xi-\alpha_{2}\right)\left(\xi-\overline{\alpha}_2\right)}
\end{equation}
with the branch cut chosen to connect $\alpha_2$ and
$\overline{\alpha}_2$ through $-\mm$. After cancelations, the
moment conditions (\ref{eq:gg_r_syst}) look similar to the case of genus 0 with one
unknown branch point $\alpha_{2}$. All dependence on $\alpha_{0}$ and $\alpha_{4}$
is in the exponentially small terms. Then we approximate $\alpha_{2}$
with $\alpha_{2}^{LT}$ which satisfies a system of two real equations
\begin{equation}
\left\{ \begin{array}{l}
\oint_{\ggh_r}\frac{f'(\xi)}{R_{LT}(\xi)}d\xi=0\\
\oint_{\ggh_r}\frac{\xi f'(\xi)}{R_{LT}(\xi)}d\xi=0,
\end{array}\right. \label{eq:alpha2_LT_system}
\end{equation}
where
\begin{equation}
R_{LT}(\xi)=\sqrt{\left(\xi-\alpha_{2}^{LT}\right)\left(\xi-\overline{\alpha}_2^{LT}\right)}
\end{equation}
with the branch cut chosen to connect $\alpha_2^{LT}$ and
$\overline{\alpha}_2^{LT}$ through $-\mm$.

One of the key advantages of these long time computations is
increased speed in exchange of precision in computing $\alpha_2$.
Solving the full system of $B$-function equations (\ref{2.1:B_system})
for $\alpha_0$, $\alpha_2$,
and $\alpha_4$ involves computing and inverting a $6$x$6$
matrix of partial derivatives for each iteration. While the long time
approximations of $\alpha_2^{LT}$ by solving (\ref{eq:alpha2_LT_system}) involve computing
and inverting a $2$x$2$ matrix of partial derivatives.

Using the same contour reduction, we compute $h'(z)$ in genus 2
for large $t$ values
\begin{equation}
h_{LT}'(z)=\frac{R_{LT}(z)}{2\pi i}\oint_{\ggh_r}\frac{
f'(\xi)}{(\xi-z)R_{LT}(\xi)}d\xi.
\end{equation}
Long-time computations of $\Im h(z)$ in genus 2 as
\begin{equation}
\Im h_{LT}(T)=\Im \int_{r_0}^{z} h'_{LT}(s)ds,
\end{equation}
where $r_0$ is some real number.

Long time computations of the singular obstruction curve are based on long
time computations of $\alpha_2^{LT}(t)$ and $\Im h_{LT}(z)$.
First, $\alpha_2^{LT}(t)$ is computed which is used
to evaluate $\Im h_{LT}(T,\alpha_2^{LT},t)$. The long-time approximation
of the singular obstruction $t_s^{LT}(x)$ then computed from
\begin{equation}
\Im h_{LT}\left(T,\alpha_2^{LT}\left(t_s^{LT}\right),t_s^{LT}\right)=0.
\end{equation}

Using our long time computations of $\alpha_2^{LT}$, $h'_{LT}(z)$,
and $\Im h_{LT}(z)$ we avoid the mentioned above difficulties with
accuracy and improve the speed of the computations in the case of
$\alpha_2^{LT}$ and $h'_{LT}(z)$.

\begin{figure}
\begin{center}
\includegraphics[height=10cm]{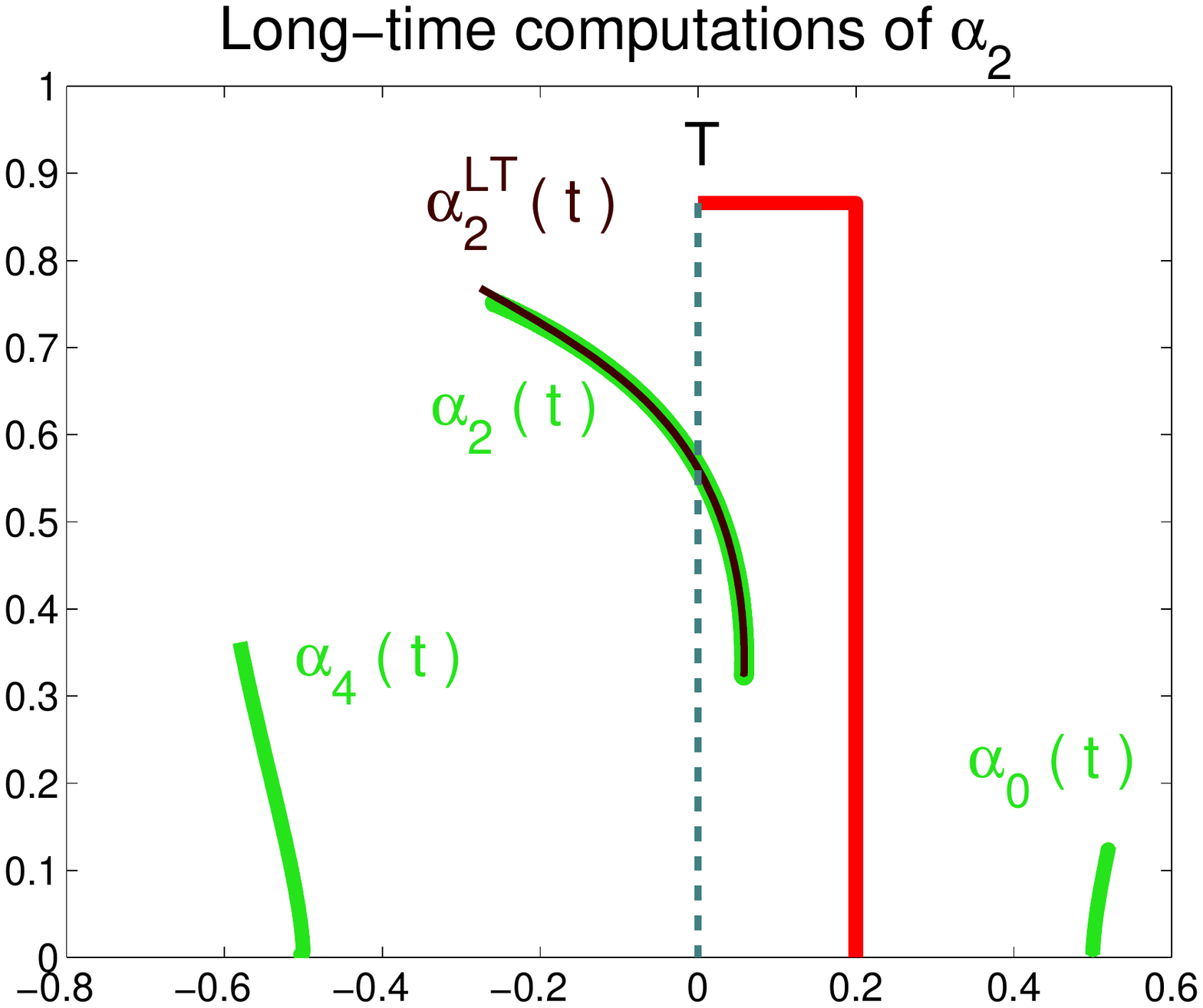}
\includegraphics[height=10cm]{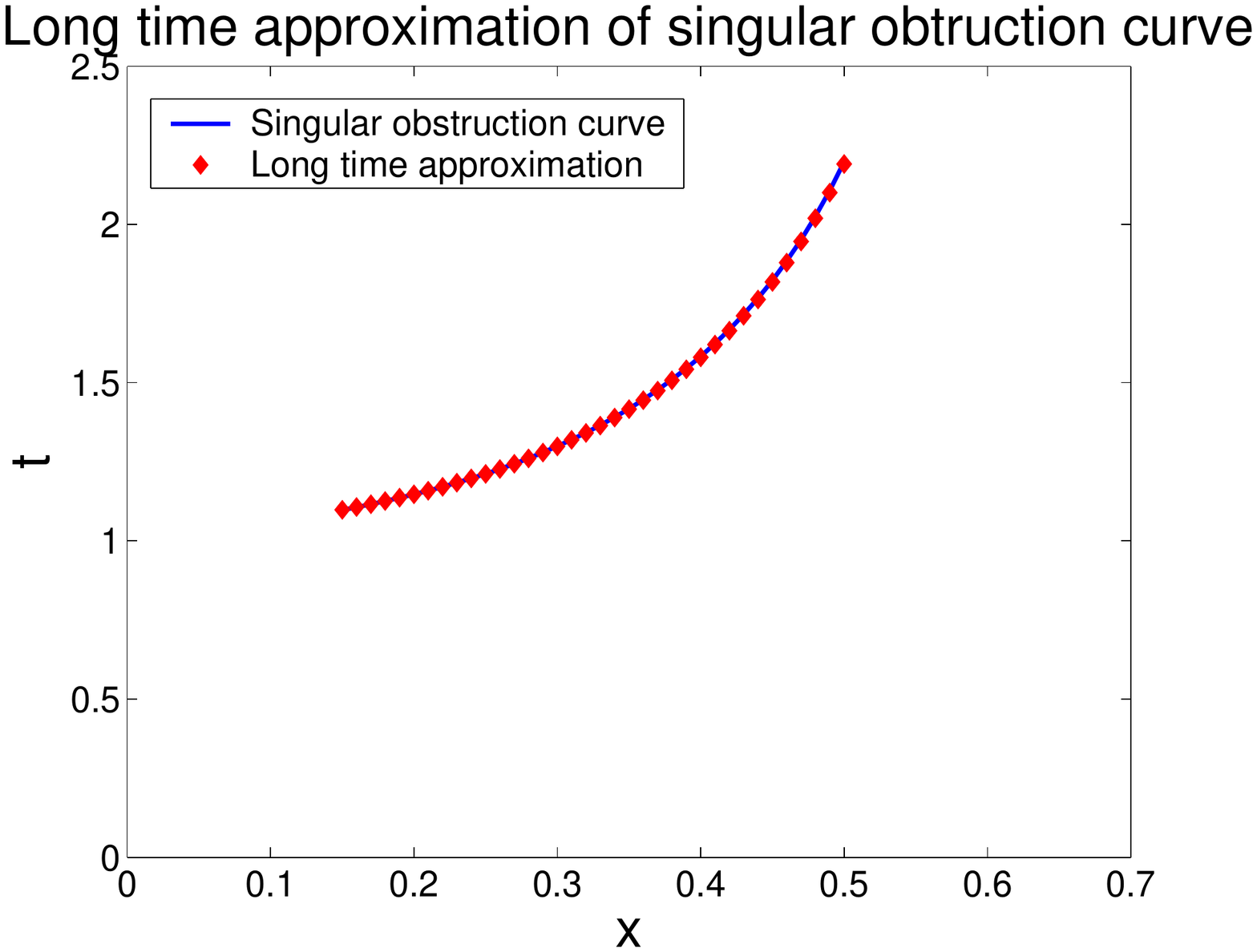}
\caption{\label{fig3:alfa_LT_t} Approximation of $\alpha_2(t)$
with $\alpha_2^{LT}(t)$ for $x=1.0$, $\mu=1$ in genus 2 (left). Approximation of the
singular obstruction $t_s(x)$ with $t_s^{LT}(x)$ for $\mu=1$ in
genus 2 (right).}
\end{center}
\end{figure}

The long time approximation $\alpha_2^{LT}$
of the correct value of $\alpha_2$ as a function of time $t$ are
presented in Figure \ref{fig3:alfa_LT_t} (left). It shows
the time evolution of $\alpha_0$, $\alpha_2$, $\alpha_4$
and $\alpha_2^{LT}$ for $x=1.0$ and $\mu=1$. $\alpha_2$ and
$\alpha_2^{LT}$ demonstrate similar and converging trajectories.
The difference convergence tends to zero as $t$ increases.

The long time approximation $t_s^{LT}(x)$ of the singular obstruction
curve $t=t_s(x)$ for $\mu=1$ is presented in Figure \ref{fig3:alfa_LT_t} (right).
It demonstrates good agreement for $t$ values as low as $1$.

\section{Discussion}\label{disc}

\subsection{First break}

\begin{figure}
\begin{center}
\includegraphics[height=7cm]{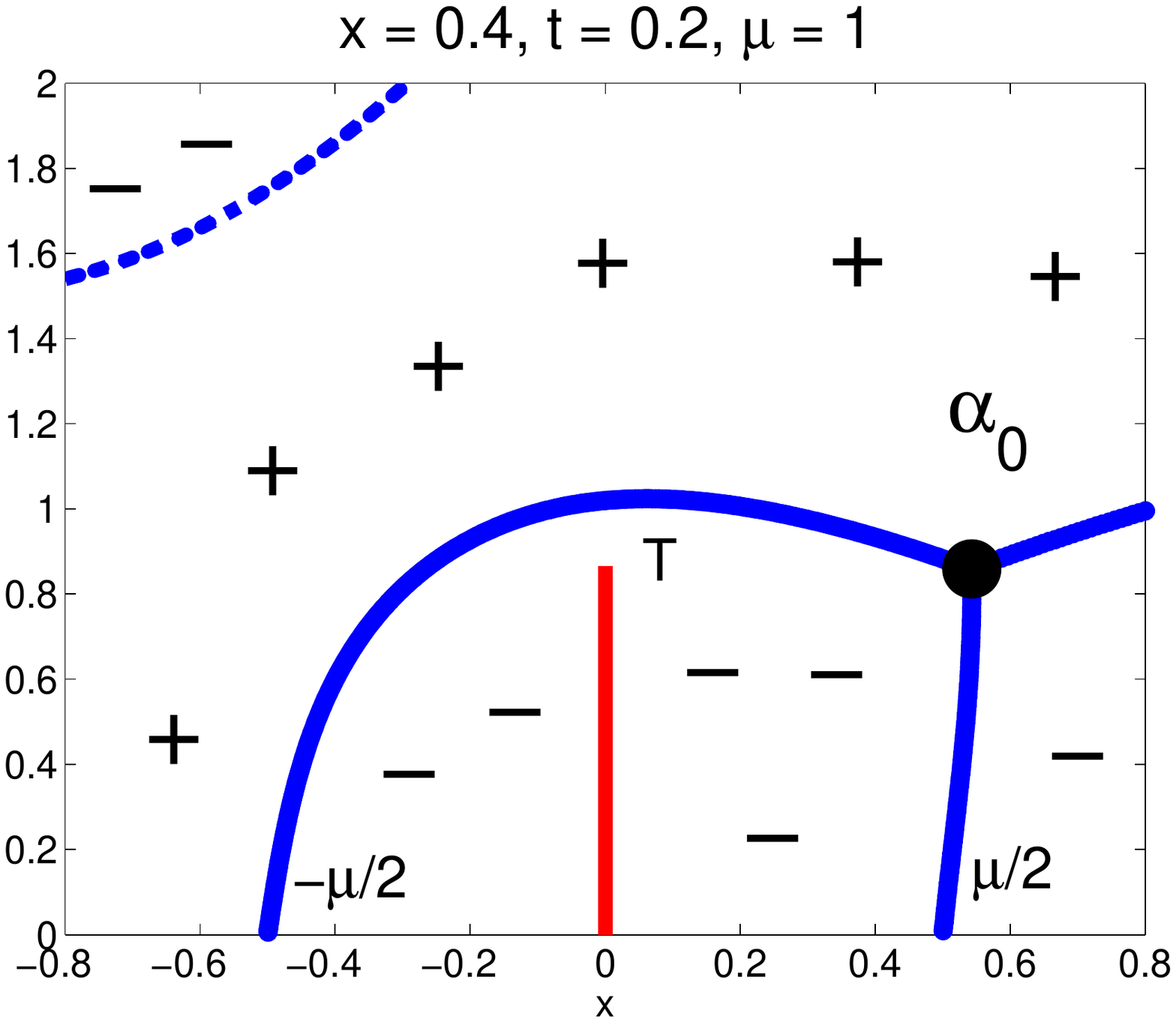}
\includegraphics[height=7cm]{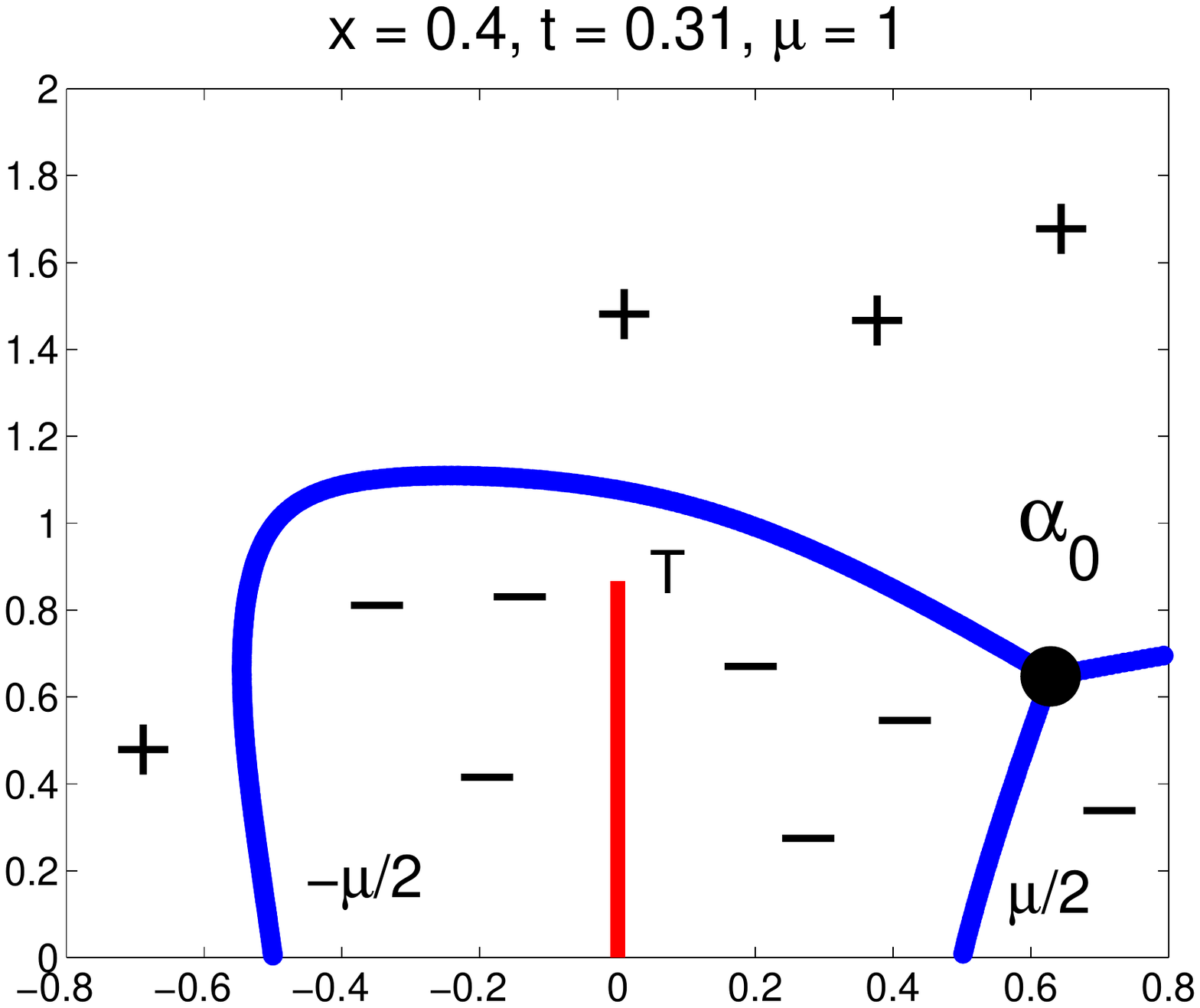}
\includegraphics[height=7cm]{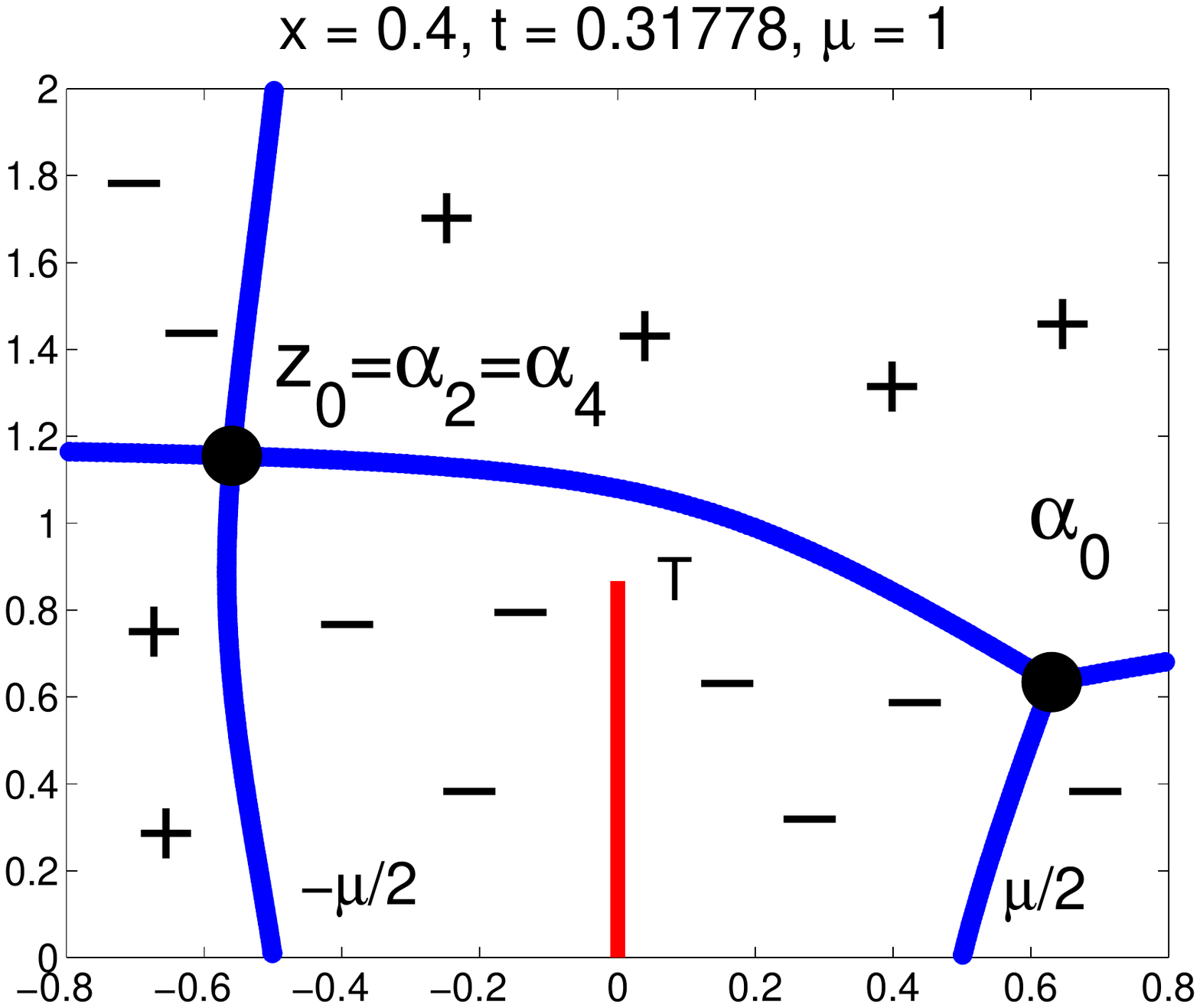}\\
\includegraphics[height=7cm]{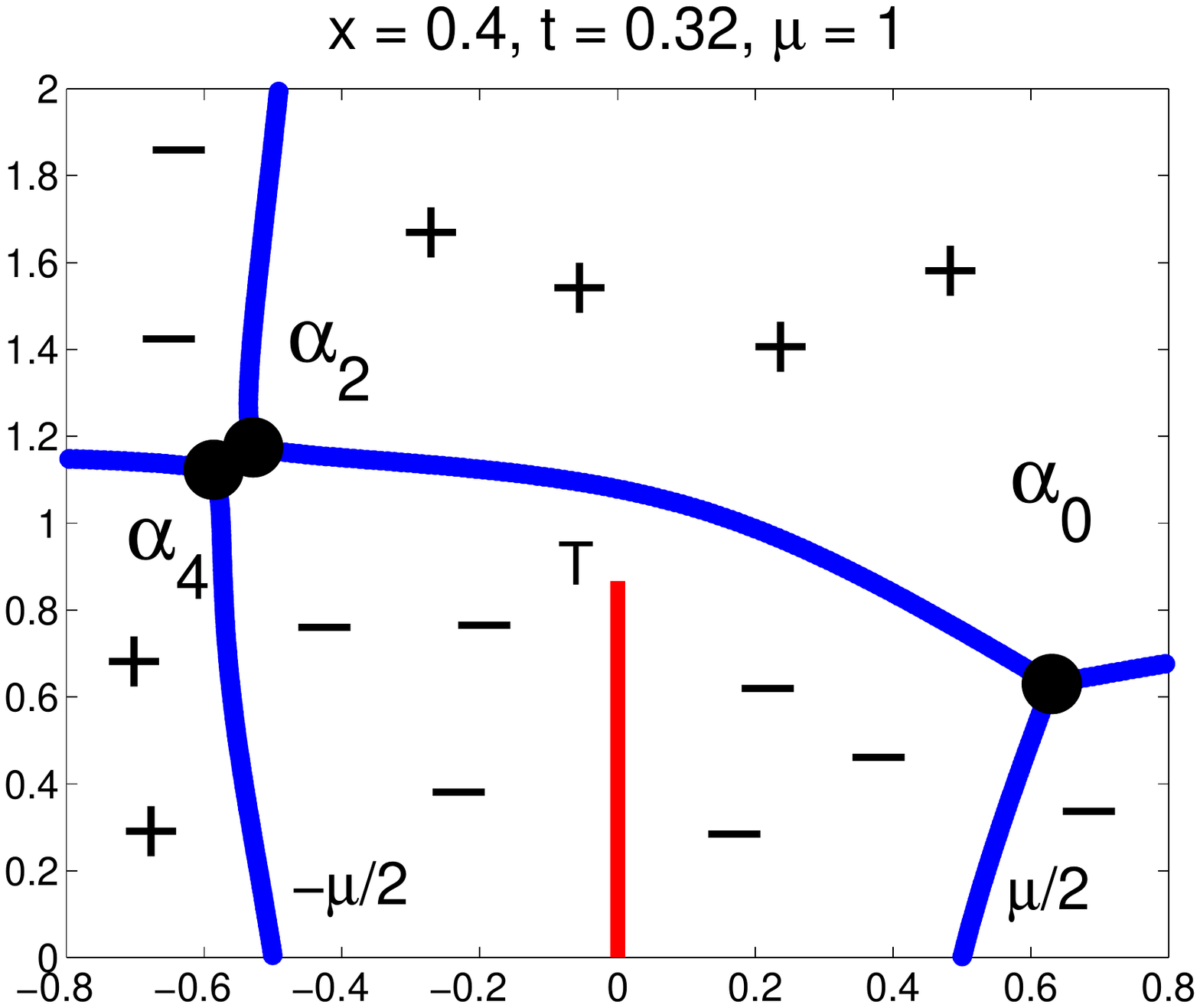}
\includegraphics[height=7cm]{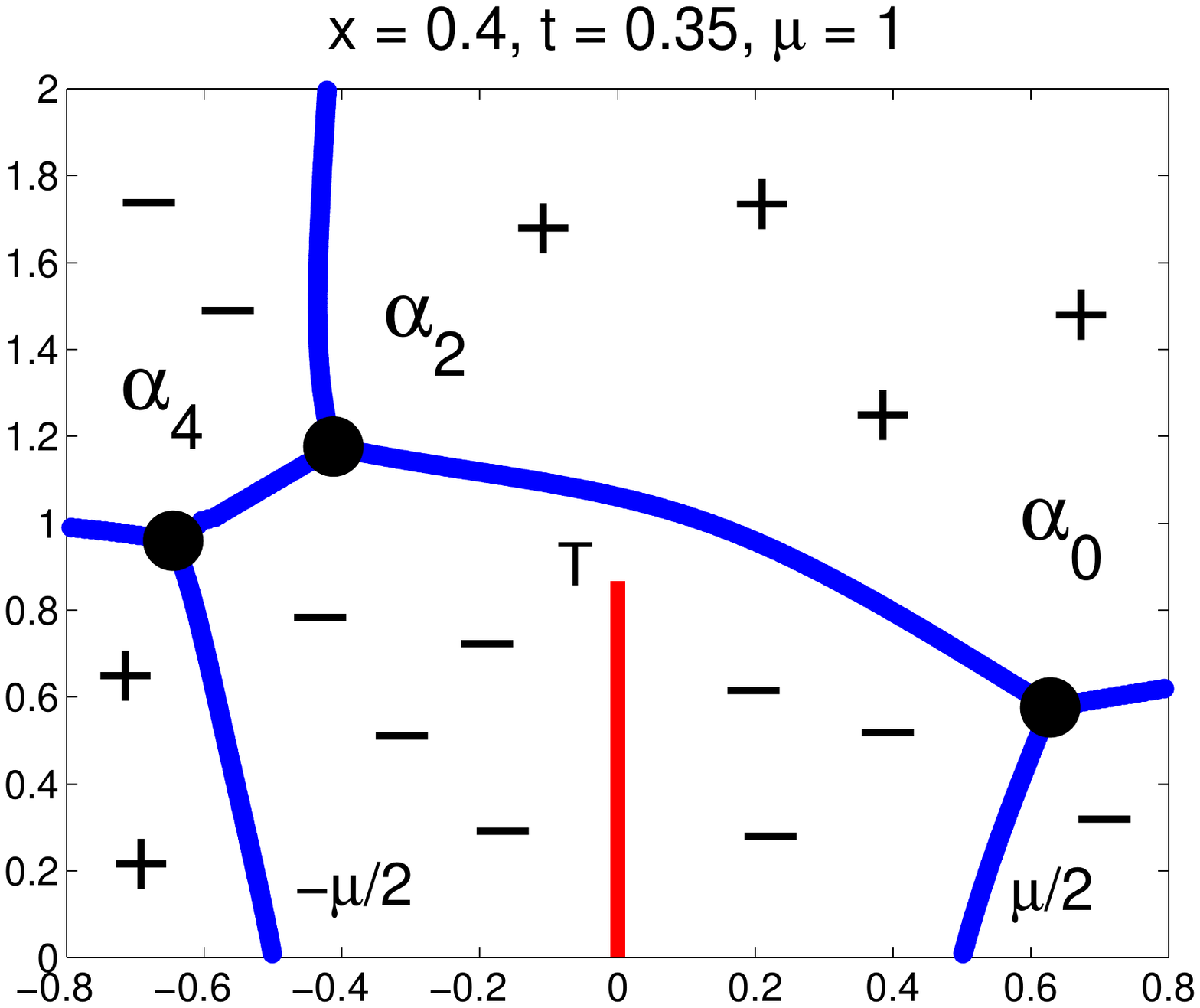}
\includegraphics[height=7cm]{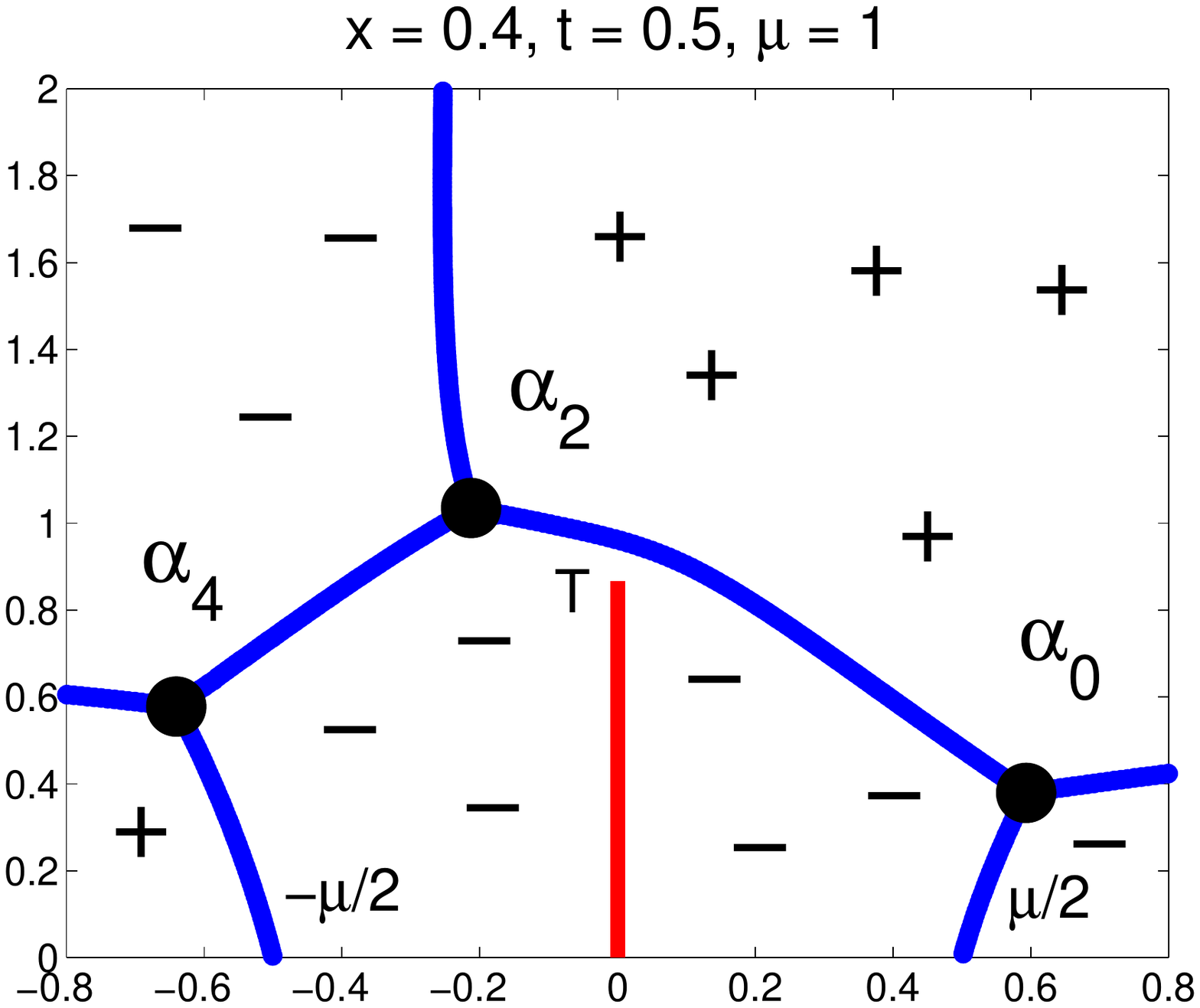}
\caption{\label{fig4:br1_mechanism} Mechanism of the first break:
$t$ evolution of the zero level curves of $\Im h$ for $x=0.4$,
$\mu=1$. This sequence illustrates the birth of the main arc
$[\alpha_2,\alpha_4]$.}
\end{center}
\end{figure}

The first break of the asymptotic solution of NLS (\ref{1.1:NLS}) with the initial conditions
(\ref{1.1:IC_original}) in the semiclassical limit $\eps\to 0$
was analytically studied by Tovbis, Venakides, and
Zhou in \cite{TVZzero_04}. They established the mechanism of the
first break for the $\mu>0$ analytically as time $t$ evolution
process. From $t=0$, as $t$ increases the genus changes from $0$ to $2$ with a new main arc
$[\alpha_2,\alpha_4]$ created in the upper half plane. In the present
work we support their proof numerically.

Figure \ref{fig4:br1_mechanism} illustrates the mechanism of the
first break from the point of view of the branch points $\alpha$'s and zero level curves
of $\Im h$. First, for small $t$ the genus is 0 and there is only one branch
point in the upper half plane $\alpha_0$. Then, a new pair of branch points $[\alpha_2,\alpha_4]$
is created, while $\alpha_0$ continues to approach $\mm$ under a modified
trajectory.

\begin{figure}
\begin{center}
\includegraphics[height=10cm]{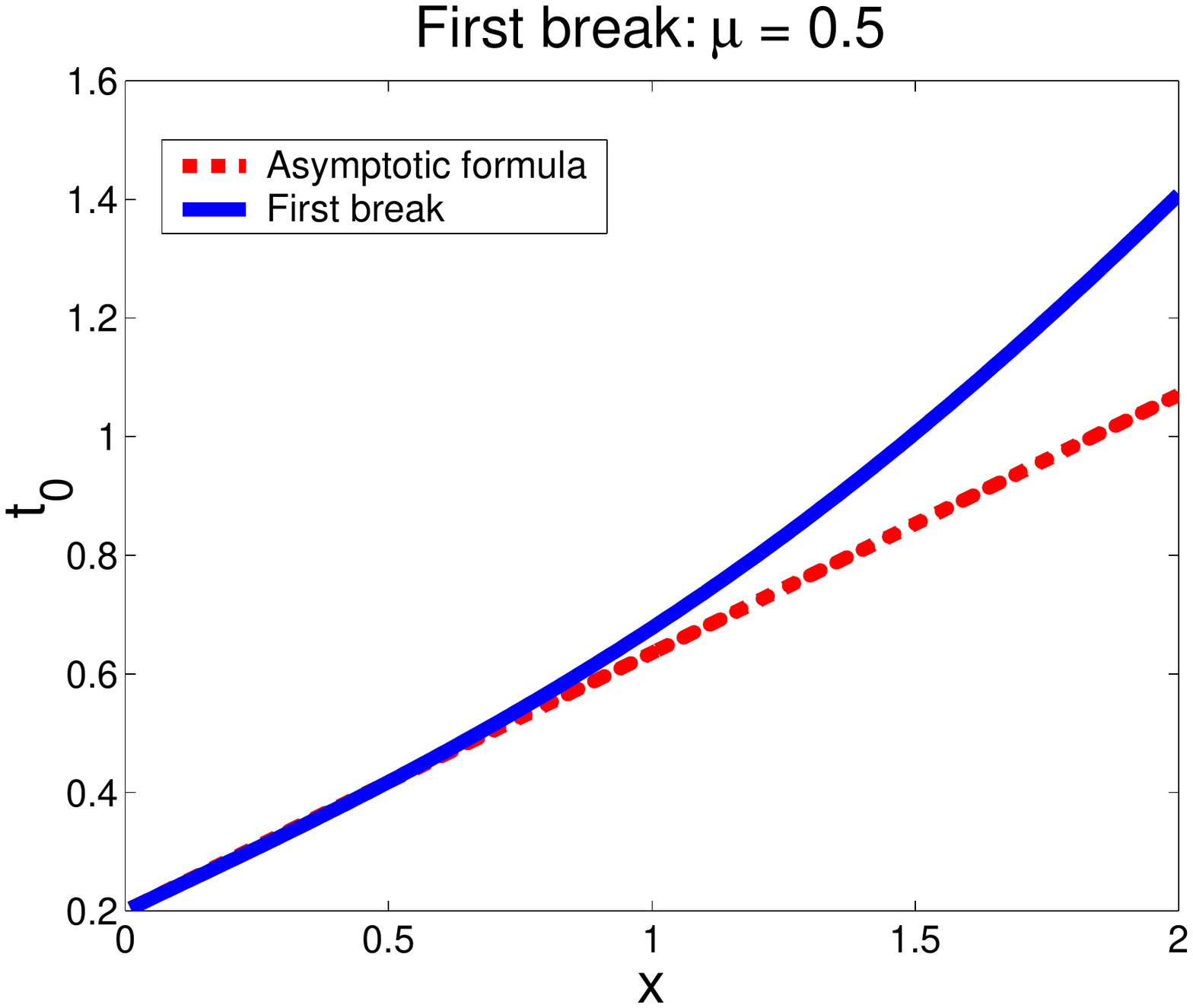}
\caption{\label{fig4:break1_asympt} First breaking curve $t_0(x)$
(solid) and its small $x$ asymptotics (\ref{eq4:br1_as})
(dashed).}
\end{center}
\end{figure}

The asymptotic behavior of the first breaking curve in the large
and small $x$ limits for $\mu\ge2$ was established in
\cite{TVZzero_04} and \cite{TVZlong2_06}
\begin{equation}
t_0(x)=\left\{\begin{array}{l}
\frac{1}{2(\mu+2)}+\frac{\cot\frac{\pi}{5}}{2\sqrt{\mu+2}}x+O\left(x^{3/2}\right),\quad
x\to 0 \\
\frac{x}{2\mu}-\frac{1}{\mu}\ln\frac{2\mu}{\mu+2T}-\frac{T/\mu}{\mu+2T}+O\left(\frac{1}{x}\right),\quad
x\to \infty
\end{array}\right. .\label{eq4:br1_as}
\end{equation}
In the soliton+radiation case $0<\mu<2$ the above large $x$
asymptotic expression produces complex answers. Our
conjecture is that the above expression is correct if one substitutes
$T=0$. This conjecture is based on the comparison of our
our leading order term in the long time limit of $\alpha_2$
(\ref{thma2b2}) with the leading term in the expression above (1.5)
in \cite{TVZlong2_06} for $\mu\ge 2$.

The small $x$ asymptotics in (\ref{eq4:br1_as}) is proved to be
valid \cite{TVZzero_04} for $0<\mu<2$. Figure
\ref{fig4:break1_asympt} demonstrates agreement of the first
breaking curve with the small $x$ asymptotic formula.

In this case, the first break is a boundary between genus 0 and
genus 2 regions. From the point of view of genus 2, the first
break is a singular event of colliding of two branch points
$\alpha_2$ and $\alpha_4$ with the main arc $[\alpha_2,\alpha_4]$
reducing to a point. Numerically we observed this phenomenon in $x$,
$t$ and $\mu$ evolutions. Figure \ref{fig4:alphas_g2_evolution} demonstrates
how the choice of parameters $x$ and $t$ for the branch points
$\alpha_0$, $\alpha_2$, $\alpha_4$ evolution correspond to the
first breaking curves.

\begin{figure}
\begin{center}
\includegraphics[height=10cm]{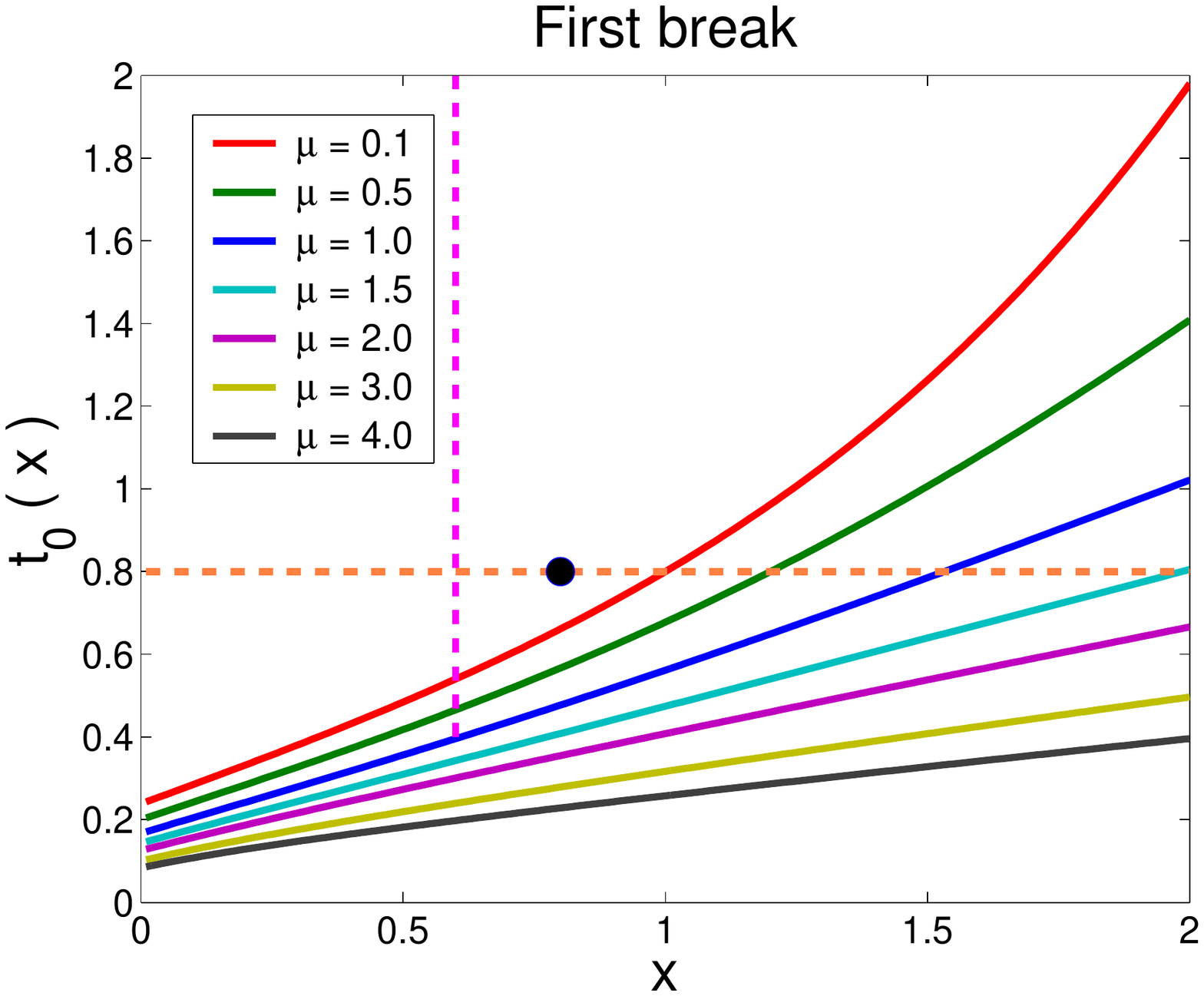}
\includegraphics[height=10cm]{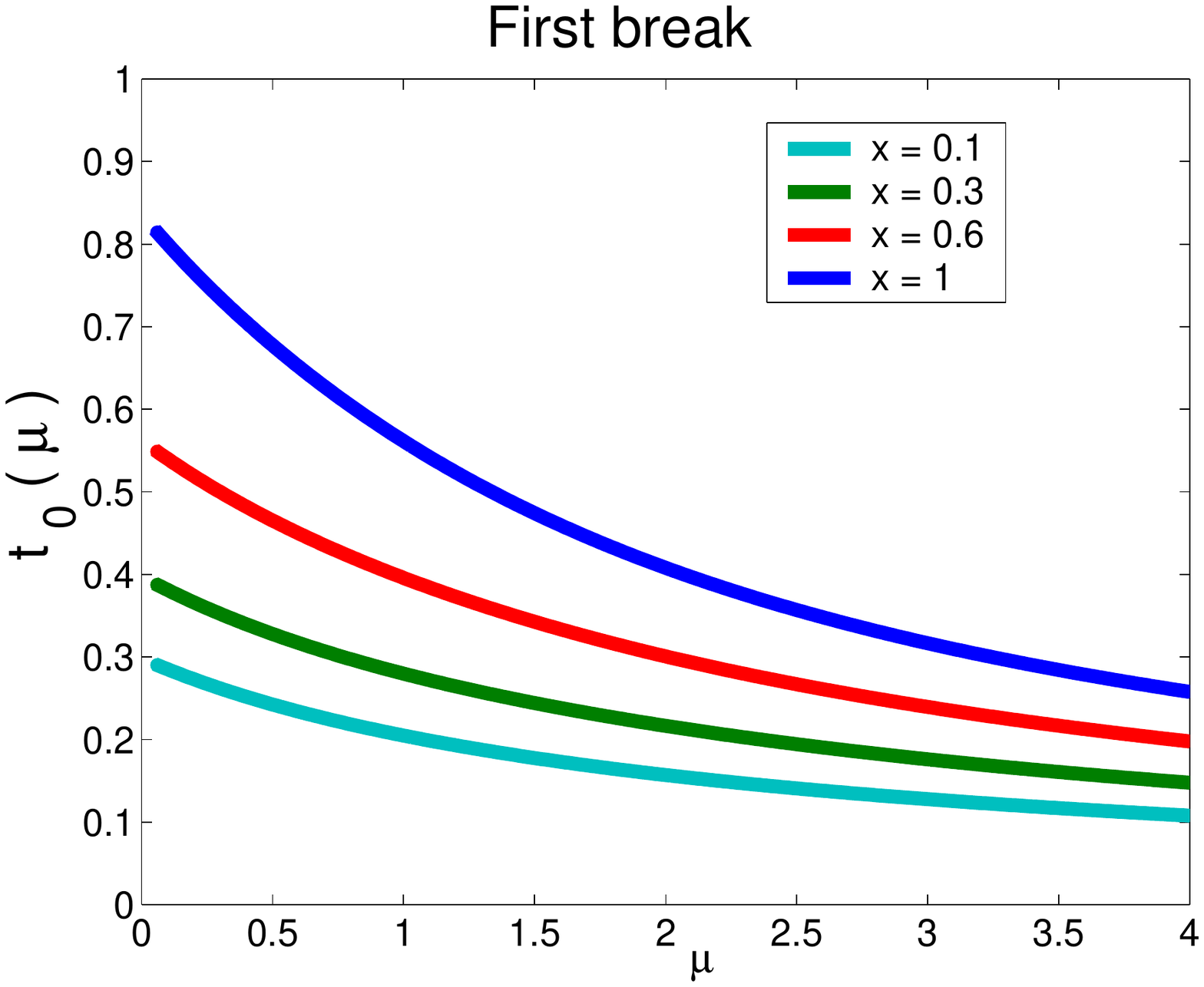}
\caption{\label{fig4:breaks1} First breaking curves: $t_0(x)$ for
several $\mu$ values (left), $t_0(\mu)$ for several $x$ values
(right).}
\end{center}
\end{figure}
\begin{figure}
\begin{center}
\includegraphics[height=7cm]{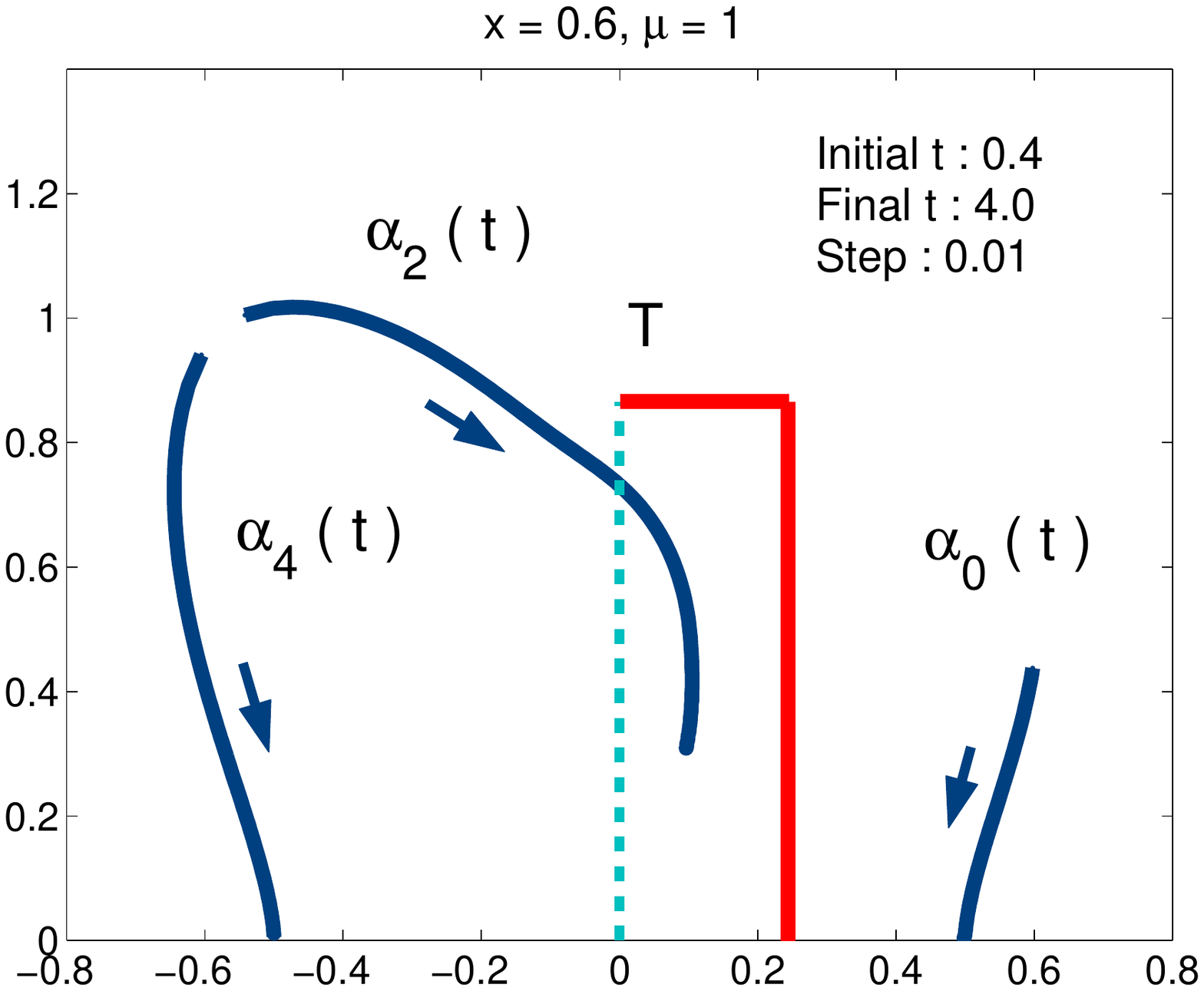}
\includegraphics[height=7cm]{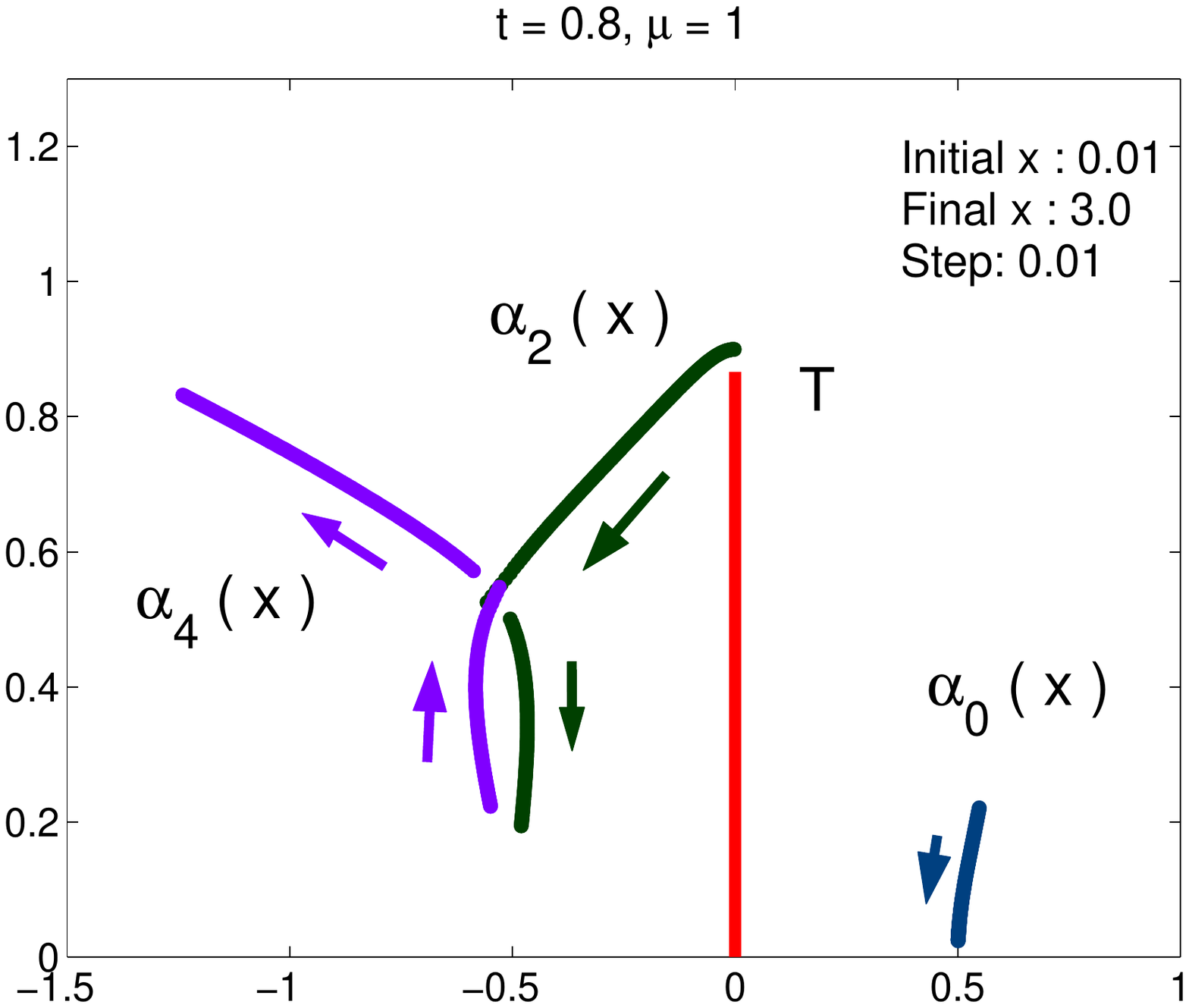}
\includegraphics[height=7cm]{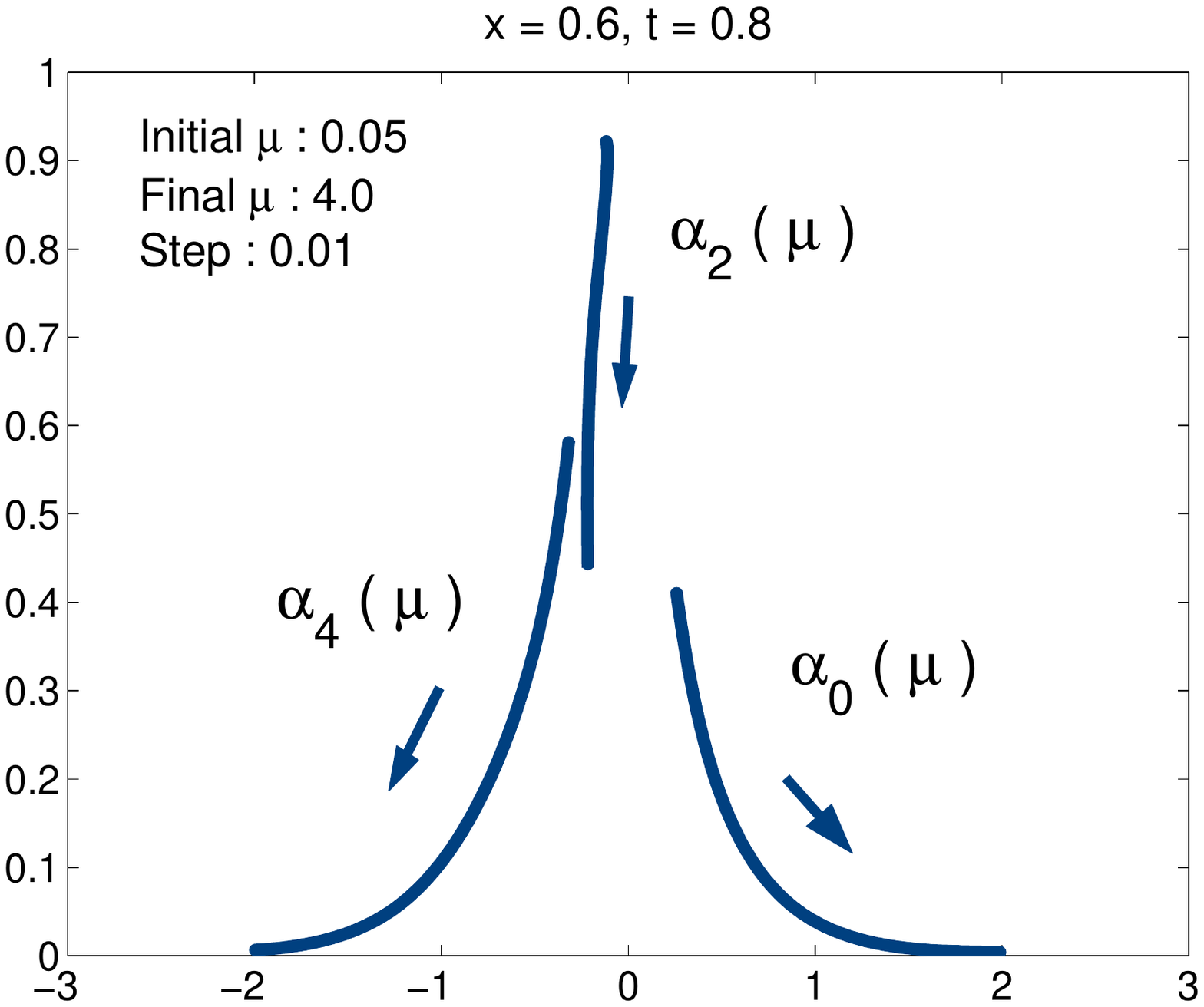}
\caption{\label{fig4:alphas_g2_evolution} Time $t$ (left), space $x$
(middle), and parameter $\mu$ (right) evolution of the branch points
$(\alpha_0,\alpha_2,\alpha_4)$ in genus 2. Collision of $\alpha_2$
and $\alpha_4$ corresponds to the first break.}
\end{center}
\end{figure}

The vertical dashed line on Figure \ref{fig4:breaks1} (left picture) at $x=0.6$ represents the $t$ evolution
for $\mu=1$ in Fig. \ref{fig4:alphas_g2_evolution} (left picture).
The line starts above the first breaking curve for $\mu=1$
so the branch points $\alpha_2$ and $\alpha_4$ do not collide.

The horizontal dashed line on Figure \ref{fig4:breaks1} (left picture) $t=0.8$ shows the $x$ evolution for
$\mu=1$ in Fig. \ref{fig4:alphas_g2_evolution} (middle picture).
This line intersects the first breaking curve at $x\approx
1.5$.

The big black dot at $x=0.8$, $t=0.8$ represents the
parameters for the $\mu$ evolution in Figure
\ref{fig4:alphas_g2_evolution} (right picture). The dot is located
in genus 2 region above all the breaking curves which is confirmed by
$\alpha$'s trajectories without collisions for $\mu$ values
between $0.05$ and $4.0$.

\begin{figure}
\begin{center}
\includegraphics[height=10cm]{br1_x.pdf}
\includegraphics[height=10cm]{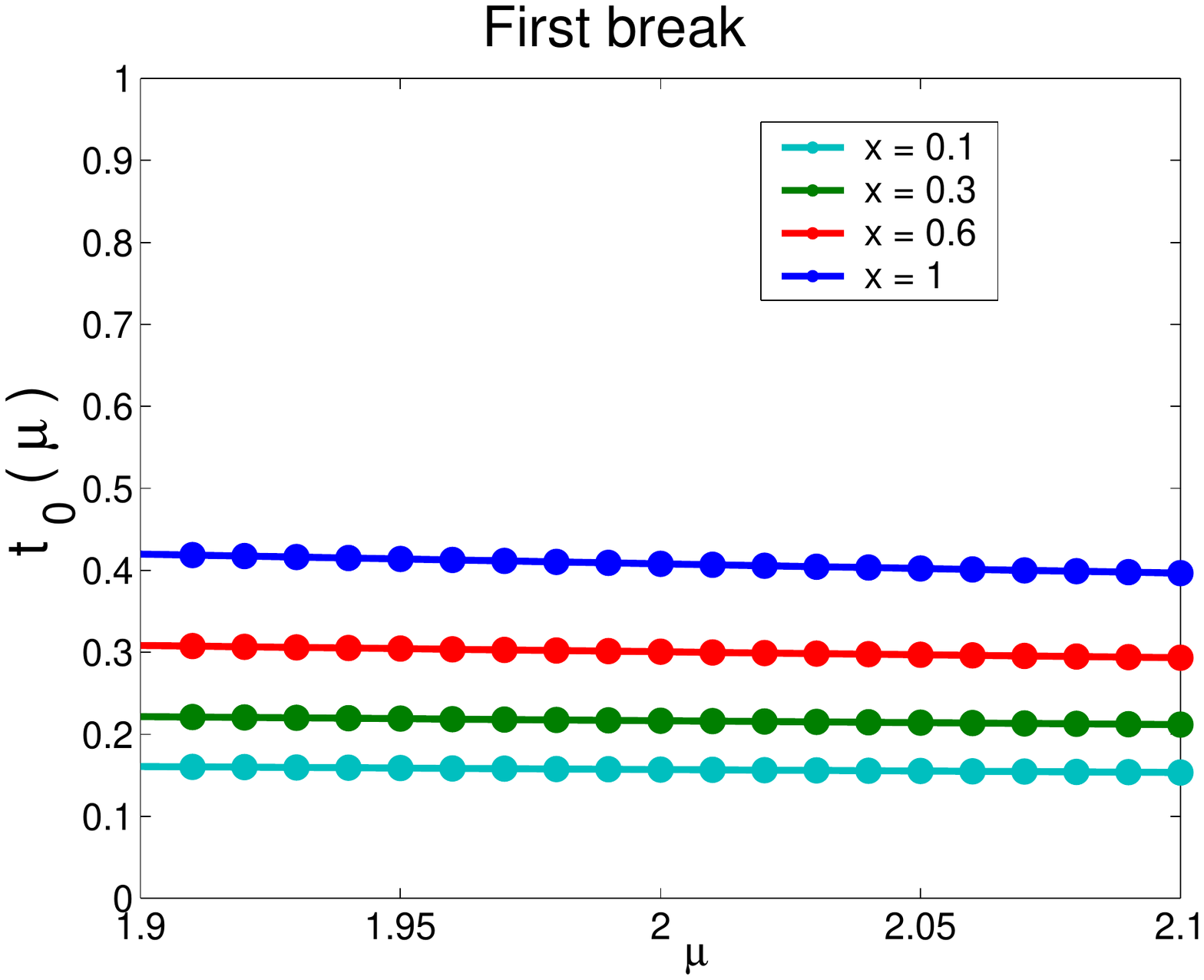}
\caption{\label{fig4:break1_x} First breaking curves $t=t_0(\mu)$
for several $x$ values.}
\end{center}
\end{figure}

Finally, we look at the first break as a function of the parameter $\mu$.
Figure \ref{fig4:break1_x} shows no sign of loss of smoothness at
$\mu=2$ which is a critical value for existence of solitons in the initial
conditions (\ref{1.1:IC_original}). This dependence is investigated in more details
in \cite{BelovVen_pre}.

\subsection{Singular obstruction}

\begin{figure}
\begin{center}
\includegraphics[height=7cm]{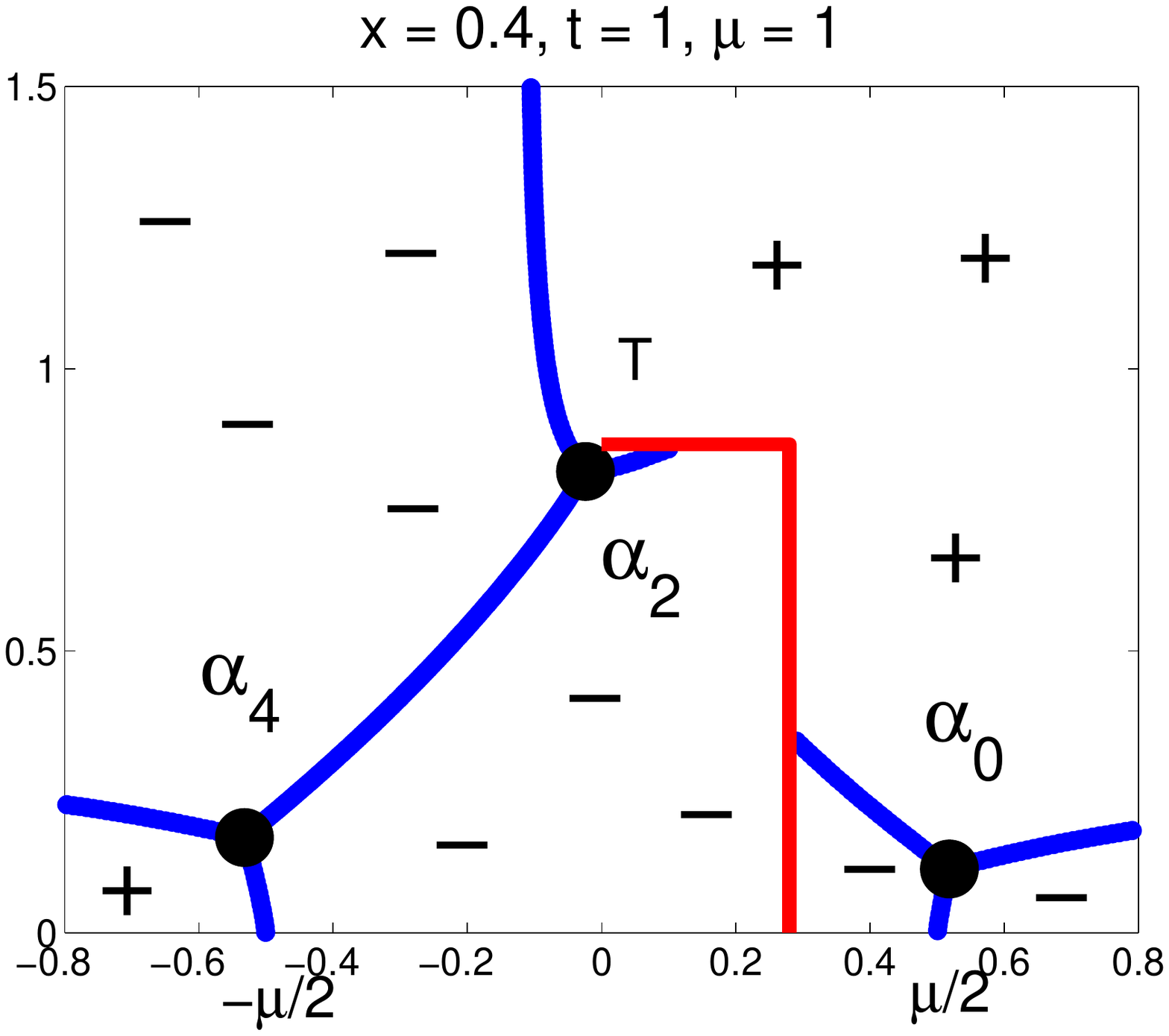}
\includegraphics[height=7cm]{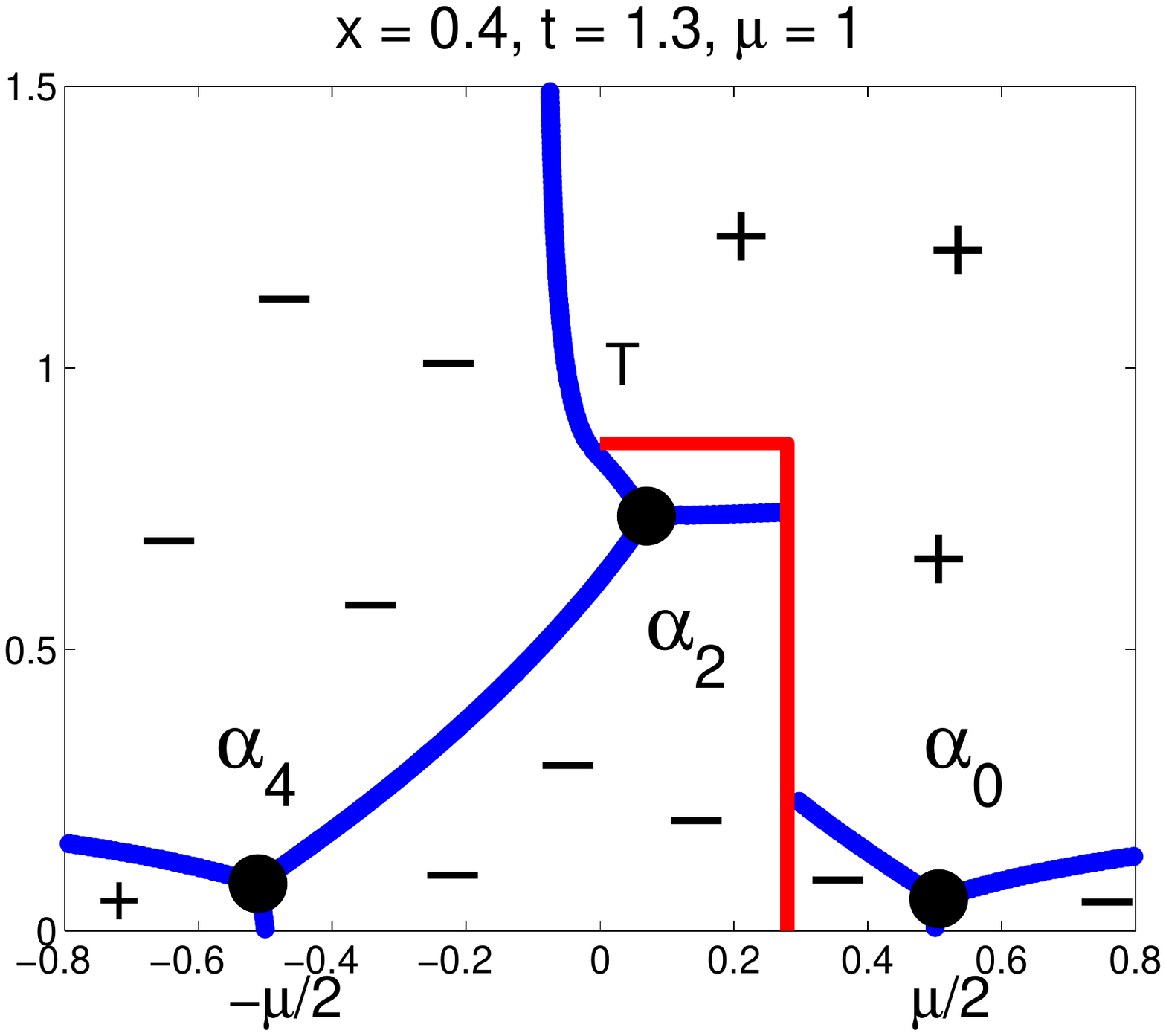}
\includegraphics[height=7cm]{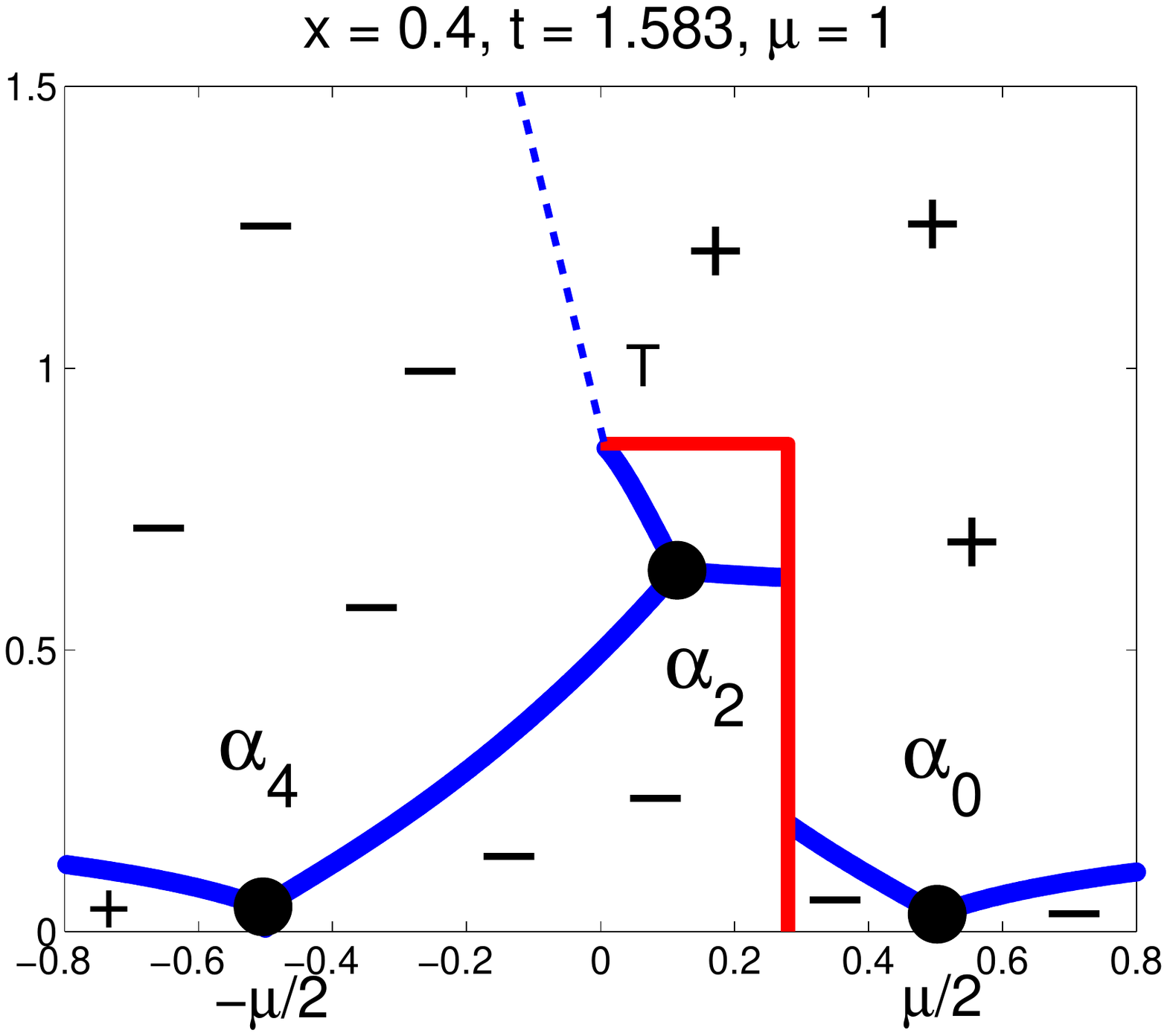}
\caption{\label{fig4:brc_mechanism} Mechanism of the singular
obstruction: time evolution of the zero level curves of $\Im h$ for
$x=0.4$, $\mu=1$ in genus 2.}
\end{center}
\end{figure}

The mechanism for the singular
obstruction is a collision of a branch of zero level curves with the
logarithmic branch point $T$ (see Fig. \ref{fig4:brc_mechanism}).
This collision closes the passage between
$\alpha_0$ and $\alpha_2$ around the logarithmic branch cut
$[0,T]$ and invalidates error estimates. Formally all the
expressions $h'(z)$, $h(z)$ are correct as solutions of RH problems however the underlying
assumptions are not valid: in Figure \ref{fig4:brc_mechanism}
(right picture) there is no path to connect $\mm$ with $-\infty$
satisfying $\Im h \ge0$ which is necessary to guarantee
capturing the leading order of the asymptotic solution of NLS.

It is not clear at this point how to extend the function $h(z)$
beyond the singular obstruction. However, the term "singular
obstruction" symbolizes the difficulties of the asymptotic analysis
rather than drastic changes in the solution. The main difficulty
is to extend a Riemann surface continuously after a
collision of a zero level curve with a logarithmic
branch point singularity.

In this paper we confirm numerically existence of the
singular obstruction and compute
two correction terms in the long time limit.
\begin{equation}
x=\ln 2+c_2\frac{\ln t}{t}+c_3\frac{1}{t}+ O\left(\frac{\ln
t}{t^{3/2}}\right), \ \ \ \ t\to\infty, \label{eq4:br2}
\end{equation}
where
\begin{equation}
c_2=-\frac{\mu}{16|T|^2},\quad \quad
c_3=-\frac{\mu\left(1+\ln\frac{4|T|^2}{\mu} \right)}{16|T|^2},
\end{equation}

\begin{figure}
\begin{center}
\includegraphics[height=10cm]{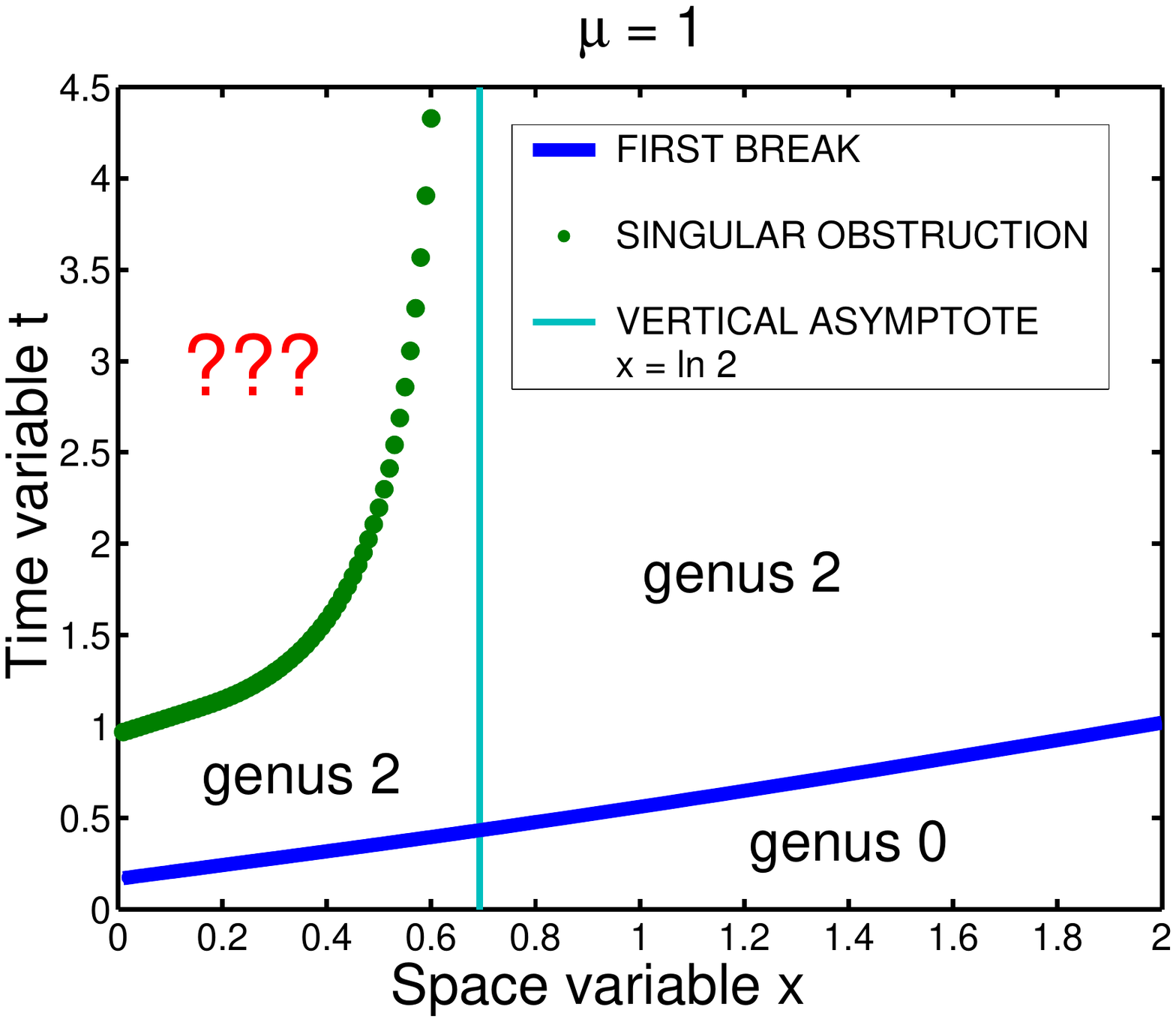}
\includegraphics[height=10cm]{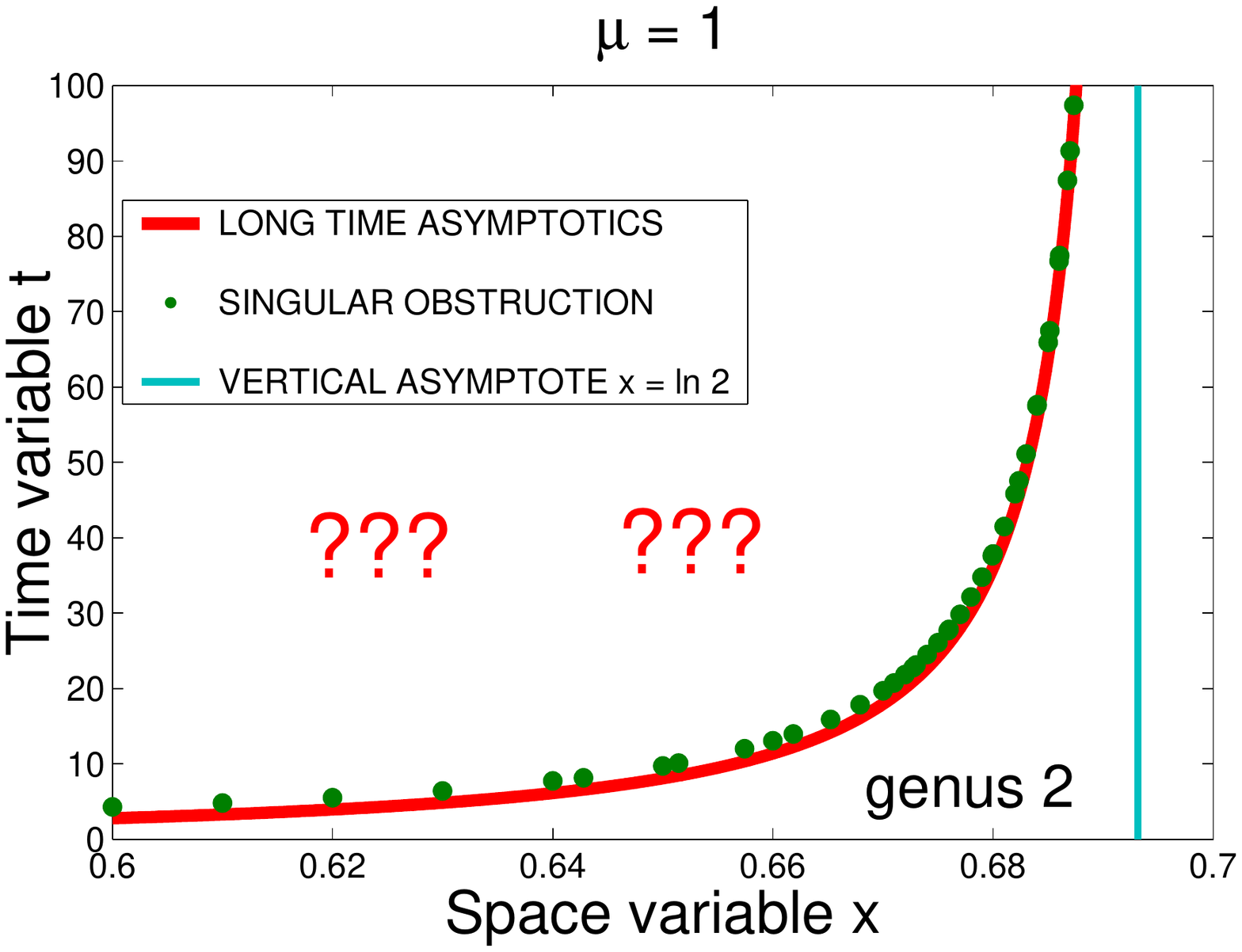}
\caption{\label{fig4:cat_break_as} Singular obstruction curve
$t=t_s(x)$ and first breaking curve $t_0(x)$ for $\mu=1$ are on the
left. The long time asymptotics of the singular obstruction curve
(\ref{eq4:br2}) and the long time computations of the singular
obstruction curve $t=t_s^{LT}(x)$ are on the right.}
\end{center}
\end{figure}
Comparison of the long time computations of the singular obstruction
with the asymptotic formula (\ref{thmbr2}) is given in Figure
\ref{fig4:cat_break_as} (right picture).

Based on numerical evidence, our conjecture that the singular
obstruction exists only for a finite interval $[0, \ln 2)$ and the
location of the vertical asymptote is independent of $\mu$.

From the asymptotic analysis of $\Im h(T,t)$ as a function of $t$ in (\ref{Im_h_at_T_asympt}),
\begin{equation}
\Im h(T,t)=|T|(x-\ln 2)+o(1), \quad t\to\infty.
\end{equation}
So $\Im h(T,t)$ has a horizontal asymptote $|T|(x-\ln 2)$ as $t\to\infty$.

Figure \ref{fig3:imag_h_at_T} suggests that for $x>\ln 2$
function $\Im h(T,t)$ has only one root while the other root
corresponding to the singular obstruction is not present.
This conjecture is supported by the asymptotics of $\Im
h(T,t)$ which is asymptotically strictly positive for $x>\ln 2$ for
all $t$ large enough.

There is no numerical evidence of any other breaks to occur before
the singular obstruction. It is an open question how to extend the current
asymptotics beyond the singular obstruction curve. More constructively,
this is a question of the new genus of the Riemann surface after collision
of the zero level curve with a logarithmic branch point.

Calculations by Lyng, Miller \cite{Lyng_07} in the case $\mu=0$
suggest that after the second break the genus is 4. However, they
are using a somewhat different approach of changing the RHP by
adjusting the reflection coefficient. Effectively, they are
changing approximations of $f(z)$ on the fly: for some region in the $(x,t)$ plane
they are using the same $f(z)$ as we are, while for other regions they are
considering a different approximation of $f(z)$ (different approximation of the reflection
coefficient at later $t$ values) and consequently a
different RHP. This seems to be equivalent to
dealing with two $f$-functions at the same time or dealing with
two sheets of the Riemann surface at the same time.

Our analysis is based on approximating $f(z)$ at $t=0$ upto order $O(\eps)$
and use this approximation for all $t>0$.
It is possible that the leading terms in the used approximation of
$f(z)$ as $\eps\to 0$, do not contain information about the second
break. Thus taking into account correction terms of the order
$O(\eps)$ and $O(\eps^2)$ is another approach to describing the
second break.

\section{Appendix}

\subsection{Higher order terms of $a_2$, $b_2$ in Theorem \ref{thm_1}}


In this section we compute higher order terms in the asymptotics of $\alpha_2=a_2+ib_2$
in the long time limit from the couple of moment conditions (\ref{M34as}).

\subsubsection{Simplification of the first moment condition (\ref{M34as})}

Consider the first moment condition with the exponentially small terms dropped
\begin{equation}
\int_{a}^{-\frac{\mu}{2}}
\frac{1}{\sqrt{\left(\xi-a\right)^2+b^2}}d\xi +\frac{1}{\pi
i}\int_{a+ib}^{a-ib}
\frac{f'(\xi)}{\sqrt{\left(\xi-a\right)^2+b^2}}d\xi
=O\left(\frac{1}{t^2}\right).\\
\label{M3as}
\end{equation}
As it was shown above, $a=O\left(\frac{\ln t}{t}\right)$ and $b=O\left(\frac{1}{\sqrt{t}}\right)$
 as $t\to\infty$. Our goal is to compute the next
two terms of asymptotics of $a$ and $b$. First, we compute the first integral
explicitly and make a change of variables in the second integral $\xi=by+a$
\begin{equation}
\left.
\ln\left(\xi-a+\sqrt{(\xi-a)^2+b^2}\right)\right|_a^{-\frac{\mu}{2}}
+\frac{1}{\pi i}\int_{i}^{-i} \frac{f'(by+a)}{\sqrt{y^2+1}}dy=O\left(\frac{1}{t^5}\right)
\label{a2_1}.
\end{equation}
Since both $a$ and $b$ are small, to estimate the
integral we decompose $f'(\xi)$ in Taylor series near $\xi=0$
while isolating the term containing $t$ explicitly
\begin{equation}
f'(\xi)=-4t\xi+f'(0)+(f''(0)+4t)\xi+\frac{f'''(0)}{2}\xi^2+O(\xi^3). \label{df_as}
\end{equation}
Next, we incorporate Schwartz reflection symmetry
$f'\left(\overline\xi\right)=\overline{f'(\xi)}$ and the fact that
all $f^{(k)}(0)$ are real for $k\ge2$, then
\begin{equation}
f'(by+a)=-4t(by+a)+C_1+iC_2\mbox{sign}\left(\frac{y}{i}\right)+C_3(by+a)+C_4(by+a)^2+O((by+a)^3),
\label{df_series}
\end{equation}
where the constants $C_j$ are obtained from (\ref{df}) as limiting
values
\begin{equation}
f^{(j)}(0)=\lim_{\varepsilon\to0^+}(f^{(j)}(0+i\varepsilon)).
\end{equation}
So
\begin{equation}
\left\{\begin{array}{l}
C_1=\ln \frac{2|T|}{\mu}-x\\
C_2=\frac{\pi}{2}\\
C_3=\frac{2}{\mu}
\end{array}\right.\label{C_values}.
\end{equation}
After plugging (\ref{df_series}) and (\ref{C_values}) back in
(\ref{a2_1}) one obtains
\[
\ln\left(\sqrt{\left(\frac{\mu}{2}+a\right)^2+b^2}-\left(\frac{\mu}{2}+a\right)\right)
-\ln b +\frac{1}{\pi i}\int_{i}^{-i}
\frac{-4t(by+a)}{\sqrt{y^2+1}}dy
\]
\begin{equation}
 +\frac{1}{\pi i}\int_{i}^{-i}
\frac{C_1+iC_2\mbox{sign}\left(\frac{y}{i}\right)+C_3(by+a)+C_4(by+a)^2+O((by+a)^3)}{\sqrt{y^2+1}}dy=0,
\label{a2_2}
\end{equation}
where after expanding all the terms, using symmetries in
the integrals, and only keeping terms upto order $O(b^2)=O\left(\frac{1}{t}\right)$, we see
\begin{equation}
\ln\left(\frac{b}{2(\frac{\mu}{2}+a)}+O(b^3)\right) +\frac{1}{\pi
i}\int_{i}^{-i} \frac{-4ta+C_1+ C_3a}{\sqrt{y^2+1}}dy +O(b^2)=0,
\label{a2_4}
\end{equation}
which simplifies to
\begin{equation}
\ln\left(\frac{b}{2(\frac{\mu}{2}+a)}\right) - (-4ta+C_1+
C_3a)+O(b^2)=0 \label{a2_5}
\end{equation}
and its final form is
\begin{equation}
\ln b+ 4ta - C_1-\ln\mu+ \left(-\frac{2}{\mu}-C_3\right)a +
O(b^2)=0 \label{a2_7}.
\end{equation}
Finally, we plug in values of $C_1$ and $C_3$ from (\ref{C_values})
\begin{equation}
\ln b+ 4ta +x-\ln 2|T| -\frac{4a}{\mu} + O(b^2)=0 \label{a2_8}.
\end{equation}
In order to extract the asymptotics of $a$ and $b$ from this
equation we need to couple it with the other moment condition.


\subsubsection{Simplification of the second moment condition (\ref{M34as})}

Consider
\begin{equation}
\int_{a}^{-\frac{\mu}{2}}
\frac{\xi-a}{\sqrt{\left(\xi-a\right)^2+b^2}}d\xi +\frac{1}{\pi
i}\int_{a+ib}^{a-ib}
\frac{(\xi-a)f'(\xi)}{\sqrt{\left(\xi-a\right)^2+b^2}}d\xi =0.
\label{M4as}
\end{equation}
In a similar manner as for the other moment condition, we evaluate
the first integral and make a change of variables $\xi=by+a$ in
the second integral
\begin{equation}
\sqrt{\left(\frac{\mu}{2}+a\right)^2+b^2}-b
 +\frac{1}{\pi
i}\int_{i}^{-i} \frac{byf'(by+a)}{\sqrt{y^2+1}}dy =0.\label{b2_2}
\end{equation}
Next, we compute the integral by expanding $f'(\xi)$ into Taylor series
(\ref{df_series})
\begin{equation}
\int_{i}^{-i} \frac{byf'(by+a)}{\sqrt{y^2+1}}dy
\end{equation}
\begin{equation}
=\int_{i}^{-i} \frac{by\left[
-4t(by+a)+C_1+iC_2\mbox{sign}\left(\frac{y}{i}\right)+C_3(by+a)
+C_4(by+a)^2+O((by+a)^3)\right]}{\sqrt{y^2+1}}dy.
\end{equation}
and using the symmetry in the integral which make many of the terms to disappear
\begin{equation}
=\left(-4tb^2+C_3b^2\right)\int_{i}^{-i} \frac{ y^2
}{\sqrt{y^2+1}}dy+2C_2bi\int_{i}^{0} \frac{ y
}{\sqrt{y^2+1}}dy+O(ab^2)
\end{equation}
\begin{equation}
=\left(-4tb^2+C_3b^2\right)\frac{\pi i}{2}+2C_2bi
+O\left(ab^2\right).
\end{equation}
Then the equation (\ref{b2_2}) reads
\begin{equation}
\sqrt{\left(\frac{\mu}{2}+a\right)^2+b^2}-b
 -2tb^2+\frac{C_3b^2}{2}+\frac{2C_2b}{\pi
} +O\left(ab^2\right)=0\label{b2_3}
\end{equation}
and after expanding the square root
\begin{equation}
\frac{\mu}{2}-2tb^2+\left(\frac{2C_2}{\pi
}-1\right)b+a+\left(\frac{1}{\mu}
 +\frac{C_3}{2}\right)b^2 +O\left(ab^2\right)=0.\label{b2_5}
\end{equation}
This equation is an extended version of the second equation in
(\ref{a2b2leading_2}). After substituting expressions for $C_2$
and $C_3$ from (\ref{C_values}) we see that the term of order the
$O(b)$ vanishes
\begin{equation}
\frac{\mu}{2}-2tb^2+a+\frac{2b^2}{\mu}
 +O\left(ab^2\right)=0.\label{b2_6}
\end{equation}

\subsubsection{Asymptotically solving the system of moment conditions (\ref{M34as})}

We want to solve the system of two asymptotic equations for $a$
and $b$
\begin{equation}
\left\{\begin{array}{l}
\ln b+ 4ta +x-\ln 2|T| -\frac{4a}{\mu} + O(b^2)=0 \\
\frac{\mu}{2}-2tb^2+a+\frac{2b^2}{\mu} +O\left(ab^2\right)=0,
\end{array}\right.\label{a2b2_1}
\end{equation}
where as we established in (\ref{eq1.4:a2b2_lead}) the leading
order solutions are
\begin{equation}
\left\{
\begin{array}{l}
a=\frac{\ln t}{8t}+O\left(\frac{1}{t}\right)\\
b=\sqrt{\frac{\mu}{4t}}+O\left(\frac{1}{t}\right).
\end{array}\right.\label{a2b2as_1}
\end{equation}
We find the correction terms by writing
\begin{equation}
\left\{
\begin{array}{l}
a=\frac{\ln t}{8t}+\frac{\varepsilon_1 \ln t}{t}\\
b=\sqrt{\frac{\mu}{4t}}+\frac{\delta_1}{\sqrt t}
\end{array}\right.\label{a2b2as_2}
\end{equation}
with functions $\varepsilon_1$, $\delta_1$ to be determined. Then
\begin{equation}
\left.
\begin{array}{l}
\ln b=-\frac{1}{2}\ln t+\ln\sqrt{\frac{\mu}{4}}+O(\delta_1)\\
b^2=\frac{\mu}{4t}+\frac{\sqrt\mu\delta_1}{t}+\frac{\delta_1^2}{t},
\end{array}\right.\label{a2b2_2}
\end{equation}
which we substitute into (\ref{a2b2_1}) and obtain
\begin{equation}
\left\{\begin{array}{l}
-\frac{1}{2}\ln t+\ln\sqrt{\frac{\mu}{4}} +O(\delta_1)+ \frac{1}{2}\ln t+4\varepsilon_1 \ln t +x-\ln 2|T| +O(a)=0 \\
\frac{\mu}{2}-\frac{\mu}{2}-2\sqrt\mu\delta_1-2\delta_1^2
+\frac{\ln t}{8t}+\frac{\varepsilon_1 \ln
t}{t}+O\left(b^2\right)=0.
\end{array}\right.\label{a2b2_3}
\end{equation}
This leads to the following asymptotic equations
\begin{equation}
\left\{\begin{array}{l}
4\varepsilon_1 \ln t+\ln\sqrt{\frac{\mu}{4}} +x-\ln 2|T|+O(\delta_1)+O(a)=0 \\
-2\sqrt\mu\delta_1+\frac{\ln t}{8t}+O\left(\frac{\varepsilon_1 \ln
t}{t}\right) +O\left(b^2\right)+O(\delta_1^2)=0,
\end{array}\right.\label{a2b2_4}
\end{equation}
which we solve for $\varepsilon_1$ and $\delta_1$
\begin{equation}
\left\{\begin{array}{l}
\varepsilon_1 = - \frac{\ln\sqrt{\frac{\mu}{4}} +x-\ln 2|T| }{4\ln t}+O\left(\frac{1}{t}\right)\\
\delta_1=\frac{\ln t}{16\sqrt\mu t}+O\left(\frac{1}{t}\right)
\end{array}\right.\label{a2b2_4}
\end{equation}
and thus
\begin{equation}
\left\{
\begin{array}{l}
a=\frac{\ln t}{8t}+\frac{\ln\frac{2|T|}{\sqrt\mu} +\ln 2-x }{4t}+O\left(\frac{\ln t }{t^2}\right)\\
b=\sqrt{\frac{\mu}{4t}}+\frac{\ln t}{16\sqrt\mu
t^{3/2}}+O\left(\frac{1}{ t^{3/2}}\right).
\end{array}\right.\label{a2b2as_2}
\end{equation}
In a similar manner we compute the next order terms for $a$ and
$b$ by introducing correction terms with unknown functions
$\varepsilon_2$ and $\delta_2$
\begin{equation}
\left\{
\begin{array}{l}
a=\frac{\ln t}{8t}+\frac{A_2}{t}+\frac{\varepsilon_2\ln t}{t^2}\\
b=\sqrt{\frac{\mu}{4t}}+\frac{B_2\ln t}{t^{3/2}}+\frac{\delta_2}{
t^{3/2}},
\end{array}\right.\label{a2b2as_2}
\end{equation}
where $A_2=\frac{1}{4}\ln{\frac{2|T|}{\sqrt\mu}}+\frac{\ln
2-x}{4}$ and $B_2=\frac{1}{16\sqrt{\mu}}$ then system
(\ref{a2b2_1}) becomes
\begin{equation}
\left\{\begin{array}{l} \frac{2B_2\ln t}{\sqrt\mu t}
+\frac{4\varepsilon_2\ln t}{t} -\frac{\ln t}{2\mu t}  +O\left(\frac{1}{t}\right)=0 \\
-\frac{2\sqrt\mu\delta_2}{t}+\frac{A_2}{t}+\frac{1}{2t}+O\left(\frac{\ln
t}{t^2}\right)=0,
\end{array}\right.\label{a2b2_6}
\end{equation}
which we solve for $\varepsilon_2$ and $\delta_2$:
\begin{equation}
\left\{\begin{array}{l} \varepsilon_2= -\frac{B_2}{2\sqrt\mu }
 +\frac{1}{8\mu } +O\left(\frac{1}{\ln t}\right) \\
\delta_2=\frac{A_2+\frac{1}{2}}{2\sqrt\mu}+O\left(\frac{\ln
t}{t}\right).
\end{array}\right.\label{a2b2_7}
\end{equation}
This provides the three leading terms of the asymptotics for $a$
and $b$
\begin{equation}
\left\{
\begin{array}{l}
a=\frac{1}{8}\frac{\ln t}{t}+\frac{\ln \frac{2|T|}{\sqrt\mu}+(\ln
2-x)}{4}\frac{1}{t}
 +\frac{3}{32\mu}\frac{\ln t }{t^{2}}+O\left(\frac{1}{t^{2}}\right) \\
b=\sqrt{\frac{\mu}{4}}\frac{1}{\sqrt{t}}+\frac{1}{16\sqrt{\mu}}\frac{\ln
t}{t^{3/2}}
 +\frac{2+\ln \frac{2|T|}{\sqrt\mu} +(\ln 2-x)}{8\sqrt{\mu}}\frac{1}{t^{3/2}}+O\left(\frac{\ln t
 }{t^{5/2}}\right),
\end{array}
\right.
\end{equation}
which completes the asymptotics of $\alpha_2$.

\subsection{Technical Lemma: asymptotics of an integral}
\begin{lemma} \label{lemma_2}
Let $f(s)$ be twice continuously differentiable on $[0,p]$ with
some $p>0$ and $f'''(s)$ exists and bounded in a small
neighborhood of $s=0$ then
\[
\int_0^pf(s)\left(\sqrt{s^2+b^2}-s\right)ds=-\frac{f(0)}{2}b^2\ln
b
\]
\begin{equation}
+\left(\frac{f(0)}{4}+\frac{f(0)\ln 2}{2}+\frac{f(p)\ln
p}{2}-\frac{1}{2}\int_0^p f'(s)\ln s ds\right)b^2+O(b^3), \quad
b\to 0^+.\label{int_as}
\end{equation}
\end{lemma}

We just outline the proof:
\begin{enumerate}
\item Split the original interval into two subintervals: small neighborhood near zero
and the rest
\begin{equation}
[0,p]=[0,\delta]\cup[\delta,p],
\end{equation}
for some $\delta\to 0^+$.

\item Show
\begin{equation}
\int_0^\delta\left(\sqrt{s^2+b^2}-s\right)ds =-\frac{b^2\ln
b}{2}+\left(\frac{1}{4}+\frac{\ln
2}{2}+\frac{\ln\delta}{2}\right)b^2
+O\left(\frac{b^4}{\delta^2}\right),
\end{equation}
\begin{equation}
\int_0^\delta s\left(\sqrt{s^2+b^2}-s\right)ds =\frac{b^2
\delta}{2}+O\left(b^3\right),
\end{equation}
\begin{equation}
\int_0^\delta s^2\left(\sqrt{s^2+b^2}-s\right)ds
=O\left(b^2\delta^2\right).
\end{equation}
\item Show
\[
\int_0^\delta f(s)\left(\sqrt{s^2+b^2}-s\right)ds
=f(0)\left[-\frac{b^2\ln b}{2}+\left(\frac{1}{4}+\frac{\ln
2}{2}+\frac{\ln\delta}{2}\right)b^2\right]
\]
\begin{equation}
+f'(0)\frac{b^2\delta}{2}
+O\left(b^3\right)+O\left(b^2\delta^2\right).
\end{equation}
\item Show
\begin{equation}
\int_\delta^p f(s)\left(\sqrt{s^2+b^2}-s\right)ds
=b^2\int_\delta^p\frac{f(s)}{2s}ds
+O\left(b^4\int_\delta^p\frac{f(s)}{s^3}ds\right).
\end{equation}
\item Show
\[
\int_\delta^p \frac{f(s)}{2s}ds =\frac{f(p)\ln
p}{2}-\frac{f(0)\ln\delta}{2} -\frac{f'(0)\delta}{2}
\]
\begin{equation}
 -\frac{1}{2}\int_0^pf'(s)\ln s\ ds +O(\delta^2),
\end{equation}
\begin{equation}
\int_\delta^p \frac{f(s)}{s^3}ds=
O\left(\frac{1}{\delta^2}\right).
\end{equation}
\item Show
\begin{equation}
\int_\delta^p f(s)\left(\sqrt{s^2+b^2}-s\right)ds
\end{equation}
\[
=b^2\left(\frac{f(p)\ln p}{2}-\frac{f(0)\ln\delta}{2}
-\frac{f'(0)\delta}{2}  -\frac{1}{2}\int_0^pf'(s)\ln s\ ds \right)
+O(b^2\delta^2)+O\left(\frac{b^4}{\delta^2}\right).
\]

\item Set $\delta=\sqrt{b}$.
\end{enumerate}

\subsection{Asymptotics of key integrals}
Here we derive the asymptotics of the integrals in (\ref{eq1.4:g'_after_split})
\begin{equation}
\left\{
\begin{array}{l}
I_1=-\frac{1}{2}\ln\left(-\mm-z\right)+\frac{\mu}{4z}+\frac{1}{2}\ln(-z)
+O(a^2),\\
I_2= \frac{b^2}{2z^2}I_1+O(ab^2),\\
I_3 = \frac{1}{4z^2}b^2\ln b +\left(\frac{1-2\ln\mu}{8z^2}
+\frac{\ln\left(1+\frac{\mu}{2z}\right)}{4z^2} \right)b^2 +O(ab),\\
I_4= O(b^4\ln b)=O\left(\frac{\ln t}{t^2}\right),\\
I_5=\frac{1}{z^2}\left(2tab^2-\frac{C_1b^2}{4}\right)+O(ab),\quad C_1=\ln \frac{2|T|}{\mu}-x,\\
I_6= O(t ab^4)=O\left(\frac{\ln t}{t^2}\right),\\
H_2 = - \frac{3tb^4}{4 z^3} + O(b^3),\\
H_k = O(tb^{k+2})=O\left(\frac{1}{t^{3/2}}\right), \quad k=3,4,\ldots
\end{array}
\right.
\end{equation}
in the long time limit $t\to\infty$ ($a\to 0$, $b\to 0$).

\subsubsection{Asymptotics of $I_1$, $I_2$}

Let $|z|>|a|$. Consider
\begin{equation}
I_1=-\frac{1}{2}\int_{a}^{-\frac{\mu}{2}} \frac{\xi}{(\xi-z)z}d\xi
\end{equation}
\begin{equation}
=-\frac{1}{2}\int_{a}^{-\frac{\mu}{2}} \left(
\frac{1}{\xi-z}+\frac{1}{z} \right)d\xi =\left.\left(
\frac{1}{2}\ln (\xi-z)-\frac{\xi}{2z} \right)\right|_a^{-\mm}
\end{equation}
and then isolate the leading order $O(1)$
\begin{equation}
= \left[ \frac{1}{2}\ln \left(-\mm-z
\right)+\frac{\mu}{4z}+\frac{1}{2}\ln (-z) \right] +\frac{1}{2}\ln
\left(1-\frac{a}{z} \right)+\frac{a}{2z}
\end{equation}
\begin{equation}
= \frac{1}{2}\ln \left(-\mm-z
\right)+\frac{\mu}{4z}+\frac{1}{2}\ln (-z) +O(a^2).
\end{equation}
The correction term is of the order $O(a^2)=O\left(\frac{\ln^2
t}{t^2}\right)$.

For the second integral
\begin{equation}
I_2= -\frac{1}{2}\left(\frac{R(z)}{\Lambda(z)}+1\right)
\int_{a}^{-\frac{\mu}{2}} \frac{\xi}{(\xi-z)z} d\xi =
-\left(\frac{R(z)}{\Lambda(z)}+1\right) I_1.
\end{equation}
Consider the front factor first
\begin{equation}
\frac{R(z)}{\Lambda(z)}+1 = 1-\sqrt{ 1+ \frac{b^2}{(z-a)^2}}
+exp.small
\end{equation}
for $|z-a|>|b|$
\begin{equation}
= - \frac{b^2}{2(z-a)^2} + O(b^4)
\end{equation}
for $|z|>|a|$
\begin{equation}
=  - \frac{b^2}{2z^2}\left[1+ O(a)\right]^2 + O(b^4)
\end{equation}
\begin{equation}
=  - \frac{b^2}{2z^2}\left[1+ O(a)\right] + O(b^4) =  -
\frac{b^2}{2z^2} + O(ab^2).
\end{equation}
Thus, since $I_1=O(1)$
\begin{equation}
I_2 =  \frac{b^2}{2z^2}I_1 + O(ab^2).
\end{equation}

\subsubsection{Asymptotics of $I_3$, $I_4$}

Consider
\begin{equation}
I_3=-\frac{1}{2}\int_{a}^{-\frac{\mu}{2}}
\left(\frac{\Lambda(\xi)}{R(\xi)} -1
\right)\frac{\xi}{(\xi-z)z}d\xi
\end{equation}
\begin{equation}
=\frac{1}{2}\int_{a}^{-\frac{\mu}{2}} \left(
\frac{\xi-a}{\sqrt{\left(\xi-a\right)^2+b^2}} + 1
\right)\left(\frac{1}{\xi-z}+ \frac{1}{z}\right)d\xi.
\end{equation}
Next, we perform integration by parts
\begin{equation}
=\frac{1}{2}\int_{a}^{-\frac{\mu}{2}} \left(\frac{1}{\xi-z}+
\frac{1}{z}\right)d\left(\sqrt{\left(\xi-a\right)^2+b^2}+(\xi-a)\right)
\end{equation}
\begin{equation}
= \frac{1}{2}\
\frac{\left(-\mm\right)\left[\sqrt{\left(\mm+a\right)^2+b^2}-\left(\mm+a\right)\right]}{\left(-\mm-z\right)z}
-\frac{ab}{2(a-z)z}  + \frac{1}{2}\int_{0}^{-\mm-a}
\frac{\sqrt{s^2+b^2}+s}{(s+a-z)^2}\ d s
\end{equation}
and we keep only terms of the order up to
$O(b^2)=O\left(\frac{1}{t}\right)$, $O(ab)=O\left(\frac{\ln
t}{t^{3/2}}\right)$
\begin{equation}
= \frac{\mu\ b^2}{4 \left(\mm+z\right)z\left[
\sqrt{\left(\mm+a\right)^2+b^2}+\left(\mm+a\right) \right]} +O(ab)
- \frac{1}{2}\int_{0}^{\mm+a} \frac{\sqrt{s^2+b^2}-s}{(s-a+z)^2}\
d s
\end{equation}
\begin{equation}
= \frac{\ b^2}{4 \left(\mm+z\right)z} +O(ab) -
\frac{1}{2}\int_{0}^{\mm+a} \frac{\sqrt{s^2+b^2}-s}{(s-a+z)^2}\ d
s. \label{eq1.4:I_3_split_expr}
\end{equation}
In the integral $a$ enters as a part of the upper limit of
integration and as a part of the integrand. Introduce a notation
\begin{equation}
J(a) = \frac{1}{2}\int_{0}^{\mm+a}
\frac{\sqrt{s^2+b^2}-s}{(s-a+z)^2}\ d s
\end{equation}
then
\begin{equation}
J(a)=J(0)+J'_a(0)\ a +O(a^2)
\end{equation}
\begin{equation}
=\frac{1}{2}\int_{0}^{\mm} \frac{\sqrt{s^2+b^2}-s}{(s+z)^2}\ d s
+a\left[
-\frac{\sqrt{\left(\mm\right)^2+b^2}-\mm}{2\left(\mm+z\right)^2}
-\frac{1}{2}2\int_{0}^{\mm} \frac{\sqrt{s^2+b^2}-s}{(s+z)^3}\ d
s\right] +O(a^2)
\end{equation}
\begin{equation}
=\frac{1}{2}\int_{0}^{\mm} \frac{\sqrt{s^2+b^2}-s}{(s+z)^2}\ d s
+O(ab^2) - a\int_{0}^{\mm} \frac{\sqrt{s^2+b^2}-s}{(s+z)^3}\ d s.
\label{eq1.4:I_3_end_expr}
\end{equation}
Both integrals have a similar structure $\int_{0}^{p}
f(s)\left(\sqrt{s^2+b^2}-s\right) d s$ which was studied in the previous
section. From Corollary \ref{lemma_2}, the first integral has asymptotics
\[
-\int_0^{\mu/2}\frac{1}{2(s+z)^2}\left(\sqrt{s^2+b^2}-s\right)ds
=\frac{1}{4z^2}b^2\ln b
\]
\begin{equation}
+\left(\frac{1-2\ln\mu}{8z^2} -\frac{1}{2\mu z}
+\frac{1}{2\mu\left(z+\mm\right)}
+\frac{\ln\left(1+\frac{\mu}{2z}\right)}{4z^2} \right)b^2 +O(b^3),
\quad b\to 0
\end{equation}
and the second integral in (\ref{eq1.4:I_3_end_expr}) has the order
\begin{equation}
\int_{0}^{\mm} \frac{\sqrt{s^2+b^2}-s}{(s+z)^3}\ d s =O(b^2\ln b).
\end{equation}
Thus, returning to (\ref{eq1.4:I_3_split_expr})
\[
I_3 = \frac{\ b^2}{4 \left(\mm+z\right)z} +O(ab)
+\frac{1}{4z^2}b^2\ln b
\]
\begin{equation}
+\left(\frac{1-2\ln\mu}{8z^2} -\frac{1}{2\mu z}
+\frac{1}{2\mu\left(z+\mm\right)}
+\frac{\ln\left(1+\frac{\mu}{2z}\right)}{4z^2} \right)b^2 +O(b^3)
+ O(ab^2\ln b)
\end{equation}
and since $O(ab)=O\left(\frac{\ln t}{t^{3/2}}\right)$, $O(ab^2\ln
b)=O\left(\frac{\ln t}{t^2}\right)$ and
$O(b^3)=O\left(\frac{1}{t^{3/2}}\right)$
\[
= -\frac{\ b^2}{2\mu \left(\mm+z\right)}+\frac{b^2}{2\mu z}
+\frac{1}{4z^2}b^2\ln b
\]
\begin{equation}
+\left(\frac{1-2\ln\mu}{8z^2} -\frac{1}{2\mu z}
+\frac{1}{2\mu\left(z+\mm\right)}
+\frac{\ln\left(1+\frac{\mu}{2z}\right)}{4z^2} \right)b^2 +O(ab).
\end{equation}
So
\begin{equation}
I_3 = \frac{1}{4z^2}b^2\ln b +\left(\frac{1-2\ln\mu}{8z^2}
+\frac{\ln\left(1+\frac{\mu}{2z}\right)}{4z^2} \right)b^2 +O(ab).
\end{equation}

Consider
\begin{equation}
I_4=\frac{1}{2}\left(\frac{R(z)}{\Lambda(z)}+1\right)
\int_{a}^{-\frac{\mu}{2}} \left( \frac{\Lambda(\xi)}{R(\xi)} -1
\right) \frac{\xi}{\xi-z}d\xi=
-\left(\frac{R(z)}{\Lambda(z)}+1\right) I_3.
\end{equation}
Similar to our calculations in Step 1
\begin{equation}
= -\left(- \frac{b^2}{2z^2} + O(ab^2) \right)I_3
\end{equation}
and since $I_3=O(b^2\ln b)$
\begin{equation}
I_4 = O(b^4 \ln b).
\end{equation}

\subsubsection{Asymptotics of $I_5$, $I_6$}

Consider
\begin{equation}
I_5= \frac{1}{2 \pi i z} \int_{a+bi}^{a-bi}
\frac{\Lambda(\xi)}{R(\xi)} f'(\xi)\ \frac{\xi}{z}d\xi
\end{equation}
\begin{equation}
= \frac{1}{2 \pi i z^2} \int_{a+bi}^{a-bi} \frac{-(\xi-a) \xi
f'(\xi)}{\sqrt{\left(\xi-a\right)^2+b^2}} d\xi
\end{equation}
changing variables $\xi=by+a$
\begin{equation}
= - \frac{1}{2 \pi i z^2} \int_{i}^{-i} \frac{ by (by+a)
f'(by+a)}{\sqrt{y^2+1}} dy
\end{equation}
decomposing $f(\xi)$ into powers of $\xi$ as in (\ref{df_series})
\[
= - \frac{1}{2 \pi i z^2} \int_{i}^{-i} \frac{ by (by+a)
}{\sqrt{y^2+1}} \left[
-4t(by+a)+C_1+iC_2\mbox{sign}\left(\frac{y}{i}\right) \right.
\]
\begin{equation}
\left.+C_3(by+a)+C_4(by+a)^2+O\left((by+a)^3\right) \right]dy,
\end{equation}
where (see (\ref{C_values}))
\begin{equation}
\left\{\begin{array}{l}
C_1=\ln \frac{2|T|}{\mu}-x\\
C_2=\frac{\pi}{2}\\
C_3=\frac{2}{\mu},
\end{array}\right.
\end{equation}
leading to
\begin{equation}
= - \frac{1}{2 \pi i z^2} \int_{i}^{-i} \frac{ (b^2y^2+aby)
}{\sqrt{y^2+1}} \left[
-4tby-4ta+C_1+iC_2\mbox{sign}\left(\frac{y}{i}\right) +O(b)
\right]dy
\end{equation}
by the symmetry argument all odd powers of $y$ vanish
\begin{equation}
= - \frac{1}{2 \pi i z^2} \int_{i}^{-i} \frac{1 }{\sqrt{y^2+1}}
\left[ -8tab^2y^2 +C_1b^2y^2 +O(ab) \right] dy
\end{equation}
\begin{equation}
= - \frac{\left( -8tab^2 +C_1b^2  \right)}{2 \pi i z^2}
\int_{i}^{-i} \frac{y^2}{\sqrt{y^2+1}} dy +O(ab).
\end{equation}
Thus since $\int_{i}^{-i} \frac{y^2}{\sqrt{y^2+1}} dy=\frac{\pi
i}{2}$
\begin{equation}
I_5= \frac{1}{z^2}\left( 2tab^2 +\frac{C_1b^2}{4}  \right)
 +O(ab).
\end{equation}

Consider
\begin{equation}
I_6= -\frac{1}{2 \pi i z}
\left(\frac{R(z)}{\Lambda(z)}+1\right)\int_{a+bi}^{a-bi}
\frac{\Lambda(\xi)}{R(\xi)} f'(\xi)\ \frac{\xi}{z}d\xi
\end{equation}
\begin{equation}
=-\left(\frac{R(z)}{\Lambda(z)}+1\right) I_5 =O(b^2) O(tab^2)
=O\left(\frac{\ln t}{t^2}\right).
\end{equation}

\subsubsection{Asymptotics of $H_k$}

Consider
\begin{equation}
H_2= - \frac{1}{2 \pi i z}\ \frac{R(z)}{\Lambda(z)}
\int_{a+bi}^{a-bi} \frac{\Lambda(\xi)}{R(\xi)} f'(\xi)
\left(\frac{\xi}{z}\right)^2 d\xi
\end{equation}
\begin{equation}
= - \frac{1}{2 \pi i z^3}\ (-1+O(b^2))\int_{a+bi}^{a-bi}
\frac{-(\xi-a) f'(\xi)}{\sqrt{\left(\xi-a\right)^2+b^2}} f'(\xi)
\xi^2 d\xi
\end{equation}
with the change of variables $\xi=by+a$
\begin{equation}
=  -\frac{1+O(b^2)}{2 \pi i z^3} \int_{i}^{-i} \frac{ by (by+a)^2
f'(by+a)}{\sqrt{y^2+1}} dy
\end{equation}
similar to computation of $I_5$
\begin{equation}
= -\frac{1+O(b^2)}{2 \pi i z^3} \int_{i}^{-i} \frac{
(b^3y^3+2ab^2y^2+a^2by) }{\sqrt{y^2+1}} \left[ -4tby-4ta +O(1)
\right]dy
\end{equation}
\begin{equation}
=  -\frac{1+O(b^2)}{2 \pi i z^3} \int_{i}^{-i}
\frac{1}{\sqrt{y^2+1}} \left[ -4tb^4y^4+O(tab^3)-4tab^3y^3
+O(ta^2b^2) + O(b^3)\right]dy.
\end{equation}
By the symmetry of the integral the $y^3$ term vanishes as well as
the term $O(tab^3)$, leaving
\begin{equation}
H_2= (1+O(b^2)) \frac{2tb^4}{\pi i z^3} \int_{i}^{-i}
\frac{y^4}{\sqrt{y^2+1}} dy +O(b^3)
\end{equation}
\begin{equation}
=  \frac{2tb^4}{\pi i z^3} \int_{i}^{-i} \frac{y^4}{\sqrt{y^2+1}}
dy +O(b^3) +O(tb^6)
\end{equation}
with a table integral $\int_{i}^{-i} \frac{y^4}{\sqrt{y^2+1}}
dy=-\frac{3\pi i}{8}$.

So
\begin{equation}
H_2 = -\frac{3 tb^4}{4 z^3} +O(b^3).
\end{equation}

Consider for $k=3,4, \ldots$
\begin{equation}
H_k= - \frac{1}{2 \pi i z}\ \frac{R(z)}{\Lambda(z)}
\int_{a+bi}^{a-bi} \frac{\Lambda(\xi)}{R(\xi)} f'(\xi)
\left(\frac{\xi}{z}\right)^k d\xi
\end{equation}
similarly to $H_2$
\begin{equation}
= - \frac{1}{2 \pi i z^{k+1}}\ (-1+O(b^2))\int_{a+bi}^{a-bi}
\frac{-(\xi-a)}{\sqrt{\left(\xi-a\right)^2+b^2}} f'(\xi) \xi^{k}
d\xi
\end{equation}
with the change of variables $\xi=by+a$
\begin{equation}
=  -\frac{1+O(b^2)}{2 \pi i z^{k+1}} \int_{i}^{-i} \frac{ by
(by+a)^k f'(by+a)}{\sqrt{y^2+1}} dy
\end{equation}
\begin{equation}
=  -\frac{1+O(b^2)}{2 \pi i z^{k+1}} \int_{i}^{-i} \frac{
O(b)O(b^k) O(tb)}{\sqrt{y^2+1}} dy.
\end{equation}

Thus
\begin{equation}
H_k=O(tb^{k+2})=O(b^k)=O(b^3).
\end{equation}

\subsection{Numerical evaluations}

\subsubsection{Numerical evaluation of $h'(z)$ and $h(z)$}

The main idea we utilize is implementing integration of
functions on Riemann surfaces
rather than in the complex plane. Such approach together with
numerical contour deformations has allowed to avoid expensive
computations of the main arcs as a preparation for any single
computations involving $h$-function. We were able to continuously
track the branch points $\alpha_2$ in genus 2 beyond colliding
with the branch cut $[0,T]$ to another sheet of a Riemann surface.
This leads to easier long time computations and allows to observe the
singular obstruction.

To compute $h'(z)$ we assume that all main and complementary arcs
are in series configuration, and the positions of the branch points
$\alpha$'s are known (see below), then
\begin{equation}
h'(z)=\frac{R(z)}{2\pi
i}\oint_{\ggh}\frac{f'(\xi)}{(\xi-z)R(\xi)}d\xi
,\label{eq3.3:general h'-function}
\end{equation}
where the point $z$ is inside of the large loop $\ggh$ (Fig.
\ref{fig3:h_cont_g2LT}) and where
\begin{equation}
f'(\xi)=-\frac{\pi i}{2}-\ln\left(\frac{\mu}{2}
-\xi\right)+\frac{1}{2}\ln\left(\xi^2-T^2\right) -x-4t\xi, \quad
\mbox{when}\ \ \Im z > 0. \label{eq3.3:general f'-function}
\end{equation}

Under the same assumptions,
we evaluate the function $h(z)$ directly rather than
integrating the derivative $h'(z)$.

\[
h(z)=\frac{R(z)}{2\pi
i}\left[\oint_{\ggh}\frac{f(\xi)}{(\xi-z)R(\xi)}d\xi +
\oint_{\ggh_c}\frac{\Omega}{(\xi-z)R(\xi)}d\xi+
\oint_{\ggh_m}\frac{W}{(\xi-z)R(\xi)}d\xi \right],
\]
where point $z$ lies inside of a large loop $\ggh$ and outside
of small loops $\ggh_c$ and $\ggh_m$ (Fig. \ref{fig3:h_cont_g2}).

The constants $\Omega$ and $W$ are solutions of the linear system

\begin{equation}
\left(
\begin{array}{cc}
\oint_{\ggh_m}\frac{1}{R(\xi)}d\xi &
\oint_{\ggh_c}\frac{1}{R(\xi)}d\xi\\
\oint_{\ggh_m}\frac{\xi}{R(\xi)}d\xi&
\oint_{\ggh_c}\frac{\xi}{R(\xi)}d\xi
\end{array}
\right)\left(
\begin{array}{c}
W\\
\Omega
\end{array}
\right)=\left(
\begin{array}{c}
-\oint_{\ggh_m}\frac{f(\xi)}{R(\xi)}d\xi\\
-\oint_{\ggh_m}\frac{\xi f(\xi)}{R(\xi)}d\xi
\end{array}
\right)
\end{equation}

The code to evaluate $h(z)$ could be viewed as a tool to support
and track evolution of RH contours. It provides valuable insight
into the behavior of the branch points and the zero level curves
of $\Im h$. This tool allows one to identify which one of several possible
scenarios of level curve evolution does occur.

\begin{remark}
The only principal difference between evaluating $h(z)$ and $g(z)$ is a requirement
for the point $z$ to be located either inside contour $\ggh$ (for $h(z)$), or outside
$\ggh$ (for $g(z)$).

Keeping in mind simple relation between $h=2g-f$, it is easy to switch between these
two functions. For example, for distant $z$ it is more efficient to use $g$.
\end{remark}

\subsubsection{Numerical evaluation of $\alpha$'s}

We compute the branch points $\left(\alpha_0,\alpha_2,\alpha_4\right)$ in genus 2
by solving the system
\begin{equation}
\left\{
\begin{array}{l}
B(\alpha_0)=0 \\
B(\alpha_2)=0 \\
B(\alpha_{4})=0
\end{array}
\right. ,\label{num_alphas_B}
\end{equation}
where
\begin{equation}
B(z)=\oint_{\ggh}\frac{f(\xi)}{(\xi-z)R(\xi)}d\xi +
\oint_{\ggh_c}\frac{\Omega}{(\xi-z)R(\xi)}d\xi+
\oint_{\ggh_m}\frac{W}{(\xi-z)R(\xi)}d\xi.
\label{eq1.3:general_B_function}
\end{equation}

Up to a constant, $B(z)$ is $h(z)$ without the $R$ factor in
front of the integrals.

\begin{remark}
We stress that $B(z)$ is in fact a function of $\alpha$'s and
$\overline{\alpha}$'s through the $R$ factor in the integrals,
that is $B(z)=B(z,\vec\alpha)$. While $B(z,\alpha)$ is an analytic
function of $z$, $B(\alpha,\alpha)$ is a non-analytic function of
$\alpha$ which depends both on $\alpha$ and $\overline{\alpha}$
through $R(\xi)$ in the denominators. So we treat the system
(\ref{num_alphas_B}) as a real $6$x$6$ system and solve it
iteratively.
\end{remark}

\subsubsection{Numerical computations of the first breaking curve}

Computations of the first break for $\mu=0$ was done by Lyng and
Miller \cite{Lyng_07}. We computed the first break for $\mu>0$. The
first break was also observed as a singular event in $x$, $t$ and
even $\mu$ evolution of the zero level curves of $\Im h$.

We treat the first breaking curve $t=t_0(x)$ as a function of $x$. For
fixed $\mu$ and $x$ we are looking for a pair $(z_0,t)$ which
satisfies the system of one complex and one real equations
\begin{equation}
\left\{ \begin{array}{l} h'(z_0,t)=0\\
\Im h(z_0,t)=0,
\end{array}\right. \label{eq:br1_syst_num}
\end{equation}
where $\Im z_0 \ge0$ and $t\ge0$. The Jacobian of this system is
singular. We start with expending $h(z)$ in powers of $z-z_0$
\begin{equation}
h(z)=A+B(z-z_0)+\frac{C}{2}(z-z_0)^2+\ldots,
\end{equation}
where $A=h(z_0)$, $B=h'(z_0)$ and $C=h''(z_0)$. Then the system
is approximated as
\begin{equation}
\left\{ \begin{array}{l}
B+C(z-z_0)=0\\
\Im \left[A+B(z-z_0)+\frac{C}{2}(z-z_0)^2\right]=0,
\end{array}\right. .
\end{equation}
Solving the first equation and substituting into the second
equation leads to
\begin{equation}
\Im \left[A-\frac{B^2}{2C}\right]=0.
\end{equation}
Thus in terms of the function $h$ the system (\ref{eq:br1_syst_num}) is replaced
with
\begin{equation}
\left\{ \begin{array}{l}
h'(z_0,t_0)=0\\
\Im \left[h(z_0,t_0)-\frac{(h'(z_0,t_0))^2}{2h''(z_0,t_0)}\right]=0
\end{array}\right. . \label{eq:num_br1}
\end{equation}
This system is solved iteratively where $z_0$ and $t_0$ are updated in turns
\begin{equation}
(z_0^{(0)},t_0^{(0)})\to z_0^{(1)}\to t_0^{(1)} \to
z_0^{(2)}\to\ldots.
\end{equation}
We use the first equation in (\ref{eq:num_br1}) to update $z_0$ and the
second equation to update $t_0$.

\subsubsection{Numerical computations of the singular obstruction}

From the numerical point of view, the singular obstruction curve
is a solution of a scalar equation
\begin{equation}
\Im h(T,x,t)=0
\end{equation}
for either $t=t_c(x)$ or $x=x_c(t)$.


\begin{thebibliography}{10}


\bibitem{BelovVen_pre}
Belov, S., Venakides, S., \emph{Smooth parametric dependence of
asymptotics of the semiclassical focusing NLS}, preprint.


\bibitem{Belok}
Belokolos, E.D., Bobenko, A.I., Enol'skii, V.Z., Its, A.R.,
Matveev, V.B., {\it Algebraic-geometric approach to nonlinear
integrable equations}, Springer-Verlag, New York 1994.


\bibitem{Bronski1}
Bronski, J.C., {\it Semiclassical eigenvalue distribution of the
Zakharov-Shabat eigenvalue problem}. Phys. D 97 (1996), no. 4,
376--397.


\bibitem{Bronski2}
Bronski, J.C., {\it  Spectral instability of the semiclassical
Zakharov-Shabat eigenvalue problem}, Advances in nonlinear
mathematics and science, Phys. D 152/153 (2001), 163--170.

\bibitem{BMM}
Bronski, J.C., McLaughlin, K.T.-R., Miller, P.D., {\it Rigorous
asymptotics for the point spectrum of the nonselfadjoint
Zakharov-Shabat eigenvalue problem with Klaus-Shaw potential}, in
preparation.

\bibitem{BTVZ_07}
Buckingham, R., Tovbis, A., Venakides, S., Zhou, X., {\it The
semiclassical focusing nonlinear Schrodinger equation. Recent
advances in nonlinear partial differential equations and
applications}, 47--80, Proc. Sympos. Appl. Math., 65, Amer. Math.
Soc., Providence, RI, 2007.

\bibitem{BV_07}
Buckingham, R., Venakides, S., {\it Long-time asymptotics of the
nonlinear Schrodinger equation shock problem}, Comm. Pure Appl.
Math. 60 (2007), no. 9, 1349--1414.


\bibitem{Ca}
Cai, D., McLaughlin, D.W., McLaughlin, K.T.R., {\it The nonlinear
Schroedinger equation as both a PDE and a dynamical system},
Handbook of dynamical systems, Vol. 2, 599-675, North-Holland,
Amsterdam, 2002.


\bibitem{CT}
Ceniceros, H., Tian, F.-R., {\it A numerical Study of the
semi-classical limit of the focusing nonlinear Schr\" odinger
equation}, Phys. Lett. A 306 (2002), no. 1, 25--34.

\bibitem{Deift_OPbook_99}
Deift, P., {\it Orthogonal polynomials and random matrices:
A Riemann-Hilbert approach}, In Courant Lecture
Notes in Mathematics, Volume 3, CIMS, New York, 1999.

\bibitem{DVZ1_97}
Deift, P., Venakides, S., Zhou, X., {\it
New results in small dispersion KdV by an extension
of the steepest descent method for Riemann-Hilbert problems},
Internat. Math. Res. Notices 6 (1997), 286-299.

\bibitem{DVZ2}
Deift, P., Venakides, S., Zhou, X., {\it
An extension of the steepest descent method for Riemann-Hilbert
problems: The small dispersion limit of the Korteweg-de Vries
(KdV) equation}, Proc. Natl. Acad. Sci., Vol. 95, (1998) pp. 450–454.

\bibitem{DZ1}
Deift, P., Zhou, X., {\it A steepest descent method for
oscillatory Riemann - Hilbert problems.  Asymptotics for the mKdV
equation}, Ann. of Math. 137 (1993), 295-370.

\bibitem{DZ2}
Deift, P., Zhou, X., {\it Asymptotics for the Painlev\'e II
equation}, Comm. Pure and Appl. Math. 48 (1995), 277-337.

\bibitem{DZ3}
Deift, P., Zhou, X., {\it Long-time asymptotics for solutions of
the NLS equation with initial data in a weighted Sobolev space},
Comm. Pure Appl. Math. 56 (2003), no. 8, 1029--1077.

\bibitem{DM}
DiFranco, J.C., Miller, P.D., {\it The semiclassical modified
nonlinear Schrodinger equation. I. Modulation theory and spectral
analysis}, Phys. D 237 (2008), no. 7, 947--997.

\bibitem{FL}
Forest, M.G., Lee, J.E., {\it Geometry and modulation theory
for the periodic nonlinear Schrodinger equation, in: C. Dafermos,
et. al. (Eds.), Oscillation Theory, Computation, and Methods of
Compensated Compactness}, Vol. 2, IMA, Springer, New York, 1986.


\bibitem{KMM_03}
Kamvissis, S.; McLaughlin, K. D. T.-R.; Miller, P. D., {\it
Semiclassical soliton ensembles for the focusing nonlinear
Schrödinger equation}, Annals of Mathematics Studies, 154.
Princeton University Press, Princeton, N.J., 2003.

\bibitem{KS}
Klaus, M., Shaw, J.K., {\it Purely imaginary eigenvalues of
Zakharov-Shabat systems}, Phys. Rev. E(3). (2002), 65(3):036607,
5pp.

\bibitem{Lax_68}
Lax, P.D., {\it Integrals of nonlinear equations of evolution and
solitary waves}, Comm. Pure Appl. Math. 21 (1968) 467--490.


\bibitem{Lyng_07}
Lyng, G., Miller, P.D., {\it The $N$-Soliton of the Focusing
Nonlinear Schr\"odinger Equation for $N$ Large},  Comm. Pure Appl.
Math. 60 (2007), no. 7, 951--1026.


\bibitem{MillerKamvissis_98}
Miller P.D., Kamvissis, S., {\it On the semiclassical limit of the
focusing nonlinear Schrödinger equation}, Phys. Lett. A 247
(1998), no. 1-2, 75--86.


\bibitem{Newell}
Newell, A.C., {\it Solitons in mathematics and physics}, CBMS-NSF
Regional Conference Series in Applied Mathematics, 48, (SIAM),
Philadelphia, PA, 1985.


\bibitem{Satsuma_74}
Satsuma, J., Yajima, N., {\it Initial value problems of
one-dimensional self-modulation of nonlinear waves in dispersive
media}, Progr. Theoret. Phys. Suppl. No. 55 (1974), 284–306.


\bibitem{Shabat_76}
Shabat, A., {\it One-dimensional perturbations of a differential
operator, and the inverse scattering problem}, Problems in
mechanics and mathmatical physics, pp. 279–296, Nauka, Moscow,
1976.


\bibitem{TV_00}
Tovbis, A., Venakides, S., {\it The eigenvalue problem for the
focusing nonlinear Schr\"odinger equation: new solvable cases},
Phys. D 146 (2000), no. 1-4, 150–164.



\bibitem{TVparam_09} Tovbis, A., Venakides, S., {\it Determinant Form of the
Complex Phase Function of the Steepest Descent Analysis of Riemann–Hilbert Problems and
Its Application to the Focusing Nonlinear Schr\"{o}dinger
Equation}, IMRN no 11 (2009), 2056--2080.


\bibitem{TVZzero_04}
Tovbis, A., Venakides, S., Zhou, X., {\it On semiclassical (zero
dispersion limit) solutions of the focusing Nonlinear
Schr\"odinger equation}, Comm. Pure Appl. Math. 57 (2004), 0877-0985.


\bibitem{TVZlong2_06}
Tovbis, A., Venakides, S., Zhou, X., {\it On the long-time limit
of semiclassical (zero dispersion limit) solutions of the focusing
nonlinear Schrodinger equation: pure radiation case}, Comm. Pure
Appl. Math. 59 (2006), no. 10, 1379--1432.


\bibitem{TVZ_07}
Tovbis, A., Venakides, S., Zhou, X., {\it Semiclassical focusing
nonlinear Schrodinger equation I: inverse scattering map and its
evolution for radiative initial data}, Int. Math. Res. Not.
(2007), no. 22, Art. ID rnm094.


\bibitem {Whitham_74}
Whitham, G.B., {\it Linear and Nolinear Waves}, Wiley, 1974.


\bibitem{ZSh_72}
Zakharov, V.E., Shabat, A.B., {\it Exact theory of two dimensional
selffocusing and onedimensional selfmodulation of waves in
nonlinear media}, Soviet Physics JETP, Vol. 34, pp. 62-69, (1972).


\bibitem{Z_98}
Zhou, X., {\it $L^2$-Sobolev space bijectivity of the scattering
and inverse scattering transforms}, Comm. Pure Appl. Math. 51
(1998), no. 7, 697--731.


\bibitem{ZZ}
Zhou, X., {\it Riemann-Hilbert problems and integrable systems},
preprint.


\end{thebibliography}
\end{document}